\documentclass[10pt,a4paper]{article}
\usepackage[utf8]{inputenc}
\usepackage{amsmath}
\usepackage{amsfonts}
\usepackage{amssymb}
\usepackage{amsmath}
\usepackage{multicol,color}
\usepackage{float}
\usepackage{graphicx}
\usepackage{amssymb}
\usepackage{theorem}
\usepackage{fancyhdr}
\usepackage{tikz}
\usepackage{enumerate}
\usepackage[margin=1.1in]{geometry}
\usepackage{pgfplots}
\usepgflibrary{shapes.geometric}
\pgfplotsset{my style/.append style={axis x line=middle, axis y line= middle, xlabel={$x$}, ylabel={$y$}, axis equal }}
\usetikzlibrary{through,calc}
\usetikzlibrary{calc,intersections}

\newcommand*{\TlYellow}[1][1]{%
\begin{tikzpicture}[every path/.style={thick,fill=lightgray},scale=#1]
\draw (0,0)    circle (0.07);
\draw (60:0.1) circle (0.07);
\draw[fill=yellow] (0.1,0) circle (0.07);}

\def\disp{\displaystyle}
\def\ve{\varepsilon}
\def\dd{\delta}

\def\lm{\lambda}
\def\O{\Omega}

\def\tilde{\widetilde}

\def\oa{\bar a}
\def\ob{\bar b}

\def\ox{\bar{x}}
\def\oy{\bar{y}}
\def\oz{\bar{z}}

\def\ou{\bar{u}}

\def\d#1{\dot{#1}}

\def\gph{\hbox{}}

\def\gg{\gamma}

\def\tto{\rightrightarrows}

\def\hat{\widehat}
\def\Hat{\widehat}
\def\tilde{\widetilde}

\def\Bar{\overline}
\def\ra{\rangle}
\def\la{\langle}
\def\ve{\varepsilon}
\def\B{I\!\!B}
\def\h{\hfill\Box}
\def\R{\mathbb{R}}
\def\N{\mathbb{N}}

\def\gph{\mbox{\rm gph}\,}

\def\dom{\mbox{\rm dom}\,}

\def\dist{\mbox{\rm dist}}

\def\bd{\mbox{\rm bd}\,}

\def\O{\Omega}

\def\ph{\varphi}
\def\emp{\emptyset}

\def\oR{\Bar{\R}}
\def\lm{\lambda}

\def\gg{\gamma}
\def\tt{\tilde t}
\def\dd{\delta}

\def\al{\alpha}

\def\be{\beta}

\def\N{I\!\!N}

\newtheorem{theorem}{Theorem}[section]

\newtheorem{proposition}[theorem]{Proposition}
\newtheorem{definition}[theorem]{Definition}
\theoremstyle{plain}{\theorembodyfont{\rmfamily}
}
\theoremstyle{plain}{\theorembodyfont{\rmfamily}
}
\theoremstyle{plain}{\theorembodyfont{\rmfamily}
}
\theoremstyle{plain}{\theorembodyfont{\rmfamily}
\newtheorem{example}[theorem]{Example}}

\theoremstyle{plain}{\theorembodyfont{\rmfamily}
}

\renewcommand{\theequation}{\thesection.\arabic{equation}}

\begin{document}
\begin{center}
{\bf OPTIMAL CONTROL OF NONCONVEX INTEGRO-DIFFERENTIAL\\ SWEEPING PROCESSES}\\[3ex]
ABDERRAHIM BOUACH\footnote{Laboratoire LMPEA, Facult\'e des Sciences Exactes et Informatique, Universit\'e Mohammed Seddik Benyahia,
Jijel, B.P. 98, Jijel 18000, Alg\'erie (abderrahimbouach@gmail.com).} TAHAR HADDAD\footnote{Laboratoire LMPEA, Facult\'e des Sciences Exactes et Informatique, Universit\'e Mohammed Seddik Benyahia, Jijel, B.P. 98, Jijel 18000, Alg\'erie (haddadtr2000@yahoo.fr).} and BORIS S. MORDUKHOVICH\footnote{Department of Mathematics, Wayne State University, Detroit, Michigan 48202, USA (aa1086@wayne.edu). Research of this author was partly supported by the US National Science Foundation under grants DMS-1512846 and DMS-1808978, by the US Air Force Office of Scientific Research under grant \#15RT0462, and by the Australian Research Council Discovery Project DP-190100555.}
\end{center}
\small{\bf Abstract.} This paper is devoted to the study, for the first time in the literature, of optimal control problems for sweeping processes governed by integro-differential inclusions of the Volterra type with different classes of control functions acting in nonconvex moving sets, external dynamic perturbations, and integral parts of the sweeping dynamics. We establish the existence of optimal solutions and then obtain necessary optimality conditions for a broad class of local minimizers in such problems. Our approach to deriving necessary optimality conditions is based on the method of discrete approximations married to basic constructions and calculus rules of first-order and second-order variational analysis and generalized differentiation. The obtained necessary optimality conditions are expressed entirely in terms of the problem data and are illustrated by nontrivial examples that include applications to optimal control models of non-regular electrical circuits.\\[1ex]
\vspace*{-0.2in}

\section{Introduction and Problem Formulation}\label{intro}
\setcounter{equation}{0}\vspace*{-0.1in}

This paper continues a series of recent publications devoted to the rather new optimal control theory for discontinuous differential inclusions governed by {\em controlled sweeping processes}. Recall that the classical (uncontrolled) sweeping process was introduced by Moreau in the 1970s (see his lecture notes \cite{mor_frict} with further references) in the form of the dissipative differential inclusion
\begin{eqnarray}\label{sp}
\dot x(t)\in-N_{C(t)}(x(t))\;\mbox{ a.e. }\;t\in[0,T],
\end{eqnarray}
where $C(t)$ is a family of convex moving sets that continuously depend on time, and where $N_C(x)$ stands for the normal cone of convex analysis. Then the theory of sweeping processes has been largely developed with a variety of applications. In the very recent survey \cite{bt}, the reader can find the abundant bibliography and discussions on these and related topics.\vspace*{-0.05in}

One of the most remarkable features of Moreau's sweeping process and its extensions is that the Cauchy problem for them has a {\em unique} solution. This excludes any optimization of
the discontinuous sweeping dynamics of type \eqref{sp}. A new view on sweeping processes was offered in \cite{chhm}, were the authors suggested to parameterize the moving sets in \eqref{sp} by control functions $C(t)=C(u(t))$ that allowed them to formulate an {\em optimal control} problem and derive first necessary optimality conditions in sweeping control theory. Since that, optimal control theory for various types of controlled sweeping processes governed by ordinary differential inclusions of type \eqref{sp} has been developed in many publications with deriving necessary optimality conditions and applications; see, e.g., \cite{ao,ac,bk,ccmn,ca3,b1,chhm3,cmn,pfs,hm,m19,mn,zeidan} and the references therein.\vspace*{-0.05in}

In contrast to the previous publications, in this paper we consider controlled sweeping processes with the dynamics governed by {\em integro-differential inclusions} of the Volterra type:
\begin{equation}\label{e:1}
\left\{
\begin{array}{l}
-\dot{x}(t) \in N_{C(u(t))}(x(t))+f_{1}(a(t),x(t))+\displaystyle\int\limits_{0}^{t}f_{2}(b(s),x(s))\mathrm{d}s\quad \text{a.e. }\;[0,T],\vspace*{0.3cm}\\
(u(\cdot),a(\cdot),b(\cdot))\in W^{1,2}([0,T],\mathbb{R}^{n+m+d}),\vspace*{0.3cm}\\
x(0)=x_0\in C(0),
\end{array}\right.
\end{equation}
on the fixed time interval $[0,T]$, where the triple $(u(\cdot),a(\cdot),b(\cdot))$ signifies {\em feasible controls} acting in the moving sets, additive perturbations, and the integral part of the sweeping dynamics, respectively. The controlled moving sets are given in the form
\begin{equation}\label{e:mset}
C(t):=C(u(t))=C+u(t),\,\,C:=\big\{x\in\mathbb{R}^{n}\;\big|\;g_{i}(x)\geq 0,\,\,\,\,i=1,\ldots,s\big\}.
\end{equation}
Since the sets $C(t)$ are generally {\em nonconvex}, the normal cone in \eqref{e:1} is understood in the generalized sense defined in Section~\ref{sec:prel}, which reduces to the one in \eqref{sp} for the case of convexity.\vspace*{-0.05in}

Given further a terminal cost function $\varphi:\mathbb{R}^{n}\to\oR:=(-\infty,\infty]$ and a running cost $l_{0}[0,T]\times\mathbb{R}^{4n+2m+2d}\to\mathbb{R}$, the {\em sweeping optimal control problem} $(P)$ consists of minimizing the Bolza-type functional
\begin{equation}\label{J0}
J_{0}[x,u,a,b]:=\varphi(x(T))+\int\limits_{0}^{T}l_{0}(t,x(t),u(t),a(t),b(t),\dot{x}(t),\dot{u}(t),\dot{a}(t),\dot{b}(t))\,dt
\end{equation}
on the set of feasible control  $(u(\cdot),a(\cdot),b(\cdot))$ and the corresponding trajectories $x(\cdot)$ of \eqref{e:1} from the space $W^{1,2}([0.T],\R^n)$. Such quadruples $(x(\cdot),u(\cdot),a(\cdot),b(\cdot))$ are called {\em feasible solutions} to $(P)$. The existence results for feasible and optimal solutions to $(P)$ are given in Section~\ref{sec:exist}. The required assumptions on the functions $\ph,l_0,g_i$ and the mappings $f_{1} : \mathbb{R}^{m+n} \to \mathbb{R}^{n} $ and $ f_{2} : \mathbb{R}^{d+n} \to \mathbb{R}^{n}$ will be formulated in Section~\ref{sec:prel}.

Observe that the integro-differential inclusion \eqref{e:1} can be represented as
\begin{equation}\label{e:2}
\left\{
\begin{array}{l}
-\dot{x}(t) \in N_{C(u(t))}(x(t))+f_{1}(a(t),x(t))+y(t)\quad\mbox{a.e.}\;\;t\in [0,T],\vspace*{0.3cm}\\
\dot{y}(t)=f_{2}(b(t),x(t))\quad \mbox{a.e.}\,\, \;t\in [0,T],\vspace*{0.3cm}\\
(u(\cdot),a(\cdot),b(\cdot))\in W^{1,2}([0,T],\mathbb{R}^{n+m+d}),\vspace*{0.3cm}\\
x(0)=x_0\in C(0),\;y(0)=0,
\end{array}\right.
\end{equation}
where $y(\cdot)\in W^{1,2}([0,T],\R^n)$ is uniquely determined under the assumptions below. Denoting now the corresponding running cost $l : [0,T]\times\mathbb{R}^{6n+2m+2d}\to\mathbb{R}$ by
\begin{equation}\label{e:running cost}
l(t,x,y,u,a,b,\dot{x},\dot{y},\dot{u},\dot{a},\dot{b}):=l_{0}(t,x,u,a,b,\dot{x},\dot{u},\dot{a},\dot{b}),
\end{equation}
we see that problem $(P)$ can be equivalently reformulated in the form: minimize
\begin{equation}\label{e:bolza}
J[x,y,u,a,b]:=\varphi(x(T))+\int\limits_{0}^{T}l(t,x(t),y(t),u(t),a(t),b(t),\dot{x}(t),\dot{y}(t),\dot{u}(t),\dot{a}(t),\dot{b}(t))\,\mathrm{d}t
\end{equation}
subject to feasible solutions of \eqref{e:2}. As follows from the sweeping inclusion \eqref{e:1} and the structure of the controlled moving sets \eqref{e:mset}, problem $(P)$ automatically involves the pointwise {\em mixed state-control constraints}
\begin{equation}\label{mixed}
g_i(x(t)-u(t))\ge 0\;\mbox{ for all }\;t\in[0,T]\;\mbox{ and }\;i=1,\ldots,s,
\end{equation}
which have been recognized among the most challenging issues even in standard optimal control theory for systems governed by smooth controlled ordinary differential equations.\vspace*{-0.05in}

The {\em uncontrolled} integro-differential sweeping process \eqref{e:1} was first introduced in \cite{bre} for the case where $f_{1}\equiv 0$, $f_{2}\equiv f_{2}(x)$, and where $C(t):=\mathbb{E}^{n}$ is the nonnegative orthant of $\R^{n}$. The motivation in \cite{bre} came from the study a one-dimensional flow of particles subject to a force field generated by the fluid itself.
Then (again uncontrolled) integro-differential sweeping process \eqref{e:1} with $f_{1}\equiv 0$ and $f_{2}\equiv f_{2}(x)\in\gamma \mathbb{B}(0,1)$ was considered in \cite{CK}, where the authors established the existence and uniqueness of solutions by developing a new penalization approach. Quite recently, the existence and uniqueness issues for a generalized uncontrolled version of \eqref{e:1} have been revisited in \cite{bou} with the extension of the results in \cite{CK} to a more general framework of \eqref{e:1} by using a new Gronwall-like differential inequality within a developed semi-discretization method. Furthermore, paper \cite{bou} contains results on the continuous dependence of solutions of \eqref{e:1} on the initial data with some applications to sweeping dynamical models arising in electronics.\vspace*{-0.05in}

To the best of our knowledge, optimal control problems for integro-differential sweeping processes have been {\em never formulated} and studied in the literature even in the case of convex and uncontrolled moving sets. The closest prototype of $(P)$ in the case of sweeping differential inclusions with no integral part was considered in \cite{b1}, where the authors established the existence of optimal solutions and obtained necessary optimality conditions for local minimizers. Furthermore, in \cite{cm} they applied the necessary optimality conditions from \cite{b1} to optimal control of the planar version of the crowd motion model in traffic equilibria. The approach to deriving necessary optimality conditions in \cite{b1} was based on the method of discrete approximations developed in \cite{m95,bori} for optimal control of Lipschitzian differential inclusions and then extended in \cite{ccmn,ca3,b1,chhm3,cmn,hm,m19,mn} to various control systems governed by discontinuous differential inclusions of the sweeping type.\vspace*{-0.05in}

Here we conduct a detailed study of the formulated optimal control problem $(P)$ for integro-differential sweeping processes with proving the {\em existence of optimal solutions} and deriving comprehensive {\em necessary optimality conditions} for a broad family of local minimizers of $(P)$, where the obtained conditions are expressed entirely in terms of the problem data. On one side, the achieved results extend those from \cite{b1} to the new class of controlled sweeping processes. On the other hand, we establish a novel necessary optimality condition of the {\em Volterra type}, which is characteristic for integro-differential sweeping control systems while being particularly useful for calculations of optimal solutions.\vspace*{-0.05in}

To reach our goals, we develop the {\em method of discrete approximations} in the new setting of dynamical systems governed by controlled sweeping integro-differential inclusions, with justifying the {\em well-posedness} of discrete approximations in the sense of establishing the $W^{1,2}$-strong approximation of feasible solutions for $(P)$ by their extended discrete counterparts as well as verifying the $W^{1,2}$-strong convergence of discrete optimal solutions to a prescribed local minimizer of $(P)$. The results obtained in this direction are of their own interest (including numerical issues), while they are exploited in the paper as a {\em driving force} to derive necessary optimality conditions in problem $(P)$ by doing this first for the discrete-time problems and then by passing to the limit from them with the diminishing discretization step. To proceed in such a way, we need---by taking into account the very structure of the integro-differential sweeping dynamics in \eqref{e:1}---to employ appropriate tools of first-order and (mainly) {\em second-order variational analysis} and {\em generalized differentiation}. It occurs that the best pick for the needed constructions are those introduced by the third author and then developed in many publications; see Section~\ref{sec:2va} for more discussions and references. Moving in this direction allows us to establish below a comprehensive collection of necessary optimality conditions for problem $(P)$ and its discrete approximations, which are of their independent benefits. The given applications to the formulated {\em optimal control models} for {\em non-regular electric circuits} with numerical calculations illustrate the efficiency of the obtained necessary optimality conditions to solve particular control problems of the integro-differential type $(P)$ that naturally appear in practical modeling.\vspace*{-0.05in}

The rest of the paper is organized as follows. Section~\ref{sec:prel} describes the {\em standing assumptions} used below and presents the required {\em preliminary results}. In Section~\ref{sec:exist} we first establish the {\em existence} of feasible and optimal solutions to problem $(P)$ and then define and discuss the notion of {\em local minimizers} for which the necessary optimality conditions are derived below.\vspace*{-0.05in}

Section~\ref{sec:disc-app} is devoted to the construction of {\em discrete approximations} for controlled integro-differential sweeping processes \eqref{e:1} and to the proof that any {\em feasible solution} to \eqref{e:1} can be {\em $W^{1,2}$-strongly approximated} by feasible solutions to discrete problems, which are piecewise linearly extended to the whole interval $[0,T]$. The obtained crucial result goes far beyond optimization and occurs to be useful as an efficient machinery of the qualitative and numerical analysis of discontinuous integro-differential inclusions of the sweeping type.\vspace*{-0.05in}

Section~\ref{sec:disc-opt} continues the discrete approximation developments of Section~\ref{sec:disc-app} while now concentrating on the approximation of the entire problem $(P)$ and its prescribed {\em local minimizer} by optimal solutions to discrete-time problems. Here we show that the constructed discrete approximations always admit optimal solutions whose extensions on $[0,T]$ {\em strongly converge} to the given local minimizer under in the $W^{1,2}$ topology.\vspace*{-0.05in}

To proceed with deriving necessary optimality conditions, we recall in Section~\ref{sec:2va} the basic {\em generalized differential constructions} of variational analysis that are needed for our study. Although all the mappings involved in the description of $(P)$ are assumed to be smooth, the unavoidable source of {\em nonsmoothness} comes from the {\em sweeping dynamics} in \eqref{e:1}, which requires the usage of appropriate {\em second-order} constructions applied to (nonconvex) {\em graphs} of the normal cones. We review the employed constructions of generalized differentiation and present calculation formulas for them expressed entirely via the given data.\vspace*{-0.05in}

In Section~\ref{sec:disc-nec} we derive {\em necessary optimality conditions} for the {\em discrete approximations} of problem $(P)$ by using the generalized differential constructions of Section~\ref{sec:2va}, their well-developed calculi, and the second-order computations. The obtained results are important for their own sake as necessary optimality conditions for discrete-time counterparts of the controlled integro-differential sweeping processes. Furthermore, the strong convergence of discrete optimal solutions established in Section~\ref{sec:disc-opt} allows us to view the the obtained necessary optimality conditions for discrete approximations as {\em suboptimality conditions} for the original problem $(P)$ governed by the sweeping integro-differential inclusions.\vspace*{-0.05in}

Section~\ref{sec:nec-sw} accumulates the developments of all the previous sections and provides {\em necessary optimality conditions} for the designated class local minimizers of the {\em original optimal control problem} $(P)$ for the integro-differential sweeping processes \eqref{e:1} with the mixed state-control constraints \eqref{mixed}. By using the method of discrete approximations together with the aforementioned tools of variational analysis and generalized differentiation, we derive a comprehensive set of necessary optimality conditions of the following two types: those which extend to the integro-differential systems the recently obtained conditions for the sweeping processes governed by differential inclusions \cite{b1}, and completely novel ones that are specific for the controlled integro-differential sweeping dynamics.\vspace*{-0.05in}

The final Section~\ref{sec:exa} is devoted to applications of the obtained necessary optimality conditions for integro-differential control problems governed by integro-differential sweeping processes to some real-life models appearing in {\em non-regular electrical circuits}. We formulate two models of this type and present a realistic example showing that the obtained results allow us to determine and fully compute optimal solutions. The novel Volterra type optimality condition occurs to be especially useful for the provided computations.

Our notation is standard in variational analysis; see, e.g., \cite{b3,rw}. Recall that we distinguish between single-valued $f\colon\R^n\to\R^m$ and set-valued $F\colon\R^n\tto\R^m$ mappings. The notation $\B(x,r)$ indicates the closed ball centered at $x$ with radius $r>0$, while the symbol $\B$ stands for the unit closed ball of the space in question. We use the sign $^*$ for the matrix transposition and denote $\N:=\{1,2,\ldots\}$.\vspace*{-0.2in}
%%%%%%%%%%%%%%%%%%%%%%%%%%%%%%%%%%%%%%%%%%%%%%%%%%%%%%%%%%%%%%%%%%%%%%%%%%%%%%%%%%%%%%%%%%%%%%%%%%%%%%%%%%%%%%%%%%%%%%%%%%%%%%%%%%%%%%%%%%%%%%%%%%%%%%%%%%%%%%%%%%%%%%%%%%%%%%%%%%%%%%%%%%%%%%%%%%%%%%%%%%%%%
%%%%%%%%%%%%%%%%%%%%%%%%%%%%%%%%%%%%%%%%%%%%%%%%%%%%%%%%%%%%%%%%%%%%%%%%%%%%%%%%%%%%%%%%%%%%%%%%%%%%%%%%%%%%%%%%%%%%%%%%%%%%%%%%%%%%%%%%%%%%%%%%%
\section{Standing Assumptions and Preliminaries}\label{sec:prel}
\setcounter{equation}{0}\vspace*{-0.1in}

Let us start this section with listing the {\em standing assumptions} imposed throughout the paper unless otherwise stated. On the other hand, in some statements we specify those assumptions from the standing ones, which are actually used in the proof of this particular result.\vspace*{-0.05in}
\begin{enumerate}
\item [($\bf\mathcal{H}_{1} $)]
There are constants $ c>0$ and $ M_{j}>0 $, $ j=1,2,3 $, and open sets $ A_{i}\subset V_{i} $ such that $ d_{H}(A_{i},\mathbb{R}^{n}\setminus V_{i})>c $ and the functions $ g_{i}(\cdot) $ as $ i=1,\ldots,s $ are twice continuously differentiable satisfying the inequalities
\begin{equation}\label{e:bound-grad}
M_{1}\leq \lVert \nabla g_{i}(x) \rVert\leq M_{2},\quad\lVert \nabla^{2} g_{i}(x) \rVert\leq M_{3}\,\,\,\text{for all}\,\,\,x\in\,V_{i},
\end{equation}
where $d_H$ signifies the standard Hausdorff distance between two sets. Furthermore,
there are positive numbers $ \beta $ and $ \rho $ for which we have the estimate
\begin{equation}\label{e:w-inverse-triangle-in}
\sum\limits_{i\in I_{\rho}(x)}\lambda_{i}\lVert \nabla g_{i}(x) \rVert\leq \beta \lVert \sum\limits_{i\in I_{\rho}(x)}\lambda_{i} \nabla g_{i}(x)  \rVert \,\,\,\text{whenever}\,\,\,x\in\,C\,\,\,\text{and}\,\,\,\lambda_{i}\geq 0
\end{equation}
with the index set for the perturbed constraints defined by
\begin{equation}\label{e:rho-index}
I_{\rho}(x):=\big\{i\in\{1,\ldots,s\}\;\big|\;g_{i}(x)\leq \rho\big\} .
\end{equation}
\item  [($ \bf\mathcal{H}_{2} $)] Let $ f_{1} : \mathbb{R}^{m+n} \to \mathbb{R}^{n} $, $ f_{2} : \mathbb{R}^{d+n}\to \mathbb{R}^{n} $
be continuous mappings satisfying the following properties:
\begin{enumerate}
\item [($ \bf\mathcal{H}_{2,1} $)] There exist nonnegative constants $ \alpha_{i} $, $i=1,2$, for which we have
\begin{equation}\label{e:growthf1-con}
\lVert f_{1}(a(t),x)\rVert\leq\lVert a(t) \rVert + \alpha_{1}\lVert x \rVert\,\,\,\text{whenever}\,\,\,t\in[0,T],\,\,x\in\underset{t\in [0,T]}\bigcup C(t),
\end{equation}
\begin{equation}\label{e:growthf2-con}
\lVert f_{2}(b(t),x)\rVert\leq \lVert b(t) \rVert +\alpha_{2}\lVert x \rVert\,\,\,\text{whenever}\,\,\,t\in[0,T],\,\,x\in\underset{t\in [0,T]}\bigcup C(t).
\end{equation}
\item [($\bf\mathcal{H}_{2,2} $)] For any real numbers $ r_{i}>0 $ and functions $ a(\cdot)$, $ b(\cdot)\in W^{1,2}([0,T],\mathbb{R}^{m+d}) $
there exist nonnegative constants $L$ and $ L_{i}$ as $i=1,2$ such that we have the estimates
\begin{equation}\label{e:Lipsf1}
\lVert f_{1}(a,x_{1})-f_{1}(a,x_{2}) \rVert\leq L\lVert x_{1}-x_{2} \rVert\,\,\,\text{for all}\,\,\,x_{1},x_{2}\in r_{1}\mathbb{B},
\end{equation}
\begin{equation}\label{e:Lipsf2}
\lVert f_{2}(b,x_{1})-f_{2}(b,x_{2}) \rVert\leq L_{1}\lVert x_{1}-x_{2} \rVert\,\,\,\text{for all}\,\,\,x_{1},x_{2}\in r_{1}\mathbb{B},
\end{equation}
\begin{equation}\label{e:Lipsf2b}
\lVert f_{2}(b_{1},x)-f_{2}(b_{2},x) \rVert\leq L_{2}\lVert b_{1}-b_{2} \rVert\,\,\,\text{for all}\,\,\,b_{1},b_{2}\in \mathbb{R}^{n},\,\,x\in r_{1}\mathbb{B}.
\end{equation}
\end{enumerate}
\item [($ \bf\mathcal{H}_{3} $)] The terminal cost $ \varphi : \mathbb{R}^{n}\to\oR$ is lower semicontinuous (l.s.c.), while and the
running cost $ l_{0} : [0,T]\times\mathbb{R}^{4n+2m+2d}\to\bar{\mathbb{R}} $ is l.s.c.\ with respect to all but time variable being continuous with respect to $t$ and being majorized by a summable function on $[0,T]$ along reference curves. Furthermore, we assume that $l_0(t,\cdot)$ is bounded from below on bounded sets for a.e.\ $t\in[0,T]$.
\end{enumerate}
It is not hard to verify that if both conditions \eqref{e:bound-grad} and \eqref{e:w-inverse-triangle-in} are satisfied, then we have the {\em positive linear independence} of the gradients  $\nabla g_i(x)$ of the active inequality constraints on $C$, and that the latter condition is equivalent in our setting to the {\em Mangasarian-Fromovitz constraint qualification} of nonlinear programming.\vspace*{-0.05in}

Putting aside the case of concavity of all the functions $g_i$ in \eqref{e:mset}, it can be seen that the sets $C$ and $C(t)$ therein are {\em nonconvex}, and hence we need to clarify in which sense the normal cone in \eqref{e:1} is understood. To proceed, for any set $\O\subset\R^n$ locally closed around $\ox\in\O$, consider the (Euclidean) distance function $\dist(x;\O):=\inf_{y\in\O}\|x-y\|$ be the projection operator $\Pi_\O\colon\R^n\tto\O$ by
\begin{equation}\label{proj}
\Pi_\O(x):=\big\{w\in\O\;\big|\;\|x-w\|=\dist(x;\O)\big\},\quad x\in\R^n.
\end{equation}
The {\em proximal normal cone} to $\O$ at $\ox$ is given by
\begin{eqnarray}\label{pnc}
N^P_\O(\ox):=\big\{v\in\R^n\big|\;\exists\,\al>0\;\mbox{ such that }\;\ox\in\Pi_\O(\ox+\al v)\big\}\;\mbox{ if }\;\ox\in\O
\end{eqnarray}
and by $N^P_\O(\ox):=\emp$ otherwise; see \cite{clsw,Mord,rw} for equivalent descriptions and further references. It is easy to see that $N^P_\O(\ox)$ is a convex cone, which may not be generally closed.\vspace*{-0.05in}

Using \eqref{pnc}, recall now the following notion, which plays an important role in variational analysis and theory of sweeping processes. Given a closed set $\O\subset\R^n$ and a number $\eta>0$, it is said that $\O$ is $\eta$-{\em prox-regular} if for all $x\in\bd\O$ and $v\in N^P_\O(x)$ with $\|v\|=1$ we have $\B(x+\eta v,\eta)\cap\O=\{x\}$, which is equivalent to
\begin{equation}\label{prox}
\la v,y-x\ra\le\dfrac{\|v\|}{2\eta}\|y-x\|^2\;\mbox{ for all }\;y\in\O,\;x\in\bd\O,\;\mbox{ and }\;v\in N^P_\O(x).
\end{equation}
We refer the reader to \cite{clsw,CT,rw} for more details and history of this and related notions. It is worth mentioning
that any closed convex subset in $\R^n$ is $\eta$-prox-regular with $\eta=\infty$.

The next result is taken from \cite[Proposition~2.9]{venl}; see also \cite{ant} for more discussions and developments.\vspace*{-0.1in}

\begin{proposition}[\bf prox-regularity of moving sets]\label{Pr:prox-reg} Let assumption $(\mathcal{H}_{1})$ be satisfied. Then for each $t\in[0,T]$ we have that the set $C(t)$ is $\eta$-prox-regular with $\eta=\dfrac{M_1}{M_3\be}$.
\end{proposition}\vspace{-0.15in}

Indeed, this statement was justified in \cite[Proposition~2.9]{venl} for the set $C$ in \eqref{e:mset}, and hence it holds for the moving set $C(t)=C+u(t)$ as a translation of $C$.

It has been well realized that the prox-regularity of closed sets ensures that all the major cones of variational analysis (Clarke, Fr\'echet, Mordukhovich) agree with each other and reduce to \eqref{pnc}; see, e.g., \cite{clsw,Mord,rw}. This allows us in \eqref{e:1} the notation `$N$' under the standing assumptions.

We close this section with the following well-known discrete version of Gronwall's inequality; see, e.g., \cite{clar} for a short proof and discussions.\vspace*{-0.1in}

\begin{proposition}[discrete Gronwall's inequality]\label{granwal}
Let $ \alpha_{i} $, $ \rho_{j} $, $ \gamma_{j} $, $ a_{j} $ be nonnegative numbers with
\begin{equation*}
e_{j+1}\leq \sigma_{j}+\rho_{j}\sum\limits_{k=0}^{j-1}e_{k}+(1+\gamma_{j})e_{j}\;\mbox{ for all }\;j\in\N.
\end{equation*}
Then whenever $ i\in\mathbb{N}$ we have the estimate
\begin{equation*}
e_{i}\leq (e_{0}+\sum\limits_{k=0}^{i-1}\sigma_{k})\exp\Big(\sum\limits_{k=0}^{i-1}(k\rho_{k}+\gamma_{k})\Big) .
\end{equation*}
\end{proposition}\vspace*{-0.2in}
%%%%%%%%%%%%%%%%%%%%%%%%%%%%%%%%%%%%%%%%%%%%%%%%%%%%%%%%%%%%%%%%%%%%%%%%%%%%%%%%%%%%%%%%%%%%%%%%%%%%%%%%%%%%%%%%%%%%%%%%%%%%%%%%%%%%%%%%%%%%%%%%%%%%%%%%%%%%%%%%%%%%%%%%%%%%%%%%%%%%%%%%%%%%%%%%%%%%%%%%%%%%%
%%%%%%%%%%%%%%%%%%%%%%%%%%%%%%%%%%%%%%%%%%%%%%%%%%%%%%%%%%%%%%%%%%%%%%%%%%%%%%%%%%%%%%%%%%%%%%%%%%%%%%%%%%%%%%%%%%%%%%%%%%%%%%%%%%%%%%%%%%%%%%%%%

\section{Existence of Optimal Solutions and Local Minimizers}\label{sec:exist}
\setcounter{equation}{0}\vspace*{-0.1in}

In this section we establish the existence of optimal solutions to $(P)$ and describe the type of local minimizers of $(P)$ used below for deriving necessary optimality conditions. But first we present the following {\em well-posedness} results for the controlled sweeping dynamics \eqref{e:1} that ensure, in particular the existence of {\em feasible} solutions $(x(\cdot),u(\cdot),a(\cdot),b(\cdot))$ in the corresponding $W^{1,2}$ spaces.\vspace{-0.1in}

\begin{proposition}[\bf well-posedness of integro-differential sweeping processes]\label{exist:solu} Fix a triple\\ $(u(\cdot),a(\cdot),b(\cdot))\in W^{1,2}([0,T],\mathbb{R}^{n+m+d})$ and consider the integro-differential sweeping process \eqref{e:1} under the assumptions in $(\mathcal{H}_{1})$ and $(\mathcal{H}_{2})$. Then this control triple generates a unique trajectory $x(\cdot)\in W^{1,2}([0,T],\mathbb{R}^{n})$ of system \eqref{e:1}. Denoting further
\begin{equation*}
\beta_{1}(t):=\max\big\{\lVert a(t) \rVert , \alpha_{1}\big\},\quad\tilde{b}(t):=2\max\big\{\beta_{1}(t),\alpha_{2}\big\}\,\,\,\text{for all}\,\,\,t\in[0,T],
\end{equation*}
\begin{equation*}
\tilde{l}:=\lVert x(0) \rVert\exp\Big(\int\limits_{0}^{T}(\tilde{b}(\tau)+1)\,d\tau\Big)+\exp\Big(\int\limits_{0}^{T}(\tilde{b}(\tau)+1)\,d\tau\Big)\int\limits_{0}^{T}\Big(\rVert\dot{u}(s)\rVert +2\beta_{1}(s)+2\int\limits_{0}^{T}\lVert b(\tau) \rVert\,d\tau\Big)\,ds,
\end{equation*}
and taking the corresponding arc $y(\cdot)$ from \eqref{e:2}, we have the estimates
\begin{equation}\label{e:3}
\lVert \dot{x}(t)+f_{1}(t,x(t))+y(t) \rVert\leq \lVert \dot{u}(t)\rVert+(1+\tilde{l})\beta_{1}(t)+\int\limits_{0}^{t}\lVert b(s) \rVert\,ds+T\alpha_{2}\tilde{l},\,\,\,\text{a.e.}\,\,\,t\in[0,T],
\end{equation}
\begin{equation}\label{e:4}
\lVert \dot{x}(t) \rVert\leq \lVert \dot{u}(t)\rVert+2(1+\tilde{l})\beta_{1}(t)+2\int\limits_{0}^{t}\lVert \tilde{b}(s) \rVert\,ds+2T\alpha_{2}\tilde{l}\,\,\,\text{a.e.}\,\,\,t\in[0,T],
\end{equation}
\begin{equation}\label{e:5}
\lVert \dot{y}(t) \rVert\leq \lVert b(t) \rVert + \alpha_{2}\tilde{l}\,\,\,\text{for all}\,\,\,t\in[0,T].
\end{equation}
\end{proposition}\vspace*{-0.15in}
{\bf Proof}. It is not hard to observe from the construction of $C(t)$ in \eqref{e:mset} that
\begin{equation}\label{haus}
\begin{array}{ll}
C(u(t))\subset C(u(s))+\lvert \upsilon(t)-\upsilon(s) \rvert\mathbb{B},\;\; \text{for all}\;\; t,s\in[0,T]\;\; \text{and}\;\; u(t)\in\mathbb{R}^{n},
\end{array}
\end{equation}
where $ \upsilon(t):=\displaystyle\int\limits_{0}^{t}\lVert \dot{u}(s) \rVert\,ds$. This implies that
\begin{equation*}
C(t) \subset C(s)+ \lvert \upsilon(t)-\upsilon(s) \rvert\mathbb{B}\;\mbox{ on }\;[0,T].
\end{equation*}
Furthermore, by the imposed condition $(\mathcal{H}_{1})$ and Proposition~\ref{Pr:prox-reg} we have that $C(t)$ is $\eta$-prox-regular for each $t\in[0,T]$. If in addition $(\mathcal{H}_{2,1})$ holds, then this ensures that
\begin{equation*}
\lVert f_{1}(a(t),x) \rVert\leq \lVert a(t) \rVert+\alpha_{1}\lVert x \rVert\leq \max\big\{ \lVert a(t) \rVert , \alpha_{1}\big\}(1+\lVert x \rVert)=\beta_{1}(t)(1+\lVert x \rVert).
\end{equation*}
Hence all the assumptions of \cite[Theorem~4.2]{bou} are satisfied, which thus yields the fulfillment of estimates \eqref{e:3}--\eqref{e:5}. Observe finally that $ \dot{x}(\cdot),\dot{y}(\cdot)\in L^{2}([0,T],\mathbb{R}^{n})$ since $\beta_{1}(\cdot)\in L^{2}([0,T],\mathbb{R}^{n})$ by the construction of $\beta_1(\cdot)$. This therefore completes the proof of the proposition. $\h$

Consider now the set-valued mapping $ F_{1} : \mathbb{R}^{3n+m}\rightrightarrows\mathbb{R}^{n} $ given by
\begin{equation}\label{F1}
F_{1}(x,y,u,a):=N_{C(u)}(x)+f_{1}(x,a)+y.
\end{equation}
It is easy to deduce from the representation of the set $C$ in \eqref{e:mset} that
\begin{equation*}
F_{1}(x,y,u,a):=\Big\{-\sum\limits_{i\in I(x-u)}\eta_{i}\nabla g_{i}(x-u)\;\Big|\;0\leq \eta_{i}\Big\}+f_{1}(x,a)+y,
\end{equation*}
(cf.\ \cite[Proposition~2.8]{venl}), where the set of {\em active constraint indices} is
\begin{equation}\label{aci}
I(y):=\big\{i\in\{1,\ldots,s\}\;\big|\;g_{i}(y)=0\big\}
\end{equation}
Define further the mappings $F : \mathbb{R}^{3n+m+d}\rightrightarrows\mathbb{R}^{n} $ and $ \tilde{f_{2}} : \mathbb{R}^{3n+m+d}\to\mathbb{R}^{n}$ by
\begin{equation*}
F(z):=F(x,y,u,a,b)=F_{1}(x,y,u,a)=N_{C(u)}(x)+f_{1}(x,a)+y,
\end{equation*}
\begin{equation*}
\tilde{f_{2}}(z):=\tilde{f_{2}}(x,y,u,a,b)=-f_{2}(b,x).
\end{equation*}
It is convenient to rewrite the Cauchy problem for the evolution system in \eqref{e:2} in terms of $z\in\mathbb{R}^{3n+m+d}$ as
\begin{equation}\label{e:7}
-\dot{z}(t)\in F(z(t))\times \tilde{f}_{2}(z(t))\times\mathbb{R}^{n}\times\mathbb{R}^{m}\times\mathbb{R}^{d}:=G(z(t))\;\mbox{ a.e. }\;t\in[0,T]
\end{equation}
with $\dot{z}(t)=(\dot{x}(t),\dot{y}(t),\dot{u}(t),\dot{a}(t),\dot{b}(t))$ and the initial conditions
\begin{equation}\label{e:8}
z(0):=(x_{0},y_0,u_0,a_0,b_0),\;\mbox{ where }\;y_0=0,\;(a_0,b_0)\in\R^m\times\R^d,\;\mbox{ and }\;g_i(x_0-u_0)\ge 0\;\mbox{ for all }\;i=1,\ldots,s.
\end{equation}\vspace*{-0.25in}

Now we are ready to obtain conditions ensuring the existence of (global) optimal solutions to the sweeping control problem $(P)$ governed by integro-differential inclusions. In what follows we'll freely use any of the three equivalent forms of $(P)$ depending on the situation.\vspace*{-0.1in}
\begin{theorem}[\bf existence of sweeping optimal solutions]\label{exist:opti}  In addition to the standing assumptions $(\mathcal{H}_{1})-(\mathcal{H}_{3})$, suppose that the running cost $l_0$ from \eqref{e:running cost} is convex with respect of velocity variables $\dot{x},\dot{u},\dot{a},\dot{b}$, and that there is a minimizing sequence $\{(x^{k}(\cdot),y^k(\cdot),u^{k}(\cdot),a^{k}(\cdot),b^{k}(\cdot))\}$ in $(P)$, which $\{(u^{k}(\cdot),a^{k}(\cdot),b^{k}(\cdot))\}$ is bounded in $W^{1,2}([0,T],\mathbb{R}^{n+m+d})$, such that the integrand $l_{0}(t,\cdot)$ is majorized by a summable function on $[0,T]$. Then problem $(P)$ admits an optimal solution belonging to the space $ W^{1,2}([0,T],\mathbb{R}^{3n+m+d})$.
\end{theorem}\vspace*{-0.15in}
{\bf Proof}. Proposition \ref{exist:solu} tells us that the set of feasible solutions to problem $(P)$ is nonempty. Let us fix a minimizing sequence of feasible solutions $z^{k}(\cdot)=(x^{k}(\cdot),y^k(\cdot), u^{k}(\cdot),a^{k}(\cdot),b^{k}(\cdot))$ for $(P)$, which is bounded in $W^{1,2}([0,T],\mathbb{R}^{3n+m+d})$. This implies, in particular, that there exists a triple $(u_0,a_0,b_0)\in\mathbb{R}^{n+m+d}$ such that $(u^{k}(0),a^{k}(0),b^{k}(0))\to(u_{0},a_{0},b_{0})$ in this space as $k\to\infty$, while the quintuple $(x_0,y_0,u_0,a_0,b_0)$ clearly satisfies \eqref{e:8}. We readily have that the sequence $\{(\dot{u}^{k}(\cdot),\dot{a}^{k}(\cdot),\dot{b}^{k}(\cdot))\}$ is bounded in $L^2([0,T],\R^{n+m+d}$. Since the latter space is reflexive, where bounded sets are weakly compact, there exists a triple $(\upsilon^{u}(\cdot),\upsilon^{a}(\cdot),\upsilon^{b}(\cdot))\in L^{2}([0,T],\mathbb{R}^{n+m+d}$ such that $(\dot{u}_{k}(\cdot),\dot{a}_{k}(\cdot),\dot{b}_{k}(\cdot))\to(\upsilon^{u}(\cdot),\upsilon^{a}(\cdot),\upsilon^{b}(\cdot))\in L^{2}([0,T],\mathbb{R}^{n+m+d}$ in the weak topology of this space along a subsequence (without relabeling). Define now the functions
\begin{equation*}
(\bar{u}(t),\bar{a}(t),\bar{b}(t))=(u_{0},a_{0},b_{0})+\int\limits_{0}^{t}(\upsilon^{u}(s),\upsilon^{a}(s),\upsilon^{b}(s))\,ds\;\mbox{ for all }\;t\in[0,T]
\end{equation*}
and observe that $(\dot{\ou}(t),\dot{\oa}(t),\dot{\ob}(t))=(\upsilon^{u}(t),\upsilon^{a}(t),\upsilon^{b}(t))$ for a.e.\ $t\in[0,T]$, and that the triple $(\bar{u}(\cdot),\bar{a}(\cdot),\bar{b}(\cdot))$ belongs to the space $W^{1,2}([0,T],\mathbb{R}^{n+m+d})$. It follows from the above and the estimates of Proposition~\ref{exist:solu} that the sequence of the corresponding trajectories $\{x^{k}(\cdot)\}$ is uniformly bounded and equicontinuous on $[0,T]$. Then the Arzel\'a-Ascoli theorem allows us to find a subsequence of $\{x^{k}(\cdot)\}$ that uniformly converges on $[0,T]$, without relabeling, to some arc $\bar{x}(\cdot)\in{\cal C}([0,T],\mathbb{R}^{n})$, which is in fact absolutely continuous on this interval.
Turning further to the $y$-component of the minimizing sequence, we deduce from the Lebesgue dominated convergence theorem and the continuity of $f_{2} $ that $ y^{k}(t)\to \bar{y}(t)$ for all $t\in[0,T]$ and that
\begin{equation*}
\dot{y}^{k}(t)\to\dot{\bar{y}}(t)\;\mbox{ for a.e. }\;t\in[0,T],\;\mbox{ where }\;\bar{y}(t):=\displaystyle\int\limits_{0}^{t}f_{2}(\bar{b}(s),\bar{x}(s))\,ds\;\mbox{ with }\;\dot{\bar{y}}(\cdot)\in L^{2}([0,T],\mathbb{R}^{n}).
\end{equation*}
It follows from \eqref{e:4} that $\{\dot{x}^{k}(\cdot)\}$ is bounded in $ L^{2}([0,T],\mathbb{R}^{n})$, and hence it weakly converges in $L^{2}([0,T],\mathbb{R}^{n})$, up to a
subsequence, to some function $w(\cdot)$ with $\dot{\bar{x}}(t)=w(t)\in L^{2}([0,T],\mathbb{R}^{n})$ for a.e.\ $t\in[0,T]$.\vspace*{-0.05in}

The next step is to verify that the limiting quintuple $\bar{z}(\cdot)=(\bar{x}(\cdot),\bar{y}(\cdot),\bar{u}(\cdot),\bar{a}(\cdot),\bar{b}(\cdot)) $ satisfies the combination of the differential inclusion and equation in \eqref{e:7}, which is equivalent to the original integro-differential inclusion in \eqref{e:1}. Since the derivative sequences $\{\dot z^k(\cdot)\}$ converges to $\dot\oz(\cdot)$ weakly in $L^2([0,T];\R^{3n+m+d})$, the classical Mazur theorem ensures the strong convergence to $\dot\oz(\cdot)$ in this space of some sequence of convex combinations of the functions $\dot z^k(t)$. Thus there is a subsequence of these convex combinations that converges to $\dot\oz(t)$ for a.e. $t\in[0,T]$. It follows from the above that there exist a function $\nu:\N\to\N$ and sequences of real numbers $\{\alpha(k)_j\}$, $j=k,\ldots,\nu(k)$, such that
$$
\alpha(k)_j\ge 0,\;\sum_{j=k}^{\nu(k)}\alpha(k)_j=1,\;\mbox{ and }\;\sum_{j=k}^{\nu(k)}\alpha(k)_j\dot{z}^j(t)\to\dot{\oz}(t)\;\mbox{ a.e. }\;t\in [0,T]
$$
as $k\to\infty$. By taking into account that $\ox(t)-\ou(t)=\disp\lim_{k\to\infty}({x}^{k}(t)-{u}^{k}(t))\in C$ we get
\begin{equation*}
\begin{aligned}
-\dot{\ox}(t)-&f_{1}\big(\oa(t),\ox(t)\big)-\oy(t)=\lim_{k\to\infty}\Big(-\sum_{j=k}^{\nu(k)}\alpha(k)_j\dot{x}^j(t)-\sum_{j=k}^{\nu(k)}\alpha(k)_j f_{1}\big(a^j(t),x^j(t)\big)-\sum_{j=k}^{\nu(k)}\alpha(k)_jy^{j}(t)\Big)\\
=&\lim_{k\to\infty}\Big(-\sum_{j=k}^{\nu(k)}\sum_{i\in I(x^j(t)-u^j(t))}\al(k)_j\eta^j_i\nabla g_i\big(x^j(t)-u^j(t)\big)\Big)\\
=&\lim_{k\to\infty}\Big(-\sum_{j=k}^{\nu(k)}\sum_{i\in I(\ox(t)-\ou(t))}\al(k)_j\eta^j_i\nabla g_i\big(x^j(t)-u^j(t)\big)\Big)\\
=&\lim_{k\to\infty}\Big(-\sum_{i\in I(\ox(t)-\ou(t))}\sum_{j=k}^{\nu(k)}\al(k)_j\eta^j_i\nabla g_i\big(x^j(t)-u^j(t)\big)\Big),
\end{aligned}
\end{equation*}
where $I(\cdot)$ is taken from \eqref{aci}, and where $\eta^{j}_{i}=0$ if $i\in I(\ox(t)-\ou(t))\backslash I(x^j(t)-u^j(t))$ due to the clear inclusion $I(x^j(t)-u^j(t))\subset I(\ox(t)-\ou(t))$ for all $j=k,\ldots,\nu(k)$ and all large $k\in\N$. Proceeding now similarly to the proof of \cite[Theorem~4.1]{b1} shows that
$$
\disp\sum^{\nu(k)}_{j=k}\al(k)_j\eta^j_i \to\be_i\;\;\text{as}\;\;k\to\infty
$$
along some subsequence (without relabeling), and that
$$
\sum_{j=k}^{\nu(k)}\al(k)_j\eta^j_i\nabla g_i(x^j(t)-u^j(t))\to\be_i\nabla g_i(\ox(t)-\ou(t))\;\; \text{as}\;\; k\to\infty.
$$
Thus we have at the representations
\begin{equation*}
\begin{aligned}
-\dot{\ox}(t)-&f_{1}\big(\oa(t),\ox(t)\big)-\oy(t)=\lim_{k\to\infty}\Big(-\sum_{i\in I(\ox(t)-\ou(t))}\sum_{j=k}^{\nu(k)}\al(k)_j\eta^j_i\nabla g_i\big(x^j(t)-u^j(t)\big)\Big)\\
&=-\sum_{i\in I(\ox(t)-\ou(t))}\be_i\nabla g_i\big(\ox(t)-\ou(t)\big)\in N_{C}\big(\ox(t)-\ou(t)\big)\;\mbox{ a.e. }\;t\in[0,T],
\end{aligned}
\end{equation*}
where $\bar{y}(t):=\displaystyle\int\limits_{0}^{t}f_{2}(\bar{b}(s),\bar{x}(s))\,ds $. This verifies the fulfillment of the inclusion in \eqref{e:7}.\vspace*{-0.05in}

To justify further the optimality of $ \bar{z}(\cdot)=(\bar{x}(\cdot),\bar{y}(\cdot),\ou(\cdot),\bar{a}(\cdot),\bar{b}(\cdot))$ in $(P)$, it is sufficient to show that
\begin{equation}\label{exis-opt}
J[\bar{x},\bar{y},\bar{u},\bar{a},\bar{b}]=J_{0}[\bar{x},\bar{u},\bar{a},\bar{b}]\leq\liminf\limits_{k\to \infty} J_{0}[x_{k},u_{k},a_{k},b_{k}]=\liminf\limits_{k\to \infty} J[x_{k},y_{k},u_{k},a_{k},b_{k}]
\end{equation}
for the Bolza functionals in \eqref{J0} and \eqref{e:bolza}. The latter readily follows from the assumptions in ($\mathcal{H}_{3}$) on the cost functions $\ph$ and $l_0$ due to the Mazur weak closure theorem and the Lebesgue dominated convergence theorem. Indeed, Mazur's theorem ensures that the weak convergence of the $\{(\dot{x}_k(\cdot),\dot{u}_k(\cdot),\dot{a}_k(\cdot),\dot{b}_k(\cdot))\}$ to $(\dot{\ox}(\cdot),\dot{\ou}(\cdot),\dot{\oa}(\cdot),\dot{\ob}(\cdot))$ in $L^2([0,T],\mathbb{R}^{2n+m+d})$ yields the $L^2$-strong convergence of convex combinations of $(\dot{x}_k(\cdot),\dot{u}_k(\cdot),\dot{a}_k(\cdot),\dot{b}_k(\cdot))$ to $(\dot{\ox}(\cdot),\dot{\ou}(\cdot),\dot{\oa}(\cdot),\dot{\ob}(\cdot))$, and thus the a.e.\ convergence of a subsequence of these convex combination on $[0,T]$ to the limiting quadruple. Employing finally the assumed convexity of the running cost $l_0$ with respect to the velocity variables verifies \eqref{exis-opt} and thus completes the proof of the theorem. $\h$

We can see that, besides our standing fairly unrestrictive assumptions, the existence of global minimizers in Theorem~\ref{exist:opti} requires the {\em convexity} of the running cost with respect to {\em velocity} variables. The obtained results is new, while this convexity phenomenon has been well recognized in the calculus of variations and optimal control of various types of dynamical systems. On the other hand, it has been also well understood in variational theory for continuous-type systems (including sweeping processes) that such problems allows a certain {\em relaxation} procedure involving the velocity {\em convexification}, which brings us to relaxed problems, where optimal solutions exist automatically and can be constructively approximated by feasible solutions to the original problems with keeping the same optimal values of the cost functionals. The reader is referred to \cite{b1,clsw,chhm3,dfm,edm,m95,bori,tols,v} for various results, discussions, and references in these directions. Following this line, we construct now the relaxed version of the optimal control problem $(P)$ for the integro-differential sweeping process under consideration.\vspace*{-0.05in}

Taking the integrand $l$ from \eqref{e:bolza} and the velocity mapping $G$ from \eqref{e:7}, define the {\em extended running cost}
\begin{equation*}
l_{G}(t,x,y,u,a,b,\upsilon):=l(t,x,y,u,a,b,\upsilon)+\delta_G(\upsilon)\;\mbox{ with }\;G=G(x,y,u,a,b),
\end{equation*}
where the indicator function $\dd_G$ of the set $G$ is given by $\dd_G(\upsilon):=0$ if $\upsilon\in G$ and $\dd_G(\upsilon):=\infty$ otherwise. Denote further by $\hat{l}_{G}(t,x,y,u,a,b,\upsilon)$ the {\em biconjugate function} to $l_{G}(t,x,y,u,a,b,\upsilon)$ with respect to the velocity variable $\upsilon=(\dot{x},\dot{u},\dot{a},\dot{b})$, i.e., given by
\begin{equation*}
\hat{l}_{G}(t,x,y,u,a,b,\upsilon):=(l_{G})^{**}(t,x,y,u,a,b,\upsilon).
\end{equation*}
Observe that $ \hat{l}_{G}(t,x,y,u,a,b,\upsilon)$ is the largest proper, convex, and l.s.c.\ function with respect to $\upsilon$, which is majorized by $l_{G}$. We clearly have that
$\hat{l}_{G}=l_{G}$ if and only if $l_{G}$ is proper, convex, and l.s.c.\ with respect to $\upsilon$.

Now we are ready to define the {\em relaxed problem} $(R)$ associated with the original optimal control problem $(P)$ for sweeping integro-differential inclusion as follows: minimize
\begin{equation}\label{R}
\hat{J}[x,y,u,a,b]:=\varphi(x(T))+\displaystyle\int\limits_{0}^{T}\hat{l}_{G}(t,x(t),y(t),u(t),a(t),b(t),\dot{x}(t),\dot{y}(t),\dot{u}(t),\dot{a}(t),\dot{b}(t))\,dt
\end{equation}
with $x(0)=x_0\in C(0)$ and $y(0)=0$. Theorem~\ref{exist:opti} ensures the existence of optimal solutions to $(R)$ under the standing assumptions made. We see furthermore that there is no difference between problems $(P)$ and $(R)$ if the original running cost $l_0$ is {\em convex with respect to the velocity} variables. In fact, there exists a deeper connection between $(P)$ and $(R)$ {\em without any convexity requirements}, which has been well recognized for particular cases of the controlled sweeping processes in \cite{edm,tols} by showing that an optimal relaxed solution can be constructively approximated by feasible original ones, and that the optimal cost values for $(P)$ and $(R)$ coincide. Our goal is to investigate this issue for more general cases of integro-differential sweeping control systems in the future research.\vspace*{-0.05in}

We conclude this section by formulating the concept of {\em local minimizers} in $(P)$ for which we derive necessary optimality conditions without any convexity assumptions.\vspace*{-0.1in}
\begin{definition}[\bf intermediate local minimizers]\label{locmin} Consider problem $(P)$ with representation \eqref{e:7} of the integro-differential inclusion, and let \eqref{R} is the relaxed problem of $(P)$. Then we say that:\\[0.08ex]
{\bf(i)} A feasible solution $\oz(\cdot)$ to $(P)$ is an {\sc intermediate local minimizer} $($i.l.m.$)$ of this problem if there exists $\ve>0$ such that $J[\oz]\le J[z]$ for any feasible solution $z(\cdot)$ to $(P)$ satisfying the $W^{1,2}$-localization condition
\begin{equation}\label{loc}
\lVert z(t)-\bar{z}(t)\rVert<\varepsilon\;\mbox{ on }\;[0,T]\;\mbox{ and }\;\int\limits_{0}^{T}\lVert\dot{z}(t)-\dot{\bar{z}}(t) \rVert^{2}\,dt<\varepsilon.
\end{equation}
{\bf(ii)} A feasible solution $\oz(\cdot)$ to $(P)$ is a {\sc relaxed intermediate local minimizer} $($r.i.l.m.$)$ of this problem if $J[\bar{z}]=\hat{J}[\bar{z}]$ and there exists $\varepsilon>0$ such that $J[\bar{z}]\le J[z]$ for any feasible solution $ z(\cdot) $ to $(P)$ satisfying the $W^{1,2}$-localization condition \eqref{loc}.
\end{definition}\vspace*{-0.15in}

Both notions in Definition~\ref{locmin} were introduced in \cite{m95} for Lipschitzian differential inclusions (see also \cite{bori,v} for subsequent studies in the Lipschitzian case), and then they were investigated and developed in \cite{ccmn,ca3,b1,chhm3,cmn,hm,mn} for various differential inclusions of the sweeping type. The introduced notions include their {\em strong} counterpart in which only the first condition from \eqref{loc} is required. In general, the defined ``intermediate" notions occupy an intermediate position between strong and weak local minimizers in variational problems.
As seen, the only difference between i.l.m.\ and r.i.l.m.\ lies in the additional requirement on the {\em local relaxation stability} $J[\bar{z}]=\hat{J}[\bar{z}]$, which is often the case (always for the Lipschitzian dynamics \cite{bori} and more as in \cite{dfm,edm,tols}) of {\em nonconvex} integro-differential systems, particularly for strong minimizers.\vspace*{-0.2in}
%%%%%%%%%%%%%%%%%%%%%%%%%%%%%%%%%%%%%%%%%%%%%%%%%%%%%%%%%%%%%%%%%%%%%%%%%%%%%%%%%%%%%%%%%%%%%%%%%%%%%%%%%%%%%%%%%%%%%%%%%%%%%%%%%%%%%%%%%%%%%%%%%%%%%%%%%%%%%%%%%%%%%%%%%%%%%%%%%%%%%%%%%%%%%%%%%%%%%%%%%%%%%
%%%%%%%%%%%%%%%%%%%%%%%%%%%%%%%%%%%%%%%%%%%%%%%%%%%%%%%%%%%%%%%%%%%%%%%%%%%%%%%%%%%%%%%%%%%%%%%%%%%%%%%%%%%%%%%%%%%%%%%%%%%%%%%%%%%%%%%%%%%%%%%%%
\section{Discrete Approximations of Integro-Differential Dynamics}\label{sec:disc-app}
\setcounter{equation}{0}\vspace*{-0.1in}

In this section we start developing the {\em method of discrete approximations} to study controlled integro-differential sweeping processes of type \eqref{e:1}. The major result obtained here deals, for the fist time in the literature, with the discontinuous integro-differential sweeping dynamics independently of its optimization as in problem $(P)$. We'll later apply it to deriving necessary optimality conditions for intermediate local minimizers of $(P)$, but so far our goal is to construct a well-posed sequence of discrete-time sweeping dynamical systems, which {\em $W^{1,2}$-strongly} approximates {\em any feasible} solution to $(P)$. We can do this under the standing assumptions in $({\cal H}_1)$ and (${\cal H}_2)$ with the quite natural additional requirement that the velocities $\dot{\ox}(\cdot)$ and $\dot{\ou}(\cdot)$ of the given feasible solution to $(P)$ are of {\em bounded variation} on $[0,T]$. The developed approach allows us to improve the corresponding results of \cite{b1} for controlled sweeping processes governed by differential inclusions and extend them to the case of integro-differential dynamics. The obtained results are {\em fully constructive} and can be treated from both {\em qualitative} and {\em numerical} viewpoints as the justification of {\em finite-dimensional} approximations of infinite-dimensional integro-differential discontinuous systems.\vspace*{-0.05in}

For each $k\in\N$, consider the {\em discrete mesh} on $[0,T]$ given by
\begin{equation}\label{mesh}
\Delta_k:=\big\{0=t^k_0<t^k_1<\ldots<t^k_k=T\big\}\;\mbox{ with }\;h^k_j:=t^{k+1}_k-t^k_j\;\mbox{ and }\;\max_{j=0,\ldots,k-1}{h^k_j}\le h_k:=\frac{T}{k}.
\end{equation}\vspace*{-0.25in}
\begin{theorem}[strong discrete approximation of feasible sweeping solutions]\label{strong:conver} Under the fulfillment of the assumptions in $({\cal H}_1)$ and $({\cal H}_2)$, fix any feasible solution $(\bar{x}(\cdot),\bar{y}(\cdot),\bar{u}(\cdot),\bar{a}(\cdot),\bar{b}(\cdot))$ to $(P)$ such that the velocity functions $\dot{\bar{x}}(\cdot)$ and $ \dot{\bar{u}}(\cdot) $ are of bounded variation on $[0,T] $, i.e., there is a number $\mu>0$ for which
\begin{equation*}
\max\big\{{\rm var}(\dot{\bar{x}}(\cdot),[0,T]),\,{\rm var}(\dot{\bar{u}}(\cdot),[0,T])\big\}\le\mu.
\end{equation*}
Then given a discrete mesh $\Delta_{k}$ in \eqref{mesh}, there exist sequences of piecewise linear functions $(x^{k}(\cdot),y^{k}(\cdot),\\u^{k}(\cdot),a^{k}(\cdot),b^{k}(\cdot))$ such that
$(x^{k}(0),y^{k}(0),u^{k}(0),a^{k}(0),b^{k}(0))=(x_{0},0,\ou(0),\oa(0),\ob(0))$ for all $k\in\mathbb{N}$, and
\begin{equation*}
x^{k}(t)=x^{k}(t^{k}_{j})+(t-t^{k}_{j})\upsilon^{k}_{j},\quad y^{k}(t)=y^{k}(t^{k}_{j})+(t-t^{k}_{j})w^{k}_{j}\;\mbox{ for }\;t\in[t^{k}_{j},t^{k}_{j+1}],\;j=0,\ldots,k-1,
\end{equation*}
where the discrete velocities $\upsilon^{k}_{j}$ and $w^{k}_{j}$ satisfy the conditions
\begin{equation*}
-\upsilon^{k}_{j}\in F_{1}(x^{k}(t^{k}_{j}),y^{k}(t^{k}_{j+1}),u^{k}(t^{k}_{j}),a^{k}(t^{k}_{j})),\quad w^{k}_{j}= f_{2}(b^{k}(t^{k}_{j}),x^{k}(t^{k}_{j})),\,\,\,j=0,\ldots,k-1,
\end{equation*}
with $F_1$ taken from \eqref{F1}. Furthermore, we have the $W^{1,2}$-norm convergence
\begin{equation*}
(x^{k}(\cdot),y^{k}(\cdot),u^{k}(\cdot),a^{k}(\cdot),b^{k}(\cdot))\to (\bar{x}(\cdot),\bar{y}(\cdot),\bar{u}(\cdot),\bar{a}(\cdot),\bar{b}(\cdot))\;\mbox{ as }\;k\to\infty.
\end{equation*}
\end{theorem}\vspace*{-0.1in}
{\bf Proof}. Let $\omega^{k}(\cdot)=(\omega^{k}_{1}(\cdot),\omega^{k}_{2}(\cdot))$ be piecewise linear functions on $[0,T] $ such that
\begin{equation*}
(\omega^{k}_{1}(t^{k}_{j}),\omega^{k}_{2}(t^{k}_{j})):=(\bar{a}(t^{k}_{j}),\bar{b}(t^{k}_{j})),\quad j=0,\ldots,k,
\end{equation*}
and let $(\vartheta^{k}_{1}(t),\vartheta^{k}_{2}(t)):=\dot{\omega}_{k}(t)$ be their velocities, which are piecewise constant and right continuous on $[0,T]$. It follows from the above constructions
that we have the convergence
\begin{equation*}
\omega^{k}(\cdot)\to(\bar{a}(\cdot),\bar{b}(\cdot))\;\mbox{ uniformly on }\;[0,T],\;\mbox{ and }\;\vartheta^{k}(\cdot)\to(\dot{\bar{a}}(\cdot),\dot{\bar{b}}(\cdot))\;\mbox{ by norm in }\;L^{2}([0,T],\mathbb{R}^{m+d})
\end{equation*}
as $k\to\infty$. Denoting $(a^{k}(t),b^{k}(t)):=\omega^{k}(t)$ for all $t\in [0,T]$, we now intend to construct sequences of piecewise linear functions $x^k(\cdot)$, $y^k(\cdot)$, and $u^k(\cdot)$ on $[0,T]$ satisfying the conclusions of the theorem. It is clearly sufficient to construct these functions at the mesh points and then extend them piecewise linearly to the whole interval $[0,T]$. For simplicity we use the notation $t_j:=t^k_j$ for the mesh points as $j=0,\ldots,k$ for each fixed $k$. Using the recurrent procedure, suppose that $x^k(t_j)$ is given and then define
\begin{equation}\label{e:9}
u^{k}(t_{j}):=x^{k}(t_{j})-\bar{x}(t_{j})+\bar{u}(t_{j}),\,\,\,j=0,\ldots,k.
\end{equation}
Remembering that the sets $F_1(z)$ in \eqref{F1} are closed and convex, we take the unique projections
\begin{equation}\label{e:10}
-\upsilon^{k}_{j}:=\Pi_{F_{1jk}}(-\dot{\ox}(t_j))\;\mbox{ where }\;F_{1jk}:=F_1(x^k(t_{j}),y^{k}(t_{j+1}),u^{k}(t_{j}),a^{k}(t_{j})),
\end{equation}
and deduce from \eqref{e:9} with taking into account the constructions of $a^{k}(\cdot)$ and $b^{k}(\cdot)$ that
\begin{equation}\label{e:13}
\begin{aligned}
&F_{1}(x^{k}(t_{j}),y^{k}(t_{j+1}),u^{k}(t_{j}),a^{k}(t_{j}))=N_{C(u^{k}(t_{j}))}(x^{k}(t_{j}))+f_{1}(a^{k}(t_{j}),x^{k}(t_{j}))+y^{k}(t_{j+1}) \notag\\
&=F_{1}(\bar{x}(t_{j}),\bar{y}(t_{j}),\bar{u}(t_{j}),\bar{a}(t_{j}))+f_{1}(\bar{a}(t_{j}),x^{k}(t_{j}))-f_{1}(\bar{a}(t_{j}),\bar{x}(t_{j}))+y^{k}(t_{j+1})-\bar{y}(t_{j}).
\end{aligned}
\end{equation}
For all $j=0,\ldots,k-1$, define the vectors and functions
\begin{equation}\label{e:11}
w^{k}_{j}:=f_{2}(\bar{b}(t_{j}),x^{k}(t_{j})),\quad x^{k}(t):=x^{k}(t_{j})+(t-t_{j})v^{k}_{j},\quad y^{k}(t)=y^{k}(t_{j})+(t-t_{j})w^{k}_{j},\,\,\,t\in [t_{j},t_{j+1}],
\end{equation}
and then use below the notation
\begin{equation*}
f_{j}^{x}(s):=\lVert \dot{\bar{x}}(t_{j})-\dot{\bar{x}}(s) \rVert\;\mbox{ and }\;f_{j}^{u}(s)=\lVert \dot{\bar{u}}(t_{j})-\dot{\bar{u}}(s) \rVert,\quad s\in [t_{j},t_{j+1}).
\end{equation*}
Select $s_{j}^{x}$ and $s_{j}^{u} $ from the subintervals $[t_{j},t_{j+1})$ so that
\begin{equation}\label{e:12}
\sup\limits_{s\in[t_{j},t_{j+1}]}f_{j}^{x}(s)\leq \lVert \dot{\bar{x}}(t_{j})-\dot{\bar{x}}(s^{x}_{j}) \rVert+2^{-k}\;\mbox{ and }\;\sup\limits_{s\in[t_{j},t_{j+1}]}f_{j}^{u}(s)\leq\lVert \dot{\bar{u}}(t_{j})-\dot{\bar{u}}(s^{u}_{j})\rVert+2^{-k}.
\end{equation}
It clearly follows from \eqref{e:10}, \eqref{e:13}, and \eqref{e:12} that
\begin{equation}\label{e:15}
\begin{aligned}
\lVert\upsilon^{k}_{j}-\dot{\bar{x}}(s) \rVert&\leq \lVert \upsilon^{k}_{j}-\dot{\bar{x}}(t_{j}) \rVert + \lVert \dot{\bar{x}}(t_{j})-\dot{\bar{x}}(s) \rVert={\rm dist}(-\dot{\bar{x}}(t_{j});F_{1}(x^{k}(t_{j}),y^{k}(t_{j}),u^{k}(t_{j}),a^{k}(t_{j})))+ f_{j}^{x}(s)\notag \\
&\leq \lVert f_{1}(\bar{a}(t_{j}),x^{k}(t_{j}))-f_{1}(\bar{a}(t_{j}),\bar{x}(t_{j}))\rVert+\lVert y^{k}(t_{j+1})-\bar{y}(t_{j})\rVert+\lVert \dot{\bar{x}}(t_{j})-\dot{\bar{x}}(s^{x}_{j}) \rVert+2^{-k}.
\end{aligned}
\end{equation}
Employing the above constructions leads us to the estimates
\begin{equation*}
\begin{aligned}
\lVert y^{k}(t_{j+1})-\bar{y}(t_{j+1}) \rVert& =\Big\|y^{k}(t_{j})+h_{j}f_{2}(\bar{b}(t_{j}),x^{k}(t_{j}))-\bar{y}(t_{j})-\int\limits_{t_{j}}^{t_{j+1}}\dot{\bar{y}}(s)\,ds\Big\|\\
&\leq \lVert y^{k}(t_{j})-\bar{y}(t_{j}) \rVert + \int\limits_{t_{j}}^{t_{j+1}}\lVert f_{2}(\bar{b}(t_{j}),x^{k}(t_{j}))-f_{2}(\bar{b}(s),\bar{x}(s)) \rVert\,ds\\
&\leq \lVert y^{k}(t_{j})-\bar{y}(t_{j}) \rVert + \int\limits_{t_{j}}^{t_{j+1}}\lVert f_{2}(\bar{b}(t_{j}),x^{k}(t_{j}))-f_{2}(\bar{b}(t_{j}),\bar{x}(t_{j})) \rVert\,ds\\
&+\int\limits_{t_{j}}^{t_{j+1}}\lVert f_{2}(\bar{b}(t_{j}),\bar{x}(t_{j}))-f_{2}(\bar{b}(t_{j}),\bar{x}(s))   \rVert\,ds +\int\limits_{t_{j}}^{t_{j+1}}\lVert f_{2}(\bar{b}(t_{j}),\bar{x}(s))-f_{2}(\bar{b}(s),\bar{x}(s))   \rVert\,ds.\\
\end{aligned}
\end{equation*}
The imposed assumptions $(\mathcal{H}_{2,2})$ readily implies that
\begin{eqnarray}\label{e:14}
\begin{array}{ll}
\lVert y^{k}(t_{j+1})-\bar{y}(t_{j+1}) \rVert& \leq \lVert y^{k}(t_{j})-\bar{y}(t_{j}) \rVert +L_{1}h_{j}\lVert x^{k}(t_{j})-\bar{x}(t_{j}) \rVert+ L_{1}\displaystyle\int\limits_{t_{j}}^{t_{j+1}}\lVert \bar{x}(t_{j})-\bar{x}(s) \rVert\,d s\\
&+ L_{2}\displaystyle\int\limits_{t_{j}}^{t_{j+1}}\lVert \bar{b}(t_{j})-\bar{b}(s) \rVert\,d s.
\end{array}
\end{eqnarray}
Since $(\bar{x}(\cdot),\bar{b}(\cdot))\in W^{1,2}([0,T];\R^{n+d})$, for all $ t_{j}\leq s\leq t_{j+1}$ we get
\begin{equation*}
\lVert\bar{x}(s)-\bar{x}(t_{j}) \rVert \leq \int\limits_{t_{j}}^{s}\lVert \dot{\bar{x}}(\tau) \rVert\,d\tau\leq \sqrt{s-t_{j}}\Big(\int\limits_{t_{j}}^{s}\lVert \dot{\bar{x}}(\tau) \rVert^{2}\,d\tau\Big)^{1/2}\leq\sqrt{h_{j}}\lVert \dot{\bar{x}} \rVert_{L^{2}([0,T],\mathbb{R}^{n})},
\end{equation*}
\begin{equation*}
\lVert\bar{b}(s)-\bar{b}(t_{j}) \rVert \leq \int\limits_{t_{j}}^{s}\lVert \dot{\bar{b}}(\tau) \rVert\,d\tau\leq \sqrt{s-t_{j}}\Big(\int\limits_{t_{j}}^{s}\lVert \dot{\bar{b}}(\tau) \rVert^{2}\,d\tau\Big)^{\dfrac{1}{2}}\leq \sqrt{h_{j}}\lVert \dot{\bar{b}} \rVert_{L^{2}([0,T],\mathbb{R}^{d})} .
\end{equation*}
Substituting the latter estimates into \eqref{e:14} gives us
\begin{eqnarray*}
\begin{array}{ll}
\lVert y^{k}(t_{j+1})-\bar{y}(t_{j+1}) \rVert \leq \lVert y^{k}(t_{j})-\bar{y}(t_{j}) \rVert +L_{1}h_{j}\lVert x^{k}(t_{j})-\bar{x}(t_{j}) \rVert+L_{1}h_{j}\sqrt{h_{j}}\lVert \dot{\bar{x}} \rVert_{L^{2}([0,T],\mathbb{R}^{n})}+L_{2}h_{j}\sqrt{h_{j}}\lVert \dot{\bar{b}} \rVert_{L^{2}([0,T],\mathbb{R}^{d})}.
\end{array}
\end{eqnarray*}
Proceeding further by induction, we obtain for all $j=0,\ldots,k-1$ that
\begin{eqnarray*}
\begin{array}{ll}
\lVert y^{k}(t_{j+1})-\bar{y}(t_{j+1}) \rVert\leq L_{1}\sum\limits_{i=0}^{j}h_{i}\lVert x^{k}(t_{i})-\bar{x}(t_{i}) \rVert+T\sqrt{h_{k}}(L_{1}\lVert \dot{\bar{x}} \rVert_{L^{2}([0,T],\mathbb{R}^{n})}+L_{2}\lVert \dot{\bar{b}} \rVert_{L^{2}([0,T],\mathbb{R}^{d})}),
\end{array}
\end{eqnarray*}
which implies in turn by using the imposed assumptions on $f_{2}$ that
\begin{equation*}
\begin{aligned}
\lVert y^{k}(t_{j+1})-\bar{y}(t_{j}) \rVert&\leq \lVert y^{k}(t_{j+1})-\bar{y}(t_{j+1}) \rVert + \lVert \bar{y}(t_{j+1})-\bar{y}(t_{j}) \rVert\\
&\leq L_{1}\sum\limits_{i=0}^{j}h_{i}\lVert x^{k}(t_{i})-\bar{x}(t_{i}) \rVert+T\sqrt{h_{k}}(L_{1}\lVert \dot{\bar{x}} \rVert_{L^{2}([0,T],\mathbb{R}^{n})}+L_{2}\lVert \dot{\bar{b}} \rVert_{L^{2}([0,T],\mathbb{R}^{d})} )+\int\limits_{t_{j}}^{t_{j+1}}\lVert \dot{\bar{y}}(s) \rVert\,ds\\
&=L_{1}\sum\limits_{i=0}^{j}h_{i}\lVert x^{k}(t_{i})-\bar{x}(t_{i}) \rVert+T\sqrt{h_{k}}(L_{1}\lVert \dot{\bar{x}} \rVert_{L^{2}([0,T],\mathbb{R}^{n})}+L_{2}\lVert \dot{\bar{b}} \rVert_{L^{2}([0,T],\mathbb{R}^{d})} ) +\int\limits_{t_{j}}^{t_{j+1}}\lVert f_{2}(\bar{b}(s),\bar{x}(s)) \rVert\,ds\\
&\leq L_{1}\sum\limits_{i=0}^{j}h_{i}\lVert x^{k}(t_{i})-\bar{x}(t_{i}) \rVert+T\sqrt{h_{k}}(L_{1}\lVert \dot{\bar{x}} \rVert_{L^{2}([0,T],\mathbb{R}^{n})}+L_{2}\lVert \dot{\bar{b}} \rVert_{L^{2}([0,T],\mathbb{R}^{d})} ) \\
&+\alpha_{2}\sqrt{h_{k}}(\lVert\bar{b}(s)\rVert_{L^{2}([0,T],\mathbb{R}^{d})}+\lVert\bar{x}(s) \rVert_{L^{2}([0,T],\mathbb{R}^{n})}) .
\end{aligned}
\end{equation*}
Substituting the latter into the previous estimate of $\lVert y^{k}(t_{j+1})-\bar{y}(t_{j+1})\rVert $ and employing $(\mathcal{H}_{2,2}) $ bring us to
\begin{eqnarray*}
\begin{aligned}
\lVert \upsilon^{k}_{j}-\dot{\bar{x}}(s) \rVert
&\leq L\lVert x^{k}(t_{j})-\bar{x}(t_{j})\rVert+L_{1}\sum\limits_{i=0}^{j}h_{i}\lVert x^{k}(t_{i})-\bar{x}(t_{i}) \rVert+T\sqrt{h_{k}}(L_{1}\lVert \dot{\bar{x}} \rVert_{L^{2}([0,T],\mathbb{R}^{n})}+L_{2}\lVert \dot{\bar{b}} \rVert_{L^{2}([0,T],\mathbb{R}^{n})} )\\
&+\alpha_{2}\sqrt{h_{k}}(\lVert\bar{b}(s)\rVert_{L^{2}([0,T],\mathbb{R}^{d})}+\lVert\bar{x}(s) \rVert_{L^{2}([0,T],\mathbb{R}^{n})}) +\lVert \dot{\bar{x}}(t_{j})-\dot{\bar{x}}(s^{x}_{j}) \rVert+2^{-k}.
\end{aligned}
\end{eqnarray*}
We get from the above relationships the following estimates:
\begin{equation*}
\begin{aligned}
\lVert x^{k}(t_{j+1})-\bar{x}(t_{j+1}) \rVert&=\lVert  x^{k}(t_{j})+h_{k}v^{k}_{j}-\bar{x}(t_{j})-\int\limits_{t_{j}}^{t_{j+1}}\dot{\bar{x}}(s)\,ds \rVert\leq\lVert x^{k}(t_{j})-\bar{x}(t_{j})
\rVert + \int\limits_{t_{j}}^{t_{j+1}}\lVert v^{k}_{j}-\dot{\bar{x}}(s)\rVert\,ds \notag\\
&\leq\lVert x^{k}(t_{j})-\bar{x}(t_{j}) \rVert + Lh_{k}\lVert x^{k}(t_{j})-\bar{x}(t_{j})\rVert+L_{1}h_{k}\sum\limits_{i=0}^{j}h_{i}\lVert x^{k}(t_{i})-\bar{x}(t_{i}) \rVert\notag\\
&+Th_{k}\sqrt{h_{k}}(L_{1}\lVert \dot{\bar{x}} \rVert_{L^{2}([0,T],\mathbb{R}^{n})}+L_{2}\lVert \dot{\bar{b}} \rVert_{L^{2}([0,T],\mathbb{R}^{d})} )\\
&+\alpha_{2}h_{k}\sqrt{h_{k}}(\lVert\bar{b}(s)\rVert_{L^{2}([0,T],\mathbb{R}^{d})}+\lVert\bar{x}(s) \rVert_{L^{2}([0,T],\mathbb{R}^{n})}) +h_{k}\lVert \dot{\bar{x}}(t_{j})-\dot{\bar{x}}(s^{x}_{j}) \rVert+h_{k}2^{-k} \notag\\
&\leq (1+Lh_{k}+L_{1}h_{k}^{2}) \lVert x^{k}(t_{j})-\bar{x}(t_{j}) \rVert  +L_{1}h_{k}^{2}\sum\limits_{i=0}^{j-1}\lVert x^{k}(t_{i})-\bar{x}(t_{i}) \rVert\notag\\
&+Th_{k}\sqrt{h_{k}}(L_{1}\lVert \dot{\bar{x}} \rVert_{L^{2}([0,T],\mathbb{R}^{n})}+L_{2}\lVert \dot{\bar{b}} \rVert_{L^{2}([0,T],\mathbb{R}^{d})} )\\
&+\alpha_{2}h_{k}\sqrt{h_{k}}(\lVert\bar{b}(s)\rVert_{L^{2}([0,T],\mathbb{R}^{d})}+\lVert\bar{x}(s) \rVert_{L^{2}([0,T],\mathbb{R}^{n})})+h_{k}\lVert \dot{\bar{x}}(t_{j})-\dot{\bar{x}}(s^{x}_{j}) \rVert+h_{k}2^{-k} .\label{1}
\end{aligned}
\end{equation*}
Apply now the discrete Gronwall's lemma from Proposition~\ref{granwal} with the given parameters
\begin{equation*}
e_{j}=\lVert x^{k}(t_{j})-\bar{x}(t_{j}) \rVert,\,\,\,\rho_{j}=L_{1}h_{k}^{2},\,\,\,\gamma_{j}=Lh_{k}+L_{1}h_{k}^{2},\;\mbox{ and }\;
\end{equation*}
\begin{equation*}
\sigma_{j}=Th_{k}\sqrt{h_{k}}(L_{1}\lVert \dot{\bar{x}} \rVert_{L^{2}([0,T],\mathbb{R}^{n})}+L_{2}\lVert \dot{\bar{b}} \rVert_{L^{2}([0,T],\mathbb{R}^{d})} )
+h_{k}\lVert \dot{\bar{x}}(t_{j})-\dot{\bar{x}}(s^{x}_{j}) \rVert+h_{k}2^{-k}
\end{equation*}
\begin{equation*}
+\alpha_{2}h_{k}\sqrt{h_{k}}(\lVert\bar{b}(s)\rVert_{L^{2}([0,T],\mathbb{R}^{d})}+\lVert\bar{x}(s) \rVert_{L^{2}([0,T],\mathbb{R}^{n})}).
\end{equation*}
This ensures the fulfilment of the relationships
\begin{equation*}
\begin{aligned}
\sum\limits_{i=0}^{j-1}\sigma_{i}&=jh_{k}\Big(T\sqrt{h_{k}}(L_{1}\lVert \dot{\bar{x}} \rVert_{L^{2}([0,T],\mathbb{R}^{n})}+L_{2}\lVert \dot{\bar{b}} \rVert_{L^{2}([0,T],\mathbb{R}^{d})} )+2^{-k}\Big)+h_{k}\sum\limits_{i=0}^{j-1}\lVert \dot{\bar{x}}(t_{j})-\dot{\bar{x}}(s^{x}_{j}) \rVert\\
&+j\alpha_{2}h_{k}\sqrt{h_{k}}(\lVert\bar{b}(s)\rVert_{L^{2}([0,T],\mathbb{R}^{d})}+\lVert\bar{x}(s) \rVert_{L^{2}([0,T],\mathbb{R}^{n})})\\
&\leq T^{2}\sqrt{h_{k}}(L_{1}\lVert \dot{\bar{x}} \rVert_{L^{2}([0,T],\mathbb{R}^{n})}+L_{2}\lVert \dot{\bar{b}} \rVert_{L^{2}([0,T],\mathbb{R}^{d})} )+T2^{-k}\\
&+h_{k}\sum\limits_{i=0}^{k-1}(\lVert \dot{\bar{x}}(t_{j})-\dot{\bar{x}}(s^{x}_{j}) \rVert+\lVert \dot{\bar{x}}(s^{x}_{j})-\dot{\bar{x}}(t_{j+1}) \rVert)+T\alpha_{2}\sqrt{h_{k}}(\lVert\bar{b}(s)\rVert_{L^{2}([0,T],\mathbb{R}^{d})}+\lVert\bar{x}(s) \rVert_{L^{2}([0,T],\mathbb{R}^{n})})\\
&\leq T^{2}\sqrt{h_{k}}(L_{1}\lVert \dot{\bar{x}} \rVert_{L^{2}([0,T],\mathbb{R}^{n})}+L_{2}\lVert \dot{\bar{b}} \rVert_{L^{2}([0,T],\mathbb{R}^{d})} )+T2^{-k}+h_{k}\mu+T\alpha_{2}\sqrt{h_{k}}(\lVert\bar{b}(s)\rVert_{L^{2}([0,T],\mathbb{R}^{d})}+\lVert\bar{x}(s) \rVert_{L^{2}([0,T],\mathbb{R}^{n})}).
\end{aligned}
\end{equation*}
Having further $\lim\limits_{k\to \infty}c_{k}:=\lim\limits_{k\to \infty} T^{2}\sqrt{h_{k}}(L_{1}\lVert \dot{\bar{x}} \rVert_{L^{2}([0,T],\mathbb{R}^{n})}+L_{2}\lVert \dot{\bar{b}} \rVert_{L^{2}([0,T],\mathbb{R}^{d})} )+T2^{-k}+h_{k}\mu+T\alpha_{2}\sqrt{h_{k}}(\lVert\bar{b}(s)\rVert_{L^{2}([0,T],\mathbb{R}^{d})}+\lVert\bar{x}(s) \rVert_{L^{2}([0,T],\mathbb{R}^{n})})=0$, gives us the conditions
\begin{equation*}
\begin{aligned}
\sum\limits_{i=0}^{j-1}(k\rho_{k}+\gamma_{k})=L_{1}\sum\limits_{i=0}^{j-1}(ih_{k}^{2}+h_{k}+h_{k}^{2})=L_{1}\Big(h_{k}^{2}\dfrac{j(j-1)}{2}+j(h_{k}+h_{k}^{2})\Big)\leq L_{1}T\Big(\dfrac{T}{2}+1+h_{k}\Big),
\end{aligned}
\end{equation*}
which imply that for all $j=0,\ldots,k-1$ we get
\begin{equation}\label{e:16}
\lVert x^{k}(t_{j})-\bar{x}(t_{j}) \rVert\leq c_{k}\exp\Big(L_{1}T(\dfrac{T}{2}+1+h_{k})\Big ),\;\mbox{ and hence}
\end{equation}
\begin{equation}\label{e:17}
\lVert y^{k}(t_{j})-\bar{y}(t_{j}) \rVert\leq L_{1}Tc_{k}\exp\Big(L_{1}T(\dfrac{T}{2}+1+h_{k})\Big )+T\sqrt{h_{k}}(L_{1}\lVert \dot{\bar{x}} \rVert_{L^{2}([0,T],\mathbb{R}^{n})}+L_{2}\lVert \dot{\bar{b}} \rVert_{L^{2}([0,T],\mathbb{R}^{d})} ),
\end{equation}
\begin{equation}\label{e:18}
\begin{aligned}
\lVert \upsilon^{k}_{j}-\dot{\bar{x}}(s) \rVert &\leq c_{k}\exp\Big(L_{1}T(\dfrac{T}{2}+1+h_{k})\Big )(L+L_{1}T)+T\sqrt{h_{k}}(L_{1}\lVert \dot{\bar{x}} \rVert_{L^{2}([0,T],\mathbb{R}^{n})}+L_{2}\lVert \dot{\bar{b}} \rVert_{L^{2}([0,T],\mathbb{R}^{d})} )\\
&+\alpha_{2}\sqrt{h_{k}}(\lVert\bar{b}(s)\rVert_{L^{2}([0,T],\mathbb{R}^{d})}+\lVert\bar{x}(s) \rVert_{L^{2}([0,T],\mathbb{R}^{n})})+ \lVert \dot{\bar{x}}(t_{j})-\dot{\bar{x}}(s^{x}_{j}) \rVert+2^{-k}.
\end{aligned}
\end{equation}
Employing the obtained conditions together with \eqref{e:16} and \eqref{e:17} tells us that
\begin{equation*}
\begin{aligned}
\lVert x^{k}(t)-\bar{x}(t) \rVert&=\lVert x^{k}(t_{j})+(t-t_{j})\upsilon^{k}_{j}-\bar{x}(t_{j})-\int\limits_{t_{j}}^{t}\dot{\bar{x}}(s)\,ds \rVert\leq \lVert x^{k}(t_{j})-\bar{x}(t_{j}) \rVert + \int\limits_{t_{j}}^{t_{j+1}}\lVert v^{k}_{j}-\dot{\bar{x}}(s)\rVert\,ds\\
&\leq c_{k}\exp\Big(L_{1}T(\dfrac{T}{2}+1)\Big) + h_{k}(\lVert \dot{\bar{x}}(t_{j})-\dot{\bar{x}}(s^{x}_{j}) \rVert+2^{-k})+h_{k}c_{k}\exp\Big(L_{1}T(\dfrac{T}{2}+1)\Big)(L+L_{1}T)\\
&+Th_{k}\sqrt{h_{k}}(L_{1}\lVert \dot{\bar{x}} \rVert_{L^{2}([0,T],\mathbb{R}^{n})}+L_{2}\lVert \dot{\bar{b}} \rVert_{L^{2}([0,T],\mathbb{R}^{d})} )+h_{k}\alpha_{2}\sqrt{h_{k}}(\lVert\bar{b}(s)\rVert_{L^{2}([0,T],\mathbb{R}^{d})}+\lVert\bar{x}(s) \rVert_{L^{2}([0,T],\mathbb{R}^{n})}),
\end{aligned}
\end{equation*}
which readily verifies the uniform convergence of the sequence $\{x^{k}(\cdot)\}$ to $\bar{x}(\cdot) $ as $k\to \infty $.\vspace*{-0.05in}

To proceed further, deduce from \eqref{e:18} for $j=0,\ldots,k-1$ that
\begin{equation*}
\begin{aligned}
\int\limits_{0}^{T}\lVert \dot{x}^{k}(t)-\dot{\bar{x}}(t) \rVert^{2}\,dt
&=\sum\limits_{j=0}^{k-1}\int\limits_{t_{j}}^{t_{j+1}}\lVert \upsilon^{k}_{j}-\dot{\bar{x}}(t) \rVert^{2}\,dt\\
&\leq 3h_{k}\sum\limits_{j=0}^{k-1}\Big(c_{k}\exp\Big(L_{1}T(\dfrac{T}{2}+1)\Big )(L+L_{1}T)+T\sqrt{h_{k}}(L_{1}\lVert \dot{\bar{x}} \rVert_{L^{2}([0,T],\mathbb{R}^{n})}+L_{2}\lVert \dot{\bar{b}} \rVert_{L^{2}([0,T],\mathbb{R}^{d})} )\Big)^{2}\\
&+3h_{k}\alpha_{2}\sqrt{h_{k}}\sum\limits_{j=0}^{k-1}\Big(\lVert\bar{b}(s)\rVert_{L^{2}([0,T],\mathbb{R}^{d})}+\lVert\bar{x}(s) \rVert_{L^{2}([0,T],\mathbb{R}^{n})}\Big)^{2}+3h_{k}\sum\limits_{j=0}^{k-1}(\lVert \dot{\bar{x}}(t_{j})-\dot{\bar{x}}(s^{x}_{j}) \rVert+2^{-k})^{2} .
\end{aligned}
\end{equation*}
Observe in addition from the constructions and assumptions above that
\begin{equation*}
\begin{aligned}
&3h_{k}\sum\limits_{j=0}^{k-1}(\lVert \dot{\bar{x}}(t_{j})-\dot{\bar{x}}(s^{x}_{j}) \rVert+2^{-k})^{2}\leq 6h_{k}\sum\limits_{j=0}^{k-1}(\lVert \dot{\bar{x}}(t_{j})-\dot{\bar{x}}(s^{x}_{j}) \rVert^{2}+4^{-k})\\
&\leq 6h_{k}\Big(\sum\limits_{j=0}^{k-1}\lVert \dot{\bar{x}}(t_{j})-\dot{\bar{x}}(s^{x}_{j}) \rVert\Big)^{2}+6T4^{-k}\leq 6h_{k}\Big(\sum\limits_{j=0}^{k-1}\lVert \dot{\bar{x}}(t_{j})-\dot{\bar{x}}(s^{x}_{j}) \rVert+\lVert \dot{\bar{x}}(s^{x}_{j})-\dot{\bar{x}}(t_{j+1}) \rVert\Big)^{2}+6T4^{-k}\leq 6h_{k}\mu^{2}+6T4^{-k},
\end{aligned}
\end{equation*}
which yields in turn the estimates
\begin{equation*}
\begin{aligned}
\int\limits_{0}^{T}\lVert \dot{x}^{k}(t)-\dot{\bar{x}}(t) \rVert^{2}\,dt
&\leq 3T\Big(c_{k}\exp\Big(L_{1}T(\dfrac{T}{2}+1)\Big )(L+L_{1}T)+T\sqrt{h_{k}}(L_{1}\lVert \dot{\bar{x}} \rVert_{L^{2}([0,T],\mathbb{R}^{n})}+L_{2}\lVert \dot{\bar{b}} \rVert_{L^{2}([0,T],\mathbb{R}^{d})} )\Big)^{2}\\
&+3T\alpha_{2}\sqrt{h_{k}}\Big(\lVert\bar{b}(s)\rVert_{L^{2}([0,T],\mathbb{R}^{d})}+\lVert\bar{x}(s) \rVert_{L^{2}([0,T],\mathbb{R}^{n})}\Big)^{2}+6h_{k}\mu^{2}+6T4^{-k}.
\end{aligned}
\end{equation*}
Remembering that $ c_{k}\to 0$, this verifies the strong convergence of the sequence $\{\dot{x}^{k}(\cdot)\}$ to $ \dot{\bar{x}}(\cdot) $ in the norm topology of $L^{2}([0,T],\mathbb{R}^{n}) $ as $k\to\infty$.

Our next step is to justify the $W^{1,2}$-strong convergence of the sequence $\{y^{k}(\cdot)\}$ to $\oy(\cdot)$. As follows from the above construction, for all $ j=0,\ldots,k-1 $ and
$t\in[t_{j},t_{j+1}]$ we have
\begin{equation*}
\begin{aligned}
\lVert \dot{y}^{k}(t)-\dot{\bar{y}}(t) \rVert&=\lVert w^{k}_{j}-\dot{\bar{y}}(t) \rVert=\lVert f_{2}(\bar{b}(t_{j}),x^{k}(t_{j}))-f_{2}(\bar{b}(t),\bar{x}(t)) \rVert\leq \lVert f_{2}(\bar{b}(t_{j}),x^{k}(t_{j}))-f_{2}(\bar{b}(t_{j}),\bar{x}(t_{j})) \rVert\\
&+\lVert f_{2}(\bar{b}(t_{j}),\bar{x}(t_{j}))-f_{2}(\bar{b}(t_{j}),\bar{x}(t)) \rVert+\lVert f_{2}(\bar{b}(t_{j}),\bar{x}(t))-f_{2}(\bar{b}(t),\bar{x}(t)) \rVert\\
&\leq L_{1}c_{k}\exp\Big(L_{1}T(\dfrac{T}{2}+1+h_{k})\Big)+\sqrt{h_{k}}(L_{1}\lVert \dot{\bar{x}} \rVert_{L^{2}([0,T],\mathbb{R}^{n})}+L_{2}\lVert \dot{\bar{b}} \rVert_{L^{2}([0,T],\mathbb{R}^{n})} ).
\end{aligned}
\end{equation*}
Combining the latter with \eqref{e:17} gives us the estimates
\begin{equation*}
\begin{aligned}
\lVert y^{k}(t)-\bar{y}(t) \rVert&=\lVert y^{k}(t_{j})+(t-t_{j})w^{k}_{j}-\bar{y}(t_{j})-\int\limits_{t_{j}}^{t}\dot{\bar{y}}(s)\,ds \rVert\leq \lVert y^{k}(t_{j})-\bar{y}(t_{j}) \rVert + \int\limits_{t_{j}}^{t_{j+1}}\lVert w^{k}_{j}-\dot{\bar{y}}(s)\rVert\,ds\\
&\leq L_{1}Tc_{k}\exp\Big(L_{1}T(\dfrac{T}{2}+1+h_{k})\Big)+T\sqrt{h_{k}}(L_{1}\lVert \dot{\bar{x}} \rVert_{L^{2}([0,T],\mathbb{R}^{n})}+L_{2}\lVert \dot{\bar{b}} \rVert_{L^{2}([0,T],\mathbb{R}^{n})} ) \\
&+L_{1}h_{k}c_{k}\exp\Big(L_{1}T(\dfrac{T}{2}+1+h_{k})\Big)+h_{k}\sqrt{h_{k}}(L_{1}\lVert \dot{\bar{x}} \rVert_{L^{2}([0,T],\mathbb{R}^{n})}+L_{2}\lVert \dot{\bar{b}} \rVert_{L^{2}([0,T],\mathbb{R}^{n})} ),
\end{aligned}
\end{equation*}
\begin{equation*}
\begin{aligned}
\int\limits_{0}^{T}\lVert \dot{y}^{k}(t)-\dot{\bar{y}}(t) \rVert^{2}\,dt&=\sum\limits_{j=0}^{k-1}\int\limits_{t_{j}}^{t_{j+1}}\lVert \dot{y}^{k}(t)-\dot{\bar{y}}(t) \rVert^{2}\,dt\\
&\leq T \Big( L_{1}c_{k}\exp\Big(L_{1}T(\dfrac{T}{2}+1+h_{k})\Big)+\sqrt{h_{k}}\Big(L_{1}\lVert \dot{\bar{x}} \rVert_{L^{2}([0,T],\mathbb{R}^{n})}+L_{2}\lVert \dot{\bar{b}} \rVert_{L^{2}([0,T],\mathbb{R}^{n})}\Big)\Big)^{2},
\end{aligned}
\end{equation*}
which therefore ensure the convergence of $ y^{k}(\cdot) $ to $ \bar{y}(\cdot) $ in $ W^{1,2}([0,T],\mathbb{R}^{n})$.\vspace*{-0.05in}

To complete the proof of the theorem, it remains to verify the $ W^{1,2}$-strong convergence of the control sequence $\{u^{k}(\cdot)\}$ to $\ou(\cdot)$. To this end, observe first that
\begin{equation}\label{e:19}
\int\limits_{0}^{T}\lVert \dot{u}^{k}(t)-\dot{\bar{u}}(t) \rVert^{2}\,dt\leq 2\int\limits_{0}^{T}\Big\lVert \dot{u}^{k}(t)-\dfrac{\bar{u}(t_{j+1})-\bar{u}(t_{j})}{h_{k}} \Big\rVert^{2}\,dt+2\int\limits_{0}^{T}\Big\lVert \dfrac{\bar{u}(t_{j+1})-\bar{u}(t_{j})}{h_{k}}-\dot{\bar{u}}(t)  \Big\rVert^{2}\,dt
\end{equation}
for all $j=0,\ldots,k-1$ and $ t\in[t_{j},t_{j+1})$. On the other hand, it follows from \eqref{e:9} that
\begin{equation*}
\begin{aligned}
&\int\limits_{0}^{T}\Big\lVert \dot{u}^{k}(t)-\dfrac{\bar{u}(t_{j+1})-\bar{u}(t_{j})}{h_{k}} \Big\rVert^{2}\,dt=\int\limits_{0}^{T}\Big\lVert \dfrac{u^{k}(t_{j+1})-u^{k}(t_{j})}{h_{k}}-\dfrac{\bar{u}(t_{j+1})-\bar{u}(t_{j})}{h_{k}} \Big\rVert^{2}\,dt\\
&=\int\limits_{0}^{T}\Big\lVert \dfrac{u^{k}(t_{j+1})-\bar{u}(t_{j+1})}{h_{k}}-\dfrac{u^{k}(t_{j})-\bar{u}(t_{j})}{h_{k}} \Big\rVert^{2}\,dt=\int\limits_{0}^{T}\Big\lVert \dfrac{x^{k}(t_{j+1})-\bar{x}(t_{j+1})}{h_{k}}-\dfrac{x^{k}(t_{j})-\bar{x}(t_{j})}{h_{k}} \Big\rVert^{2}\,dt\\
&=\int\limits_{0}^{T}\Big\lVert \dfrac{x^{k}(t_{j+1})-x^{k}(t_{j})}{h_{k}}-\dfrac{\bar{x}(t_{j+1})-\bar{x}(t_{j})}{h_{k}} \Big\rVert^{2}\,dt
\leq 2\int\limits_{0}^{T}\Big\lVert \dfrac{x^{k}(t_{j+1})-x^{k}(t_{j})}{h_{k}}-\dot{\bar{x}}(t)\Big\rVert^{2}\,dt+2\int\limits_{0}^{T}\Big\lVert \dot{\bar{x}}(t)-\dfrac{\bar{x}(t_{j+1})-\bar{x}(t_{j})}{h_{k}} \Big\rVert^{2}\,dt\\
&=2\int\limits_{0}^{T}\lVert \dot{x}^{k}(t)-\dot{\bar{x}}(t) \rVert^{2}\,dt+2\int\limits_{0}^{T}\Big\lVert \dot{\bar{x}}(t)-\dfrac{\bar{x}(t_{j+1})-\bar{x}(t_{j})}{h_{k}} \Big\rVert^{2}\,dt\quad\mbox{and}
\end{aligned}
\end{equation*}
\begin{equation*}
\begin{aligned}
&\Big\lVert \dot{\bar{x}}(t)-\dfrac{\bar{x}(t_{j+1})-\bar{x}(t_{j})}{h_{k}} \Big\rVert\leq\lVert \dot{\bar{x}}(t)-\dot{\bar{x}}(t_{j})\rVert+\Big\lVert \dot{\bar{x}}(t_{j})-\dfrac{\bar{x}(t_{j+1})-\bar{x}(t_{j})}{h_{k}} \Big\rVert\\
&\leq \lVert \dot{\bar{x}}(t)-\dot{\bar{x}}(t_{j})\rVert+\dfrac{1}{h_{k}} \int\limits_{t_{j}}^{t_{j+1}}\lVert\dot{\bar{x}}(t_{j})-\dot{\bar{x}}(t)\rVert\,dt \leq\lVert\dot{\bar{x}}(t_{j})-\dot{\bar{x}}(s_{j}^{x})\rVert+\dfrac{1}{h_{k}} \int\limits_{t_{j}}^{t_{j+1}}\lVert\dot{\bar{x}}(t_{j})-\dot{\bar{x}}(s_{j}^{x})\rVert\,dt+2^{-k}\\
&\leq \lVert\dot{\bar{x}}(t_{j})-\dot{\bar{x}}(s_{j}^{x})\rVert+\lVert\dot{\bar{x}}(s_{j}^{x})-\dot{\bar{x}}(t_{j+1})\rVert+2^{-k}.
\end{aligned}
\end{equation*}
Hence we arrive at the relationships
\begin{equation*}
\begin{aligned}
&\int\limits_{0}^{T}\Big\lVert \dot{\bar{x}}(t)-\dfrac{\bar{x}(t_{j+1})-\bar{x}(t_{j})}{h_{k}} \Big\rVert^{2}\,dt\leq\int\limits_{0}^{T}\Big(\lVert\dot{\bar{x}}(t_{j})-\dot{\bar{x}}(s_{j}^{x})\rVert+\lVert\dot{\bar{x}}(s_{j}^{x})-\dot{\bar{x}}(t_{j+1})\rVert+2^{-k}\Big)^{2}\,dt\\
&\leq 2\int\limits_{0}^{T}\Big(\lVert\dot{\bar{x}}(t_{j})-\dot{\bar{x}}(s_{j}^{x})\rVert+\lVert\dot{\bar{x}}(s_{j}^{x})-\dot{\bar{x}}(t_{j+1})\rVert\Big)^{2}\,dt+T4^{-k}=2\sum\limits_{j=0}^{k-1}
\int\limits_{t_{j}}^{t_{j+1}}\Big(\lVert\dot{\bar{x}}(t_{j})-\dot{\bar{x}}(s_{j}^{x})\rVert+\lVert\dot{\bar{x}}(s_{j}^{x})-\dot{\bar{x}}(t_{j+1})\rVert\Big)^{2}\,dt\\
&+T4^{-k+1}\leq 2 h_{k}\Big(\sum\limits_{j=0}^{k-1}\big(\lVert\dot{\bar{x}}(t_{j})-\dot{\bar{x}}(s_{j}^{x})\rVert+\lVert\dot{\bar{x}}(s_{j}^{x})-\dot{\bar{x}}(t_{j+1})\rVert\big)\Big)^{2}+T4^{-k+1}\leq 2 h_{k}\mu^{2}+T4^{-k+1},
\end{aligned}
\end{equation*}
which clearly yield the estimate
\begin{equation}\label{e:20}
\int\limits_{0}^{T}\Big\lVert \dot{u}^{k}(t)-\dfrac{\bar{u}(t_{j+1})-\bar{u}(t_{j})}{h_{k}} \Big\rVert^{2}\,dt\leq 2\int\limits_{0}^{T}\lVert \dot{x}^{k}(t)-\dot{\bar{x}}(t) \rVert^{2}\,dt+8 h_{k}\mu^{2}+2T4^{-k+1}.
\end{equation}
In the same way we obtain the condition
\begin{equation}\label{e:21}
\int\limits_{0}^{T}\Big\lVert \dot{\bar{u}}(t)-\dfrac{\bar{u}(t_{j+1})-\bar{u}(t_{j})}{h_{k}} \Big\rVert^{2}\,dt\leq 2 h_{k}\mu^{2}+T4^{-k+1} .
\end{equation}
Substituting \eqref{e:20} and \eqref{e:21} into \eqref{e:19} gives us
\begin{equation*}
\int\limits_{0}^{T}\lVert \dot{u}^{k}(t)-\dot{\bar{u}}(t) \rVert^{2}\,dt\leq 4\int\limits_{0}^{T}\lVert \dot{x}^{k}(t)-\dot{\bar{x}}(t) \rVert^{2}\,dt+16 h_{k}\mu^{2}+4T4^{-k+1}+4 h_{k}\mu^{2}+2T4^{-k+1},
\end{equation*}
which readily ensures the strong convergence of the sequence $\{\dot u^{k}(\cdot)\}$ to $\dot{\ou}(\cdot)$ in the norm topology of $ L^{2}([0,T],\mathbb{R}^{n})$ with $ u^{k}(0)=\ou(0)+x^{k}(0)-x_0=\ou(0)$. Finally, for all $ t\in [0,T]$ we get
\begin{equation*}
\begin{aligned}
\lVert u^{k}(t)-\bar{u}(t) \rVert\leq \int\limits_{0}^{t}\lVert \dot{u}^{k}(s)-\dot{\bar{u}}(s) \rVert\,ds\leq \sqrt{T}\Big(\int\limits_{0}^{T}\lVert \dot{u}^{k}(t)-\dot{\bar{u}}(t) \rVert^{2}\,ds\Big)^{1/2} \to 0\;\mbox{ as }\;k\to \infty,
\end{aligned}
\end{equation*}
and thus the sequence $\{u^{k}(\cdot)\}$ converges to the feasible control $\bar{u}(\cdot) $ strongly in $  W^{1,2}([0,T],\mathbb{R}^{n})$ as $k\to\infty$. This completes the proof of the theorem. $\h$\vspace*{-0.2in}

\section{Discrete Approximations of Local Optimal Solutions}\label{sec:disc-opt}
\setcounter{equation}{0}\vspace*{-0.1in}

In this section we continue developing the method of discrete approximations for controlled integro-differential sweeping processes while in a different framework. Our goal is to (strongly in $W^{1,2}$) approximate a prescribed relaxed intermediate local minimizer of the optimal control problem $(P)$ by optimal solutions to discretized sweeping control systems. Given an r.i.l.m.\
$\bar{z}(\cdot)=(\bar{x}(\cdot), \bar{y}(\cdot) ,\bar{u}(\cdot),\bar{a}(\cdot),\bar{b}(\cdot))$ of $(P)$ with the number $\ve>0$ from Definition~\ref{locmin} (ii), consider the mesh $\Delta_k$ defined in \eqref{mesh} and construct the sequence of discrete-time sweeping optimal control problems $(P_k)$, $k\in\N$, as follows:
\begin{equation*}
\begin{aligned}
\text{minimize}\,\,\, J_{k}(z^{k})&:=\varphi(x^{k}_k)+h_{k}\sum\limits_{j=0}^{k-1}l\Big(t^{k}_{j},x^{k}_{j},y^{k}_{j},u^{k}_{j},a^{k}_{j},b^{k}_{j},\dfrac{x^{k}_{j+1}-x^{k}_{j}}{h_{k}},\dfrac{y^{k}_{j+1}-y^{k}_{j}}{h_{k}},
\dfrac{u^{k}_{j+1}-u^{k}_{j}}{h_{k}},\dfrac{a^{k}_{j+1}-a^{k}_{j}}{h_{k}},\dfrac{b^{k}_{j+1}-b^{k}_{j}}{h_{k}}\Big)\\
&+\dfrac{1}{2}\sum\limits_{j=0}^{k-1}\int\limits_{t^{k}_{j}}^{t^{k}_{j+1}}\Big(\Big\lVert \dfrac{x^{k}_{j+1}-x^{k}_{j}}{h_{k}}-\dot{\bar{x}}(t) \Big\rVert^{2}
+\Big\lVert \dfrac{y^{k}_{j+1}-y^{k}_{j}}{h_{k}}-\dot{\bar{y}}(t) \Big\rVert^{2}\Big)\,dt\\
&+\dfrac{1}{2}\sum\limits_{j=0}^{k-1}\int\limits_{t^{k}_{j}}^{t^{k}_{j+1}}\Big(\Big\lVert \dfrac{u^{k}_{j+1}-u^{k}_{j}}{h_{k}}-\dot{\bar{u}}(t) \Big\rVert^{2}+\Big\lVert \dfrac{a^{k}_{j+1}-a^{k}_{j}}{h_{k}}-\dot{\bar{a}}(t) \Big\rVert^{2}+\Big\lVert \dfrac{b^{k}_{j+1}-b^{k}_{j}}{h_{k}}-\dot{\bar{b}}(t) \Big\rVert^{2}\Big)\,dt
\end{aligned}
\end{equation*}
over the collections $ z^{k}=(x_{0}^{k},\ldots,x^{k}_k,y_{0}^{k},\ldots,y^{k}_k,u^{k}_{0},\ldots,u^{k}_{k},a_{0}^{k},\ldots,a^{k}_{k},b_{0}^{k},\ldots,b^{k}_{k}) $ subject to the constraints:\\
\begin{equation}\label{e:5.0}
x^{k}_{j+1}\in x^{k}_{j}-h_kF_h(x^k_j,y^k_j,u^k_j,a^k_j,b^k_j),\quad j=0,\ldots,k-1,
\end{equation}
where the discrete velocity mapping $F_h$ is given in the form
\begin{equation}\label{Fh0}
F_h(x^k_j,y^k_j,u^k_j,a^k_j,b^k_j):=F_{1}(x^{k}_{j},y^{k}_{j+1},u_{j}^{k},a^{k}_{j})=N_{C(u^k_j)}(x^k_j)+f_1(x^k_j,a^k_j)+y^k_{j+1}
\end{equation}
due to the definition of $F_1$ in \eqref{F1}, and where $(x^{k}_{0},u^{k}_{0},a^{k}_{0})=(x_{0},\bar{u}(0),\bar{a}(0))$,
\begin{equation}\label{e:5.1}
y^{k}_{j+1}=y^{k}_{j}+h_{k}f_{2}(b^{k}_{j},x^{k}_{j})\;\mbox{ with }\;(y^{k}_{0},b^{k}_{0})=(0,\bar{b}(0)),
\end{equation}
\begin{equation}\label{e:5.1a}
g_{i}(x^{k}_{k}-u^{k}_{k})\geq 0,\,\,\,\,\,i=1,\ldots,s,
\end{equation}
\begin{equation}\label{e:5.2}
\lVert (x^{k}_{j},y^{k}_{j},u^{k}_{j},a^{k}_{j},b^{k}_{j})-(\bar{x}(t^{k}_{j}),\bar{y}(t^{k}_{j}),\bar{u}(t^{k}_{j}),\bar{a}(t^{k}_{j}),\bar{b}(t^{k}_{j})) \rVert\leq\dfrac{\varepsilon}{2},\,\,\,\,\, j=0,\ldots,k-1,
\end{equation}
\begin{equation}\label{e:5.3}
\begin{aligned}
&\sum\limits_{j=0}^{k-1}\int\limits_{t^{k}_{j}}^{t^{k}_{j+1}}\Big(\Big\lVert \dfrac{x^{k}_{j+1}-x^{k}_{j}}{h_{k}}-\dot{\bar{x}}(t) \Big\rVert^{2}+\Big\lVert \dfrac{y^{k}_{j+1}-y^{k}_{j}}{h_{k}}-\dot{\bar{y}}(t) \Big\rVert^{2}\Big)\,dt +\\
&\sum\limits_{j=0}^{k-1}\int\limits_{t^{k}_{j}}^{t^{k}_{j+1}}\Big(\Big\lVert \dfrac{u^{k}_{j+1}-u^{k}_{j}}{h_{k}}-\dot{\bar{u}}(t) \Big\rVert^{2}+\Big\lVert \dfrac{a^{k}_{j+1}-a^{k}_{j}}{h_{k}}-\dot{\bar{a}}(t) \Big\rVert^{2}+\Big\lVert \dfrac{b^{k}_{j+1}-b^{k}_{j}}{h_{k}}-\dot{\bar{b}}(t) \Big\rVert^{2}\Big)\,dt\leq \dfrac{\varepsilon}{2}.
\end{aligned}
\end{equation}\vspace*{-0.1in}

To proceed further, we first have to make sure that problems $(P_k)$ admit optimal solution.\vspace*{-0.1in}
\begin{proposition}[\bf existence of optimal solutions to discrete approximations]\label{discri:opti} Suppose that the standing assumptions $(\mathcal{H}_{1})-(\mathcal{H}_{3})$ are satisfied around the given r.i.l.m.\ $\oz(\cdot)$. Then each problem $(P_k)$ has an optimal solution provided that $k\in\N$ is sufficiently large.
\end{proposition}\vspace*{-0.1in}
{\bf Proof}. Observe that the $W^{1,2}$-strong convergence results of Theorem~\ref{strong:conver} and the construction of problems $(P_k)$ allow us to conclude that the sets of feasible solutions to $(P_k)$ are nonempty whenever $k$ is sufficiently large. Now we fix $k\in\N$ and show that set of feasible solutions to $(P_k)$ is bounded. Indeed, pick a sequence $z^{\nu}=(x^{\nu}_{0},\ldots,x^{\nu}_{k},y^{\nu}_{0},\ldots,y^{\nu}_{k},u^{\nu}_{0},\ldots,u^{\nu}_{k},a^{\nu}_{0},\ldots,a^{\nu}_{k},b^{\nu}_{0},\ldots,b^{\nu}_{k})$ of feasible solutions to
$(P_{k})$ that converges to some $z=(x_{0},\ldots,x^{k},y_{0},\ldots,y^{k},u_{0},\ldots,u^{k},a_{0},\ldots,a^{k},b_{0},\ldots,b^{k})$ as $ \nu\to 0 $ and show that $z $ is feasible to $ (P_{k}) $ as well. Observe that $ g_{i}(x_{j}-u_{j})=\lim\limits_{\nu\to\infty}g_{i}(x^{\nu}_{j}-u^{\nu}_{j})\geq 0 $ for all $ i=1,\ldots,s $ and $ j=0,\ldots,k-1 $, which ensures that $ x_{j}\in C(u_{j}) $. It follows from the closed graph property of the normal cone mapping in \eqref{e:5.0} with $F_1$ taken from \eqref{F1} and the continuity of the functions $ f_{1} $ and $ f_{2} $ that
\begin{equation*}
-\dfrac{x_{j+1}-x_{j}}{h_{k}}-f_{1}(a_{j},x_{j})-y_{j+1}\in N_{C}(x_{j}-u_{j}),\,\,\,y_{j+1}=y_{j}+h_{k}f_{2}(b_{j},x_{j}),
\end{equation*}
or equivalently, that $ x_{j+1}\in x_{j}+h_{k}F_{1}(x_{j},y_{j},u_{j},a_{j}) $ for all $ j=0,\ldots,k-1$, which verifies the claimed closedness of the feasible solution set. To conclude the proof of the proposition, we notice that the latter set is bounded due to the constraints in \eqref{e:5.2}. Thus the existence of optimal solutions in $(P_k)$ is ensured by the classical Weierstrass theorem in finite-dimensional spaces. $\h$\vspace*{-0.03in}

Now we are ready to establish the desired strong convergence of extended discrete optimal solutions of $(P_k)$ to the prescribed r.i.l.m.\ of the original problem $(P)$.\vspace*{-0.1in}
\begin{theorem}[\bf strong convergence of extended discrete optimal solutions]\label{disc:conver} Let $\bar{z}(\cdot)=(\bar{x}(\cdot), \bar{y}(\cdot) ,\\ \bar{u}(\cdot),\bar{a}(\cdot),\bar{b}(\cdot))$ be an r.i.l.m.\ for problem $(P)$. In addition to the standing assumptions $ (\mathcal{H}_{1})-(\mathcal{H}_{3})$ imposed along $\oz(\cdot)$, suppose that the cost functions $ \varphi $ and $ l_{0}\equiv l $ are continuous at $ \bar{x}(T) $ and at $(t,\bar{z}(t),\dot{\bar{z}}(t)) $ for a.e.\ $ t\in[0,T] $, respectively. Then any sequence of optimal solutions $\bar{z}^{k}(\cdot)=(\bar{x}^{k}(\cdot),\bar{y}^{k}(\cdot),\bar{u}^{k}(\cdot),\bar{a}^{k}(\cdot),\bar{b}^{k}(\cdot)) $ to $ (P_{k}) $, being piecewise linearly extended to $[0,T] $, converges to $ \bar{z}(\cdot) $ strongly $ W^{1,2}([0,T],\mathbb{R}^{3n+m+d}) $ as $k\to\infty$.
\end{theorem}\vspace*{-0.1in}
{\bf Proof}. We know from Proposition~\ref{discri:opti} that each problem $(P_k)$ admits optimal solutions $\oz^k(\cdot)$ for large $k\in\N$. Extend any $\oz^k(\cdot)$ piecewise linearly to the continuous-time interval $[0,T]$. We aim at verifying that
\begin{equation}\label{e:5.4}
\lim\limits_{k\to\infty}\int\limits_{0}^{T}\lVert(\dot{\bar{x}}^{k}(t),\dot{\bar{y}}^{k}(t)
,\dot{\bar{u}}^{k}(t),\dot{\bar{a}}^{k}(t),\dot{\bar{b}}^{k}(t))-(\dot{\bar{x}}(t),\dot{\bar{y}}(t)
,\dot{\bar{u}}(t),\dot{\bar{a}}(t),\dot{\bar{b}}(t))\rVert^{2}\,dt=0,
\end{equation}
which clearly yields the convergence of the quintuple $(\bar{x}^{k}(\cdot),\bar{y}^{k}(\cdot),\bar{u}^{k}(\cdot),\bar{a}^{k}(\cdot),\bar{b}^{k}(\cdot)) $ to  $(\bar{x}(\cdot),\bar{y}(\cdot),\bar{u}(\cdot),\bar{a}(\cdot),\bar{b}(\cdot)) $ in the norm topology of $ W^{1,2}([0,T],\mathbb{R}^{3n+m+d})$ as $k\to\infty$.\vspace*{-0.05in}

To prove \eqref{e:5.4}, we argue by contraposition and suppose that it fails, i.e., the limit in \eqref{e:5.4}, along a subsequence (without relabeling) equals to some $\xi>0$. The
weak compactness of the unit ball in $L^{2}([0,T],\mathbb{R}^{3n+m+d})$ allows us to find $ (\upsilon^{x}(\cdot),\upsilon^{y}(\cdot),\upsilon^{u}(\cdot),\upsilon^{a}(\cdot),\upsilon^{b}(\cdot)) $ such that the sequence of derivatives $\{\dot{\bar{z}}^{k}(\cdot):=(\dot{\bar{x}}^{k}(\cdot),\dot{\bar{y}}^{k}(\cdot),\dot{\bar{u}}^{k}(\cdot),\dot{\bar{a}}^{k}(\cdot),\dot{\bar{b}}^{k}(\cdot))\}$ converges weakly to a quintuple $(\upsilon^{x}(\cdot),\upsilon^{y}(\cdot),\upsilon^{u}(\cdot),\upsilon^{a}(\cdot),\upsilon^{b}(\cdot))\in L^{2}([0,T],\mathbb{R}^{3n+m+d})$ as $k\to\infty$. Defining further $\tilde{z}(\cdot)=(\tilde{x}(\cdot),\tilde{y}(\cdot),\tilde{u}(\cdot),\tilde{a}(\cdot),\tilde{b}(\cdot))$ by
\begin{equation*}
\tilde{z}(t):=(x_{0},0,\bar{u}(0),\bar{a}(0),\bar{b}(0))+\int\limits_{0}^{t}(\upsilon^{x}(s),\upsilon^{y}(s),\upsilon^{u}(s),\upsilon^{a}(s),\upsilon^{b}(s))\,ds,\quad t\in[0,T],
\end{equation*}
we get that $\dot{\tilde{z}}(t)=(\upsilon^{x}(t),\upsilon^{y}(t),\upsilon^{u}(t),\upsilon^{a}(t),\upsilon^{b}(t))$ for a.e.\ $t\in[0,T]$. Arguing now as in the proof of Theorem~\ref{exist:opti} shows that the arc $\tilde{z}(\cdot)$ is feasible to the original problem $(P)$ and hence to the relaxed one $(R)$.\vspace*{-0.05in}

Next we check that $\tilde{z}(\cdot)$ satisfies the localization conditions in \eqref{loc} relative to $\bar{z}(\cdot)$. Indeed, the first condition in \eqref{loc} follows directly from the passage to the limit in \eqref{e:5.2} as $k\to\infty$. To justify the second condition in \eqref{loc}, we pass to the limit in \eqref{e:5.3} due to the established weak convergence of the derivatives $\bar{z}^{k}(\cdot)\to\dot{\tilde{z}}(\cdot)$ and the lower semicontinuity of the norm function in $L^2$. This tells us that
\begin{equation*}
\begin{aligned}
&\int\limits_{0}^{T}\lVert (\dot{\tilde{x}}(t),\dot{\tilde{y}}(t),\dot{\tilde{u}}(t),\dot{\tilde{a}}(t),\tilde{\bar{b}}(t))-(\dot{\bar{x}}(t),\dot{\bar{y}}(t),\dot{\bar{u}}(t),\dot{\bar{a}}(t),\dot{\bar{b}}(t)) \rVert\\
&\leq
\liminf\limits_{k\to \infty}\sum\limits_{j=0}^{k-1}\int\limits_{t^{k}_{j}}^{t^{k}_{j+1}}\Big(\Big\lVert \dfrac{x^{k}_{j+1}-x^{k}_{j}}{h_{k}}-\dot{\bar{x}}(t) \Big\rVert^{2}+\Big\lVert \dfrac{y^{k}_{j+1}-y^{k}_{j}}{h_{k}}-\dot{\bar{y}}(t) \Big\rVert^{2}\Big)\,dt\\
&+\sum\limits_{j=0}^{k-1}\int\limits_{t^{k}_{j}}^{t^{k}_{j+1}}\Big(\Big\lVert \dfrac{u^{k}_{j+1}-u^{k}_{j}}{h_{k}}-\dot{\bar{u}}(t) \Big\rVert^{2}+\Big\lVert \dfrac{a^{k}_{j+1}-a^{k}_{j}}{h_{k}}-\dot{\bar{a}}(t) \Big\rVert^{2}+\Big\lVert \dfrac{b^{k}_{j+1}-b^{k}_{j}}{h_{k}}-\dot{\bar{b}}(t) \Big\rVert^{2}\Big)\,dt\leq\dfrac{\varepsilon}{2},
\end{aligned}
\end{equation*}
and thus we are done with \eqref{loc}. Furthermore, it follows from the construction of the relaxed problem $(R)$ in \eqref{R} due to the convexity of $\Hat{l}_G$ in the velocity variables, the established weak convergence of the extended discrete derivatives, and applications of the Mazur theorem as in the proof of Theorem~\ref{strong:conver} that
\begin{equation*}
\begin{aligned}
&\int\limits_{0}^{T}\hat{l}_{G}(t,\tilde{x}(t),\tilde{y}(t),\tilde{u}(t),\tilde{a}(t),\tilde{b}(t),\dot{\tilde{x}}(t),\dot{\tilde{y}}(t),\dot{\tilde{u}}(t),\dot{\tilde{a}}(t),\dot{\tilde{b}}(t))\,dt\\
&\leq \liminf\limits_{k\to \infty}h_{k}\sum\limits_{j=0}^{k-1}l\Big(t^{k}_{j},\bar{x}^{k}_{j},\bar{y}^{k}_{j},\bar{u}^{k}_{j},\bar{a}^{k}_{j},\bar{b}^{k}_{j},\dfrac{\bar{x}^{k}_{j+1}-\bar{x}^{k}_{j}}{h_{k}},\dfrac{\bar{y}^{k}_{j+1}
-\bar{y}^{k}_{j}}{h_{k}},\dfrac{\bar{u}^{k}_{j+1}-\bar{u}^{k}_{j}}{h_{k}},\dfrac{\bar{a}^{k}_{j+1}-\bar{a}^{k}_{j}}{h_{k}},\dfrac{\bar{b}^{k}_{j+1}-\bar{b}^{k}_{j}}{h_{k}}\Big).
\end{aligned}
\end{equation*}
We also observe in this way, with taking into account the above definition of $\xi$, that
\begin{equation}\label{e:5.5}
\begin{aligned}
\hat{J}[\tilde{z}]+\dfrac{\xi}{2}&=\varphi(\tilde{x}(T))+\int\limits_{0}^{T}\hat{l}_{G}\Big(t,\tilde{x}(t),\tilde{y}(t),\tilde{u}(t),\tilde{a}(t),\tilde{b}(t),\dot{\tilde{x}}(t),\dot{\tilde{y}}(t),
\dot{\tilde{u}}(t),\dot{\tilde{a}}(t),\dot{\tilde{b}}(t)\Big)\,dt+\dfrac{\xi}{2}\\
&\leq\liminf\limits_{k\to \infty} J_{k}(\bar{z}^{k}).
\end{aligned}
\end{equation}
Employing Theorem~\ref{strong:conver} gives us a sequence $\{(x^{k}(\cdot),y^{k}(\cdot),u^{k}(\cdot),a^{k}(\cdot),b^{k}(\cdot))\}$ of feasible solutions to $(P_{k})$  that strongly
$W^{1,2}$-approximates the r.i.l.m.\ $ (\bar{x}(\cdot),\bar{y}(\cdot),\bar{u}(\cdot),\bar{a}(\cdot),\bar{b}(\cdot))$, which is a feasible solution to $(P)$. The imposed continuity assumptions on $ \varphi $ and $ l $ yield
\begin{equation}\label{e:5.6}
\lim\limits_{k\to \infty} J_{k}[x^{k},y^{k},u^{k},a^{k},b^{k}]=J[\bar{x},\bar{y},\bar{u},\bar{a},\bar{b}].
\end{equation}
On the other hand, it follows from the optimality of $ \bar{z}^{k}(\cdot):= (\bar{x}^{k}(\cdot),\bar{y}^{k}(\cdot),\bar{u}^{k}(\cdot),\bar{a}^{k}(\cdot),\bar{b}^{k}(\cdot))$ in $(P_k)$ that
\begin{equation}\label{e:5.7}
J_{k}[\bar{z}^{k}]\leq J_k[x^{k},y^{k},u^{k},a^{k},b^{k}]\;\mbox{ for each }\;k\in\mathbb{N}.
\end{equation}
Combining finally the relationships in \eqref{e:5.5}--\eqref{e:5.7}, we conclude that
\begin{equation*}
\tilde{J}[\tilde{x},\tilde{y},\tilde{u},\tilde{a},\tilde{b}]<\tilde{J}[\tilde{x},\tilde{y},\tilde{u},\tilde{a},\tilde{b}]+\dfrac{\xi}{2}\leq J[\bar{z}]=\hat{J}[\bar{z}].
\end{equation*}
Due to the choice of $\xi>0$ above, the latter clearly contradicts the fact that $ \bar{z}(\cdot) $ is an r.i.l.m.\ for problem $(P)$ and thus verifies the
limiting condition \eqref{e:5.4}. This completes the proof of the theorem. $\h$ \vspace*{-0.2in}
%%%%%%%%%%%%%%%%%%%%%%%%%%%%%%%%%%%%%%%%%%%%%%%%%%%%%%%%%%%%%%%%%%%%%%%%%%%%%%%%%%%%%%%%%%%%%%%%%%%%%%%%%%%%%%%%%%%%%%%%%%%%%%%%%%%%%%%%%%%%%%%%%%%%%%%%%%%%%%%%%%%%%%%%%%%%%%%%%%%%%%%%%%%
%%%%%%%%%%%%%%%%%%%%%%%%%%%%%%%%%%%%%%%%%%%%%%%%%%%%%%%%%%%%%%%%%%%%%%%%%%%%%%%%%%%%%%%%%

\section{Generalized Differentiation and Second-Order Calculations}\label{sec:2va}
\setcounter{equation}{0}\vspace*{-0.1in}

This section briefly recalls some tools of first-order and second-order generalized differentiation in variational analysis, which are instrumental in deriving necessary optimality conditions for the discrete-time and continuous-time sweeping control problems formulated above. The reader can find more details and references in the books \cite{Mord,b3,rw} for the first-order and in \cite{Mord,b3} for the second-order issues. Observe that, although the initial data of problem $(P)$ and its discrete approximations are smooth and/or convex, the intrinsic source of nonsmoothness comes from the sweeping dynamics and its discretization, which lead us to nonconvex graphical sets and require the usage of robust generalized differentiation with adequate properties. It occurs that the most appropriate constructions for these purposes are the robust nonconvex notions introduced by the third author, while their convexification fails the needed results, particularly of the second-order.\vspace*{-0.05in}

Given a nonempty set $\O\subset\R^n$ locally closed around $\ox\in\O$, consider the associated projection operator \eqref{proj} and recall that the (basic/limiting/Mordukhovich) {\em normal cone} to $\O$ at $\ox$ is defined by
\begin{eqnarray}\label{e:Mor-nc}
N_\O(\ox):=\big\{v\in\R^n\big|\;\exists\,x_k\to\ox,\;w_k\in\Pi_\O(x_k),\;\al_k\ge 0\;\mbox{ s.t. }\;\al_k(x_k-w_k)\to v\;\mbox{ as }\;k\to\infty\big\}
\end{eqnarray}
with $N_\O(\ox):=\emp$ for $\ox\notin\O$. If $\O$ is prox-regular, the normal cone \eqref{e:Mor-nc} agrees with the proximal normal cone \eqref{pnc}, and they both reduce to the normal cone of convex analysis when $\O$ is convex. Furthermore, the convex closure of \eqref{e:Mor-nc} gives us the Clarke normal cone \cite{clsw}. Note that, despite the nonconvexity of \eqref{e:Mor-nc}, this normal cone and the associated subdifferential and coderivative constructions for functions and multifunctions are {\rm robust} and enjoy {\em full calculi} based on variatiional/extremal principles of variational analysis. This is not the case for the proximal constructions and may also fail for Clarke's ones without imposing additional interiority-type assumptions.\vspace*{-0.05in}

Considering a set-valued mapping/multifunction $ F : \mathbb{R}^{n}\rightrightarrows \mathbb{R}^{q} $ with $\dom F:=\{x\in\R^n\;|\;F(x)\ne\emp\}$ and taking a point $ (\bar{x},\bar{y})$ from the graph
\begin{equation*}
\gph\,F:=\big\{(x,y)\in\mathbb{R}^{n}\times \mathbb{R}^{q}\;\big|\;y\in F(x)\big\},
\end{equation*}
the {\em coderivative} $ D^{\ast}F(\bar{x},\bar{y}) : \mathbb{R}^{q}\rightrightarrows \mathbb{R}^{n} $ of $ F $ at $ (\bar{x},\bar{y}) $ is defined by
\begin{equation}\label{cod}
D^{\ast}F(\bar{x},\bar{y})(u):=\big\{ \upsilon\in\mathbb{R}^{n}\;\big|\;(\upsilon,-u)\in N_{{\rm\small gph}\,F}(\bar{x},\bar{y})\}\;\mbox{ for all }\; u\in\mathbb{R}^{q},
\end{equation}
where we skip $\oy$ if $F$ is single-valued. If in the latter case $F$ is continuously differentiable around $ \bar{x} $, then
\begin{equation*}
D^{\ast}F(\bar{x})(u):=\big\{\nabla F(\bar{x})^{\ast}u\big\}\;\mbox{ for all }\; u\in\mathbb{R}^{q}.
\end{equation*}
Given further an extended-real-valued function $ \varphi:\mathbb{R}^{n}\to \bar{\mathbb{R}}$  finite at $ \bar{x} $, the (first-order) {\em subdifferential} of $ \varphi $ at $ \bar{x} $ is defined geometrically by
\begin{equation}\label{sub}
\partial \varphi(\bar{x}) : = \big\{ \upsilon\in\mathbb{R}^{n}\;\big|\; (\upsilon,-1)\in N_{{\rm\small epi}\,\ph}(\bar{x},\varphi(\bar{y}))\big\}
\end{equation}
while admitting various analytical representations given in \cite{b3,rw}. Following \cite{m92}, the {\em second-order subdifferential} (generalized Hessian) of $ \varphi $ at $ \bar{x} $ relative to $ \bar{y}\in \varphi(\bar{x})$ is defined by
\begin{equation}\label{2nd}
\partial^{2}\varphi(\bar{x},\bar{y})(u):=\big(D^{\ast}\partial\varphi\big)(\bar{x},\bar{y})(u)\;\mbox{ for all }\;u\in\mathbb{R}^{q}.
\end{equation}
If $ \varphi $ is twice continuously differentiable around $ \bar{x} $, then we get
\begin{equation*}
\partial^{2}\varphi(\bar{x})(u):=\big\{\nabla^{2}\varphi(\bar{x})(u)\big\}\;\mbox{ for all }\; u\in\mathbb{R}^{q}
\end{equation*}
via the classical (symmetric) Hessian of $\ph$ at $\ox$. Note that replacing the limiting normal cone in \eqref{cod} and \eqref{2nd} by its convexification (and hence by Clarke's normal cone) to graphical sets {\em dramatically enlarges} the corresponding constructions and makes them unusable for further analysis and applications. In particular, it follows from \cite[Theorem~3.5]{r} that Clarke's normal cone is a linear {\em subspace} of the maximal dimension if the graph of $F$  is a ``Lipschitzian manifold," which is the case of, e.g., subgradient mappings $F=\partial\ph$ in \eqref{sub} generated by {\em prox-regular} functions; see \cite{Mord,b3,rw} for more discussions. In contrast, the second-order subdifferential \eqref{2nd} provides an adequate machinery of second-order variational analysis in applications to controlled sweeping processes. We largely use below the following theorem, which gives us a {\em precise computation} of the velocity mapping associated with the integro-differential sweeping dynamics of our study entirely in terms of the given data of $(P)$. Due to the very structure of the sweeping dynamics, the obtained result contains the computation of the second-order subdifferential of the indicator function in \eqref{2nd}.\vspace*{-0.1in}
\begin{theorem}[\bf second-order calculation for integro-differential sweeping dynamics]\label{Th:co-nc} Define the set-valued mapping $ F_{h} : \mathbb{R}^{3n+m+d}\rightrightarrows \mathbb{R}^{n} $ by
\begin{equation}\label{Fh}
F_{h}(x,y,u,a,b):=F_{1}(x,y,u,a)+hf_{2}(b,x)=N_{C}(x-u)+f_{1}(a,x)+y+hf_{2}(b,x),
\end{equation}
where $F_1$ is given in \eqref{F1} with $f_1$, $f_2$, and $C$ taken from \eqref{e:1} and \eqref{e:mset}, respectively, and $h>0$. Pick any $ u\in \mathbb{R}^{n} $ with $ x-u \in C $, take $ w\in N_{C}(x-u)+f_{1}(a,x)+y+hf_{2}(b,x)$, and suppose that $f_1,f_2\in{\cal C}^1$, $g:=(g_1,\ldots,g_s)\in{\cal C}^2$ around the corresponding points with full rank of the Jacobian matrix $\nabla g(x-u)$. Let $\lambda:=(\lambda_{1},\ldots,\lambda_{s})\in\R^s$ be a unique vector with nonnegative components satisfying the equation
\begin{equation}\label{e:4.1}
-\nabla g(x-u)^{*}\lambda=w-f_{1}(a,x)-y-hf_{2}(b,x).
\end{equation}
Then the coderivative of $F_h$ is calculated by the formula
\begin{equation}\label{e:4.2}
\begin{aligned}
D^{*}F_{h}(x,y,u,a,b,w)(z)=&\Big\{ \nabla_{x} f_{1}(a,x)^{*}z+h\nabla_{x} f_{2}(b,x)^{*}z-\Big(\sum\limits_{i=1}^{s}\lambda_{i}\nabla_{x}^{2}g_{i}(x-u)\Big)z-\nabla_{x}g(x-u)^{*}\sigma,\\
&\Big(\sum\limits_{i=1}^{s}\lambda_{i}\nabla_{u}^{2}g_{i}(y-u)\Big)z+\nabla_{u}g(y-u)^{*}\sigma,\nabla_{a} f_{1}(a,x)^{*}z,h\nabla_{b} f_{2}(b,x)^{*}z\Big\}
\end{aligned}
\end{equation}
for all $z\in\dom D^*N_{C}(x-u,w-f_{1}(a,x)-y-hf_{2}(b,x))$, where the coderivative domain is represented as
\begin{equation}\label{e:4.3}
\begin{aligned}
&\dom D^*N_{C}(x-u,w-f_{1}(a,x)-y-hf_{2}(b,x))=\big\{z\in\R^n\;\big|\;\exists\,\lambda\in\R^s_+\;\mbox{ such that }\\
&-\nabla g(x-u)\lambda=w-f_{1}(a,x)-y-hf_{2}(b,x),\;\lambda_{i}\langle \nabla\,g_{i}(x-u),z\rangle=0\;\mbox{ for all }\;i=1,\ldots,s\big\},
\end{aligned}
\end{equation}
and where we have in \eqref{e:4.2} that $\sigma_{i}=0$ if either $g_{i}(x-u)>0$ or $\lambda_{i}=0$ with $\langle \nabla\,g_{i}(x-u),z\rangle>0 $, and that $\sigma_{i}\geq 0$ if $g_{i}(x-u)=0$, $\lambda_{i}=0$ with $\langle \nabla\,g_{i}(x-u),z\rangle<0$.
\end{theorem}\vspace*{-0.1in}
{\bf Proof}. It follows the lines in the proof of \cite[Theorem~6.2]{b1} but applying now to the different form of the velocity mapping \eqref{Fh} employed in this paper. We skip the details for brevity while mentioning that the coderivative sum and chain rules taken from \cite[Theorems~1.62 and 1.66]{Mord} allow us to derive the claimed equalities in \eqref{e:4.2} and \eqref{e:4.3} to the second-order calculations developed in \cite[Theorem~3.3]{hos}.$\h$\vspace*{-0.2in}
%%%%%%%%%%%%%%%%%%%%%%%%%%%%%%%%%%%%%%%%%%%%%%%%%%%%%%%%%%%%%%%%%%%%%%%%%%%%%%%%%%%%%%%%%%%%%%%%%%%%%%%%%%%%%%%%%%%%%%%%%%%%%%%%%%%%%%%%%%%%%%%%%%%%%%%%%%%%%%%%%%%%%%%%%%%%%%%%%%%%%%%%%%%%%%%%%%
%%%%%%%%%%%%%%%%%%%%%%%%%%%%%%%%%%%%%%%%%%%%%%%%%%%%%%%%%%%%%%%%%%%%%%%%%%%%%%%%%%%%%%%%%%%%%%%%%%%%%%%%%%%%%%%%%%%%%%%%%%%%%%%%%%%%%%%%%%%%%%%%%%%%%%

\section{Necessary Optimality Conditions for Discrete Problems}\label{sec:disc-nec}
\setcounter{equation}{0}\vspace*{-0.1in}

In this section we obtain necessary optimality conditions for (global) optimal solutions to the discrete approximation problems $(P_k)$ constructed in Section~\ref{sec:disc-opt} for each $k\in\N$. Due to the $W^{1,2}$-strong convergence result established in Theorem~\ref{disc:conver}, the obtained optimality conditions for problems $(P_k)$ with sufficiently large approximation indices $k$ can be viewed as necessary {\em suboptimality} conditions for the original continuous-time optimal control problem $(P)$ with the integro-differential sweeping dynamics \eqref{e:1}. In practice, such conditions for $(P_k)$ provide sufficient information for designing numerical algorithms to solve the original problem $(P)$. Nevertheless, in Section~\ref{sec:nec-sw} we furnish the rigorous  procedure to derive precise necessary optimality conditions for relaxed intermediate local minimizers in $(P)$ by passing to the limit as $k\to\infty$ from the necessary optimality conditions for $(P_k)$ obtained below.\vspace*{-0.05in}

In what follows we present two sets of necessary optimality conditions in problems $(P_k)$ for each fixed $k\in\N$. The first theorem concerns a class of more general problems given in form $(P_k)$ with an arbitrary closed-graph velocity mapping $F_h$. The obtained results can be treated as extended {\em discrete-time Euler-Lagrange conditions} for discrete approximations of integro-differential sweeping control problems with the adjoint systems described via the basic normal cone \eqref{e:Mor-nc} to graphs of velocity mappings.\vspace*{-0.1in}
\begin{theorem}[\bf extended Euler-Lagrange conditions for discrete optimal solutions]\label{discr:N.C} Fix $k\in\N$ to be sufficiently large and pick any optimal solution $\bar{z}^{k}:=(x_{0},\ldots,\bar{x}^{k}_{k},y_{0},\ldots,\bar{y}^{k}_{k},\bar{u}_{0}^{k},\ldots,\bar{u}^{k}_{k},a_{0},\ldots,\bar{a}^{k}_{k},b_{0},\ldots,\bar{b}^{k}_{k})$ to problem $(P_k)$ governed by an arbitrary closed-graph multifunction in \eqref{e:5.0}, while not taken from \eqref{Fh0}. Suppose that the cost functions $\ph$ and $l_0$ from \eqref{J0} are locally Lipschitzian around the corresponding components of the optimal solution for all $t\in\Delta_k$ $($with the index $t$ dropped below$)$, and that the mappings $f_1,f_2,g_i$ are continuously differentiable around the optimal points. Denote
\begin{equation}\label{e:6.1a}
\theta^{xk}_{j}:= \int\limits_{t^{k}_{j}}^{t^{k}_{j+1}}\Big( \dfrac{\bar{x}^{k}_{j+1}-\bar{x}^{k}_{j}}{h_{k}}-\dot{\bar{x}}(t)\Big)\,dt,\;\theta^{uk}_{j}:=\int\limits_{t^{k}_{j}}^{t^{k}_{j+1}}
\Big(\dfrac{\bar{u}^{k}_{j+1}-\bar{u}^{k}_{j}}{h_{k}}-\dot{\bar{u}}(t)\Big)\,dt,\;\theta^{ak}_{j}=\int\limits_{t^{k}_{j}}^{t^{k}_{j+1}}
\Big(\dfrac{\bar{a}^{k}_{j+1}-\bar{a}^{k}_{j}}{h_{k}}-\dot{\bar{a}}(t)\Big)\,dt,
\end{equation}
\begin{equation}\label{e:6.1}
\theta^{bk}_{j}:=\int\limits_{t^{k}_{j}}^{t^{k}_{j+1}}\Big(\dfrac{\bar{b}^{k}_{j+1}-\bar{b}^{k}_{j}}{h_{k}}-\dot{\bar{b}}(t)\Big)\,dt,\;\mbox{ and }\;\theta^{yk}_{j}:= \int\limits_{t^{k}_{j}}^{t^{k}_{j+1}}\Big(\dfrac{\bar{y}^{k}_{j+1}-\bar{y}^{k}_{j}}{h_{k}}-\dot{\bar{y}}(t)\Big)\,dt,\quad j=0,\ldots,k-1.
\end{equation}
Then there are $\lambda^{k}\geq 0 $, $ \alpha^{k}=(\alpha^{k}_{1},\ldots,\alpha^{k}_{s})\in\mathbb{R}^{s}_{+} $, $ p^{k}_{j}=(p^{xk}_{j},p^{yk}_{j},p^{dk}_{j},p^{uk}_{j},p^{ak}_{j},p^{bk}_{j})\in \mathbb{R}^{4n+m+d}$, $ j=0,\ldots,k $,
\begin{equation*}
(w^{xk}_{j},w^{uk}_{j},w^{ak}_{j},w^{bk}_{j},v^{xk}_{j},v^{uk}_{j},v^{ak}_{j},v^{bk}_{j})\in \partial l_{0}\Big(\bar{x}^{k}_{j},\bar{u}^{k}_{j},\bar{a}^{k}_{j},\bar{b}^{k}_{j},\dfrac{\bar{x}^{k}_{j+1}-\bar{x}^{k}_{j}}{h_{k}},\dfrac{\bar{u}^{k}_{j+1}-\bar{u}^{k}_{j}}{h_{k}},\dfrac{\bar{a}^{k}_{j+1}
-\bar{a}^{k}_{j}}{h_{k}},\dfrac{\bar{b}^{k}_{j+1}-\bar{b}^{k}_{j}}{h_{k}}\Big)
\end{equation*}
such that the following conditions are satisfied:
\begin{equation}\label{e:6.2}
\lambda^{k}+\lVert \alpha^{k} \rVert+\sum\limits_{j=0}^{k-1}\lVert p^{kx}_{j} \rVert +\sum\limits_{j=0}^{k}\lVert p^{kd}_{j} \rVert +\lVert p^{ky}_{0} \rVert + \lVert p^{ku}_{0} \rVert+\lVert p^{ka}_{0} \rVert+\lVert p^{kb}_{0} \rVert \neq 0,
\end{equation}
\begin{equation}\label{e:6.3}
\alpha^{k}_{i}g_{i}(\bar{x}^{k}_{k}-\bar{u}^{k}_{k})=0,\,\,\,i=1,\ldots,s,
\end{equation}
\begin{equation}\label{e:6.4}
-p^{x}_{k}\in\lambda^{k}\partial\varphi(\bar{x}^{k}_{k})-\sum\limits_{i=1}^{s}\alpha^{k}_{i}\nabla_{x}g_{i}(\bar{x}^{k}-\bar{u}^{k}),\;\;
-p^{uk}_{k}=\sum\limits_{i=1}^{s}\alpha_{i}^k\nabla g_{i}(\bar{x}^{k}_{k}-\bar{u}^{k}_{k}),
\end{equation}
\begin{equation}\label{e:6.6}
p^{yk}_{k}=0,\,\,\,p^{ak}_{k}=0,\,\,\,p^{bk}_{k}=0,
\end{equation}
\begin{equation}\label{e:6.7}
p^{y}_{j+1}=\lambda^{k}h_{k}^{-1}\theta^{y}_{j}+h_{k}^{-1}p^{d}_{j+1},\;\;
p^{u}_{j+1}=\lambda^{k}(v^{u}_{j}+ h_{k}^{-1}\theta^{u}_{j}),\,\,\, j=0,\ldots,k-1,
\end{equation}
\begin{equation}\label{e:6.9}
p^{a}_{j+1}=\lambda^{k}(v^{a}_{j}+ h_{k}^{-1}\theta^{a}_{j}),\;\mbox{ and }\;
p^{b}_{j+1}=\lambda^{k}(v^{b}_{j}+ h_{k}^{-1}\theta^{b}_{j}),\,\,\, j=0,\ldots,k-1.
\end{equation}
Furthermore, for all $j=0,\ldots,k-1$ we have the extended discrete-time Euler-Lagrange inclusions:
\begin{equation}\label{e:6.11}
\Big(\dfrac{p^{x}_{j+1}-p^{x}_{j}}{h_{k}}-\lambda^{k}w^{x}_{j}+h_{k}^{-1}\xi_{j}^{x}p^{d}_{j+1},\dfrac{p^{y}_{j+1}-p^{y}_{j}}{h_{k}}-\dfrac{p^{d}_{j+1}}{h_{k}},\dfrac{p^{u}_{j+1}-p^{u}_{j}}{h_{k}}
-\lambda^{k}w^{u}_{j},\dfrac{p^{a}_{j+1}-p^{a}_{j}}{h_{k}}-\lambda^{k}w^{a}_{j},
\end{equation}
\begin{equation*}
\dfrac{p^{b}_{j+1}-p^{b}_{j}}{h_{k}}-\lambda^{k}w^{b}_{j}+h_{k}^{-1}\xi_{j}^{b}p^{d}_{j+1},
p^{x}_{j+1}-\lambda^{k}(v^{x}_{j}+h_{k}^{-1}\theta^{x}_{j})\Big)\notag\in N_{{\rm\small gph}\,F_h}\Big(\bar{x}^{k}_{j},\bar{y}^{k}_{j},\bar{u}^{k}_{j},\bar{a}^{k}_{j},\bar{b}^{k}_{j},\dfrac{\bar{x}^{k}_{j+1}-\bar{x}^{k}_{j}}{-h_{k}}\Big),
\end{equation*}
where the vectors $\xi^{x}_{j}$ and $\xi^{b}_{j}$ are defined by
\begin{equation}\label{xi}
\xi^{x}_{j}:=\nabla_{x}f_{2}(\bar{b}^{k}_{j},\bar{x}^{k}_{j}),\;\mbox{ and }\;\xi^{b}_{j}:=\nabla_{b}f_{2}(\bar{b}^{k}_{j},\bar{x}^{k}_{j}),\quad j=0,\ldots,k-1.
\end{equation}
\end{theorem}\vspace*{-0.1in}
{\bf Proof}. Denote $z:=(x_{0},\ldots,x^{k},y_{0},\ldots,y^{k},u_{0},\ldots,u^{k},a_{0},\ldots,a^{k},b_{0},\ldots,b^{k},X^{k}_{0},\ldots,X^{k}_{k-1},Y_{0},Y_{1},\ldots,Y^{k},\\
U_{0},\ldots,U_{k-1},A_{0},A_{1},\ldots,A_{k-1},B_{0},B_{1},\ldots,B_{k-1}) $, where $ x_{0} $ is fixed, and consider the following problem of mathematical programming (MP): minimize
\begin{equation*}
\begin{aligned}
\varphi_{0}(z):=\varphi(x^{k}_{k})+h_{k}\sum\limits_{j=0}^{k-1}l(x_{j},y_{j},u_{j},a_{j},b_{j},X_{j},Y_{j},U_{j},A_{j},B_{j})+\dfrac{1}{2}\sum\limits_{j=0}^{k-1}\int\limits_{t_{j}}^{t_{j+1}}\lVert (X_{j},Y_{j},U_{j},A_{j},B_{j})-\dot{\bar{z}}(t) \rVert^{2}\,dt
\end{aligned}
\end{equation*}
subject to the equality, inequality, and many geometric constraints of the graphical type defined by
\begin{equation*}
e_{j}^{x}(z):=x_{j+1}-x_{j}-h_{k}X_{j}=0,\;\;e_{j}^{y}(z):=y_{j+1}-y_{j}-h_{k}Y_{j}=0,\quad j=0,\ldots,k-1,
\end{equation*}
\begin{equation*}
e_{j}^{u}(z):=u_{j+1}-u_{j}-h_{k}U_{j}=0,\;\;e_{j}^{a}(z):=a_{j+1}-a_{j}-h_{k}A_{j}=0,\quad j=0,\ldots,k-1,
\end{equation*}
\begin{equation*}
e_{j}^{b}(z):=b_{j+1}-b_{j}-h_{k}B_{j}=0,\;\;d_{j}(z):=Y_{j}-f_{2}(b_{j},x_{j})=0,\quad j=0,\ldots,k-1,
\end{equation*}
\begin{equation*}
c_{i}(z):=-g_{i}(x^{k}-u^{k})\leq 0,\,\,\,i=1,\ldots,s,
\end{equation*}
\begin{equation*}
\phi_{j}(z):=\lVert(x_{j},y_{j},u_{j},a_{j},b_{j})-\bar{z}(t_{j}) \rVert-\dfrac{\varepsilon}{2}\leq 0,\quad j=0,\ldots,k-1,
\end{equation*}
\begin{equation*}
\phi_{k}(z):=\sum\limits_{j=0}^{k-1}\int\limits_{t_{j}}^{t_{j+1}}\lVert (X_{j},Y_{j},U_{j},A_{j},B_{j})-\dot{\bar{z}}(t)\rVert^{2}\,dt-\dfrac{\varepsilon}{2}\leq 0,
\end{equation*}
\begin{equation*}
z\in\Xi_{j}:=\big\{z\;\big|\;-X_{j} \in F_{h}(x_{j},y_{j},u_{j},a_{j},b_{j})\big\},\quad j=0,\ldots,k-1,
\end{equation*}
\begin{equation*}
z\in\Xi_{k}=\big\{z\;\big|\;x_{0}\;\mbox{ is fixed },\;y_{0}=0,\,\,\,(u_{0},a_{0},b_{0})=(\bar{u}(0),\bar{a}(0),\bar{b}(0))\big\}.
\end{equation*}
Applying now the necessary optimality conditions from \cite[Theorem~6.5 (ii)]{b3} together with the intersection rule for basic normals taken from \cite[Corollary~2.17]{b3}
to the chosen solution $\oz^k$ of problem $(MP)$, observe first that the inequality constraints in $(MP)$ defined by the functions $\phi_j$ for $j=0,\ldots,k$ are inactive for sufficiently large $k\in\N$ due to the $W^{1,2}$-strong convergence $\oz^k(\cdot)\to\oz(\cdot)$ established above in Theorem~\ref{disc:conver}. All of this allows us to find $ \lambda^{k}\geq 0 $, $ \alpha^{k}=(\alpha^{k}_{1},\ldots,\alpha^{k}_{s})\in \mathbb{R}^{s}_{+}$, $ p_{j}=(p^{x}_{j},p^{y}_{j},p^{u}_{j},p^{a}_{j},p^{b}_{j},p^{d}_{j})\in \mathbb{R}^{4n+m+d} $ as $j=1,\ldots,k$, and
\begin{equation*}
z^{\ast}_{j}=\big(x^{\ast}_{0j},\ldots,x^{\ast}_{kj},y^{\ast}_{0j},\ldots,y^{\ast}_{kj},u^{\ast}_{0j},\ldots,u^{\ast}_{kj},a^{\ast}_{0j},\ldots,a^{\ast}_{kj},b^{\ast}_{0j},\ldots,b^{\ast}_{kj},X^{\ast}_{0j},
\ldots,X^{\ast}_{(k-1)j},Y^{\ast}_{0j},\ldots,Y^{\ast}_{(k-1)j},
\end{equation*}
\begin{equation*}
U^{\ast}_{0j},\ldots,U^{\ast}_{(k-1)j},A^{\ast}_{0j},\ldots,A^{\ast}_{(k-1)j},B^{\ast}_{0j},\ldots,B^{\ast}_{(k-1)j}\big),\,\,\,j=0,\ldots,k-1,
\end{equation*}
which are not zero simultaneously, such that the following conditions are satisfied :
\begin{equation}\label{e:6.12}
z^{\ast}_{j}\in N_{\Xi_{j}}(\bar{z}^{k}),\,\,\,j=0,\ldots,k,
\end{equation}
\begin{equation}\label{e:6.14}
\alpha^{k}_{i} c_{i}(\bar{z}^{k})=0,\,\,\,i=1,\ldots,s,
\end{equation}
\begin{equation}\label{e:6.13}
-\sum\limits_{j=0}^{k}z^{\ast}_{j}\in \lambda^{k}\partial\varphi_{0}(\bar{z}^{k})+\sum\limits_{i=1}^{s}\alpha^{k}_{i}\nabla c_{i}(\bar{z}^{k}) + \nabla e(\bar{z}^{k})^*p,
\end{equation}
where $ e(z):=(e^{x}_{0}(z),\ldots,e^{x}_{k-1}(z),e^{y}_{0}(z),\ldots,e^{y}_{k-1}(z),e^{u}_{0}(z),\ldots,e^{u}_{k-1}(z),e^{a}_{0}(z),\ldots,e^{a}_{k-1}(z),
e^{b}_{0}(z),\ldots,e^{b}_{k-1}(z),\\d_{0}(z),\ldots,d_{k-1}(z))\in \mathbb{R}^{k(4n+m+d)}$ and $ p=(p_{1},\ldots,p_{k})\in \mathbb{R}^{k(4n+m+d)} $. The structures of the constraint sets $\Xi_j$ provide an equivalent expressions of the inclusions in \eqref{e:6.12} as
\begin{equation}\label{e:6.15}
\big(x^{\ast}_{jj},y^{\ast}_{jj},u^{\ast}_{jj},a^{\ast}_{jj},b^{\ast}_{jj},-X^{\ast}_{jj}\big)\in N_{{\rm\small gph}\,F_{h}}\Big(\bar{x}_{j},\bar{y}_{j},\bar{u}_{j},\bar{a}_{j},\bar{b}_{j},\dfrac{\bar{x}_{j+1}-\bar{x}_{j}}{-h_{k}}\Big),\,\,\,j=0,\ldots,k-1,
\end{equation}
while the other components of $ z^{\ast}_{j} $ are zero. Similarly we have that the vectors $ (x^{\ast}_{0k},y^{\ast}_{0k},u^{\ast}_{0k},a^{\ast}_{0k},b^{\ast}_{0k}) $ determined by the normal cone to $ \Xi_{k} $ might be the only nonzero components of $ z^{\ast}_{k} $. This gives us the representation
\begin{equation*}
\begin{aligned}
-\sum\limits_{j=0}^{k}z^{\ast}_{j}=&(-x^{\ast}_{00}-x^{\ast}_{0k},-x^{\ast}_{11},\ldots,-x^{\ast}_{(k-1)(k-1)},0,-y^{\ast}_{00}-y^{\ast}_{0k},-y^{\ast}_{11},\ldots,-y^{\ast}_{(k-1)(k-1)},0,\\
&-u^{\ast}_{00}-u^{\ast}_{0k},-u^{\ast}_{11},\ldots,-u^{\ast}_{(k-1)(k-1)},0,-a^{\ast}_{00}-a^{\ast}_{0k},-a^{\ast}_{11},\ldots,-a^{\ast}_{(k-1)(k-1)},0,\\
&-b^{\ast}_{00}-b^{\ast}_{0k},-b^{\ast}_{11},\ldots,-b^{\ast}_{(k-1)(k-1)},0,-X^{\ast}_{00},\ldots,-X^{\ast}_{(k-1)(k-1)},0,\ldots,0).
\end{aligned}
\end{equation*}
For the other terms in \eqref{e:6.13} we get
\begin{equation*}
\begin{aligned}
\sum\limits_{i=1}^{s}\alpha^{k}_{i}\nabla c_{i}(\bar{z}^{k})=(0,\ldots,-\sum\limits_{i=1}^{d}\alpha_{i}\nabla_{x} g_{i}(\bar{x}^{k}-\bar{u}^{k}),0,\ldots,0,0,\ldots,\sum\limits_{i=1}^{d}\alpha_{i}\nabla_{u} g_{i}(\bar{x}^{k}-\bar{u}^{k}_{k}),0,\ldots,0,0\ldots,0),
\end{aligned}
\end{equation*}
\begin{equation*}
\begin{aligned}
\nabla e(\bar{z}^{k})^*p=&\big(-p^{x}_{1}-\xi_{0}^{x}p^{d}_{1},p^{x}_{1}-p^{x}_{2}-\xi_{1}^{x}p^{d}_{2},\ldots,p^{x}_{k-1}-p^{x}_{k}-\xi_{k-1}^{x}p^{d}_{k},p^{x}_{k},\\
&-p^{y}_{1},p^{y}_{1}-p^{y}_{2},\ldots,p^{y}_{k-1}-p^{y}_{k},p^{y}_{k},-p^{u}_{1},p^{u}_{1}-p^{u}_{2},\ldots,p^{u}_{k-1}-p^{u}_{k},p^{u}_{k},\\
&-p^{a}_{1},p^{a}_{1}-p^{a}_{2},\ldots,p^{a}_{k-1}-p^{a}_{k},p^{a}_{k},\\
&-p^{b}_{1}-\xi_{0}^{b}p^{d}_{1},p^{b}_{1}-p^{b}_{2}-\xi_{1}^{b}p^{d}_{2},\ldots,p^{b}_{k-1}-p^{b}_{k}-\xi_{k-1}^{b}p^{d}_{k},p^{b}_{k},\\
&-h_{k}p^{x}_{1},\ldots,-h_{k}p^{x}_{k},-h_{k}p^{y}_{1}+p_{1}^{d},\ldots,-h_{k}p^{y}_{k}+p^{d}_{k},-h_{k}p^{u}_{1},\ldots,-h_{k}p^{u}_{k},\\
&-h_{k}p^{a}_{1},\ldots,-h_{k}p^{a}_{k},-h_{k}p^{b}_{1},\ldots,-h_{k}p^{b}_{k}\big).
\end{aligned}
\end{equation*}
Using the summation structure of the function $\ph_0$ and applying there the subdifferential sum rule from \cite[Theorem~2.19]{b3} give us the inclusion
\begin{equation*}
\begin{aligned}
\partial \varphi_{0}(\bar{z}^{k})\subset\partial \varphi(\bar{x}^{k}_{k})+h_{k}\sum\limits_{j=0}^{k-1}\partial l\big(\bar{x}_{j},\bar{y}_{j},\bar{u}_{j},\bar{a}_{j},\bar{b}_{j},\bar{X}_{j},\bar{Y}_{j},\bar{U}_{j},\bar{A}_{j},\bar{B}_{j}\big)+\sum\limits_{j=0}^{k-1}\nabla \rho_{j}(\bar{z}^{k}),
\end{aligned}
\end{equation*}
where $ \rho_{j}(y):=\dfrac{1}{2}\displaystyle\int\limits_{t_{j}}^{t_{j+1}}\lVert (X_{j},Y_{j},U_{j},A_{j},B_{j})-\dot{\bar{z}}(t) \rVert^{2}\,dt$. It is easy to see that
\begin{equation*}
\begin{aligned}
&\nabla \rho_{j}(\bar{z}^{k})=(0,\ldots,0,0,\ldots,0,0,\ldots,0,0,\ldots,0,0\ldots,0,0,\ldots,0,\theta^{x}_{j},0,\ldots,0,0\ldots,0,\theta^{y}_{j},0,\ldots,0,0,\ldots,0,\theta^{u}_{j},0,\ldots,0,\\
& 0,\ldots,0,\theta^{a}_{j},0,\ldots,0,0,\ldots,0,\theta^{b}_{j},0,\ldots,0),
\end{aligned}
\end{equation*}
\begin{equation*}
\begin{aligned}
&\partial l(\bar{x}_{j},\bar{y}_{j},\bar{u}_{j},\bar{a}_{j},\bar{b}_{j},\bar{X}_{j},\bar{Y}_{j},\bar{U}_{j},\bar{A}_{j},\bar{B}_{j})
=(0,\ldots,0,w^{x}_{j},0,\ldots,0,0,\ldots,0,0,\ldots,0,w^{u}_{j},0,\ldots,0,0,\ldots,0,w^{a}_{j},0,\ldots,0,0,\ldots,\\&0,
w^{b}_{j},0,\ldots,0,0,\ldots,0,v^{x}_{j},0,\ldots,0,0,\ldots,0,0,
\ldots,0,v^{u}_{j},0,\ldots,0,0,\ldots,0,v^{a}_{j},0,\ldots,0,0,\ldots,0,v^{b}_{j},0,\ldots,0).
\end{aligned}
\end{equation*}
Since $\partial \varphi(\bar{x}^{k}_{k})$ can be written in the form $(0,\ldots,0,\upsilon^{k},0,\ldots,0,0,\ldots,0,0,\ldots,0,0,\ldots,0)$, we arrive at the representation of the term $ \lambda^{k}\partial\varphi_{0}(\bar{z}^{k}) $ in \eqref{e:6.13} as
\begin{equation*}
\lambda^{k}(h_{k}w^{x}_{0},\ldots,h_{k}w^{x}_{k-1},\upsilon^{k},0,\ldots,0,h_{k}w^{u}_{0},\ldots,h_{k}w^{u}_{k-1},0,
h_{k}w^{a}_{0},\ldots,h_{k}w^{a}_{k-1},0,h_{k}w^{b}_{0},\ldots,h_{k}w^{b}_{k-1},0,
\end{equation*}
\begin{equation*}
\theta^{x}_{0}+h_{k}v^{x}_{0},\ldots,\theta^{x}_{k-1}+h_{k}v^{x}_{k-1},\theta^{y}_{0},\ldots,\theta^{y}_{k-1},\theta^{u}_{0}+h_{k}v^{u}_{0},\ldots,\theta^{u}_{k-1}+h_{k}v^{u}_{k-1},\theta^{a}_{0}+h_{k}v^{a}_{0},\ldots,\theta^{a}_{k-1}+h_{k}v^{a}_{k-1},
\end{equation*}
\begin{equation*}
\theta^{b}_{0}+h_{k}v^{b}_{0},\ldots,\theta^{b}_{k-1}+h_{k}v^{b}_{k-1}) .
\end{equation*}
Combining all the above together decomposes the inclusion in \eqref{e:6.13} into the following equalities:
\begin{equation*}
-x^{\ast}_{00}-x^{\ast}_{0k}=\lambda^{k}h_{k}w^{x}_{0}-p^{x}_{1} -\xi_{0}^{x}p^{d}_{1},
\end{equation*}
\begin{equation}\label{e:6.16}
-x^{\ast}_{jj}=\lambda^{k}h_{k}w^{x}_{j}+p^{x}_{j}-p^{x}_{j+1}-\xi_{j}^{x}p^{d}_{j+1},\,\,\,j=1,\ldots,k-1,
\end{equation}
\begin{equation}\label{e:6.17}
0=\lambda v^{k}-\sum\limits_{i=1}^{s}\alpha^{k}_{i}\nabla_{x}g_{i}(\bar{x}^{k}-\bar{u}^{k})+p^{x}_{k},
\end{equation}
\begin{equation}\label{e:6.19}
-y^{\ast}_{00}=-p^{y}_{1} +p^{d}_{1},\quad-y^{\ast}_{jj}=   p^{y}_{j}-p^{y}_{j+1}+p^{d}_{j+1},\,\,\,j=1,\ldots,k-1,\quad p^{y}_{k}=0,
\end{equation}
\begin{equation}\label{e:6.20}
-u^{\ast}_{00}-u^{\ast}_{0k}=\lambda^{k}h_{k}w^{u}_{0}-p^{u}_{1},
\end{equation}
\begin{equation}\label{e:6.21}
-u^{\ast}_{jj}=\lambda^{k}h_{k}w^{u}_{j}+p^{u}_{j}-p^{u}_{j+1},\,\,\,j=1,\ldots,k-1,
\end{equation}
\begin{equation}\label{e:6.22}
0=\sum\limits_{i=1}^{s}\alpha^{k}_{i}\nabla_{u}g_{i}(\bar{x}^{k}-\bar{u}^{k})+p^{u}_{k},
\end{equation}
\begin{equation}
-a^{\ast}_{00}-a^{\ast}_{0k}=\lambda^{k}h_{k}w^{a}_{0}-p^{a}_{1},
\end{equation}
\begin{equation}\label{e:6.23}
-a^{\ast}_{jj}=\lambda^{k}h_{k}w^{a}_{j}+p^{a}_{j}-p^{a}_{j+1},\,\,\,j=1,\ldots,k-1,\quad 0=p^{a}_{k},
\end{equation}
\begin{equation*}
-b^{\ast}_{00}-b^{\ast}_{0k}=\lambda^{k}h_{k}w^{b}_{0}-p^{b}_{1}-\xi_{j}^{b}p^{d}_{1},
\end{equation*}
\begin{equation}\label{e:6.24}
-b^{\ast}_{jj}=\lambda^{k}h_{k}w^{b}_{j}+p^{b}_{j}-p^{b}_{j+1}-\xi_{j}^{b}p^{d}_{j+1},\,\,\,j=1,\ldots,k-1,\quad 0=p^{b}_{k},
\end{equation}
\begin{equation}\label{e:6.25}
-X^{\ast}_{jj}=\lambda^{k}(h_{k}v^{x}_{j}+\theta^{x}_{j})-h_{k}p^{x}_{j+1},\,\,\,j=0,\ldots,k-1,
\end{equation}
\begin{equation}\label{e:6.26}
0=\lambda^{k}\theta^{y}_{j}-h_{k}p^{y}_{j+1}+p^{d}_{j+1},\,\,\,j=0,\ldots,k-1,
\end{equation}
\begin{equation}\label{e:6.27}
0=\lambda^{k}(h_{k}v^{u}_{j}+\theta^{u}_{j})-h_{k}p^{u}_{j+1},\,\,\,j=0,\ldots,k-1 ,
\end{equation}
\begin{equation}\label{e:6.28}
0=\lambda^{k}(h_{k}v^{a}_{j}+\theta^{a}_{j})-h_{k}p^{a}_{j+1},\,\,\,j=0,\ldots,k-1,
\end{equation}
\begin{equation}\label{e:6.29}
0=\lambda^{k}(h_{k}v^{b}_{j}+\theta^{b}_{j})-h_{k}p^{b}_{j+1},\,\,\,j=0,\ldots,k-1.
\end{equation}
We deduce from  \eqref{e:6.16}, \eqref{e:6.19}, \eqref{e:6.21}, \eqref{e:6.23}, \eqref{e:6.24}, and \eqref{e:6.25} that
\begin{equation*}
\begin{array}{ll}
\Big(\dfrac{p^{x}_{j+1}-p^{x}_{j}}{h_{k}}-\lambda^{k}w^{x}_{j}+h_{k}^{-1}\xi_{j}^{x}p^{d}_{j+1},\dfrac{p^{y}_{j+1}-p^{y}_{j}}{h_{k}}-\dfrac{p^{d}_{j+1}}{h_{k}},\dfrac{p^{u}_{j+1}-p^{u}_{j}}{h_{k}}
-\lambda^{k}w^{u}_{j},\dfrac{p^{a}_{j+1}-p^{a}_{j}}{h_{k}}-\lambda^{k}w^{a}_{j},\\
\dfrac{p^{b}_{j+1}-p^{b}_{j}}{h_{k}}-\lambda^{k}w^{b}_{j}+\xi_{j}^{b}\dfrac{p^{d}_{j+1}}{h_{k}},
p^{x}_{j+1}-\lambda^{k}(v^{x}_{j}+h_{k}^{-1}\theta^{x}_{j})\Big)=h_{k}^{-1}(x^{\ast}_{jj},y^{\ast}_{jj},u^{k\ast}_{jj},a^{\ast}_{jj},b^{\ast}_{jj},-X^{\ast}_{jj}).
\end{array}
\end{equation*}
Then the conditions in \eqref{e:6.11}  follow from \eqref{e:6.15}, and the conditions in \eqref{e:6.3}, \eqref{e:6.4}, \eqref{e:6.7}, \eqref{e:6.9} follow from \eqref{e:6.14}, \eqref{e:6.17}, \eqref{e:6.22}, \eqref{e:6.26}, \eqref{e:6.27}, \eqref{e:6.28}, and \eqref{e:6.29}, respectively.\vspace*{-0.05in}

To verify finally the nontriviality condition \eqref{e:6.2} of the theorem, denote $ p_{0}:=(x^{\ast}_{0k},y^{\ast}_{0k},u^{\ast}_{0k},a^{\ast}_{0k},b^{\ast}_{0k}) $ and suppose by contraposition that  $ \lambda^{k}=0 $, $ \alpha^{k}=0 $, $ p^{y}_{0}=0 $, $ p^{u}_{0}=0 $, $ p^{a}_{0}=0 $, $ p^{b}_{0}=0 $, $ p^{x}_{j}=0 $ for $ j=0,\ldots,k-1 $, and $ p^{d}_{j}=0 $ for $ j=0,\ldots,k $. Then we get the implications
\begin{equation*}
\begin{aligned}
&\eqref{e:6.17}\Longrightarrow p^{x}_{k}=0,\;\mbox{ i.e., }\;p^{x}_{j}=0,\,\,\, j=0,\ldots,k,\\
&\eqref{e:6.16}\Longrightarrow x^{\ast}_{jj}=0\;\mbox{ and }\;\eqref{e:6.25}\Longrightarrow X^{\ast}_{jj}=0,\,\,\,j=0,\ldots,k-1,\\
&\eqref{e:6.27}, \eqref{e:6.28}, \eqref{e:6.29}, \eqref{e:6.26}\Longrightarrow p^{u}_{j}=0, p^{a}_{j}=0, p^{b}_{j}=0, p^{y}_{j}=0 ,\,\,\,j=0,\ldots,k,\\
&\eqref{e:6.19}, \eqref{e:6.21}, \eqref{e:6.23}, \eqref{e:6.24}\Longrightarrow y^{\ast}_{jj}=0, u^{\ast}_{jj}=0,  a^{\ast}_{jj}=0,  b^{\ast}_{jj}=0,\,\,\, j=0,\ldots,k-1.
\end{aligned}
\end{equation*}
We know that all the components of $ z^{\ast}_{j} $ different from $ (x^{\ast}_{jj},y^{\ast}_{jj},u^{\ast}_{jj},a^{\ast}_{jj},b^{\ast}_{jj},X^{\ast}_{jj}) $ are zero for $ j=0,\ldots,k-1  $. Hence  $ z^{\ast}_{j}=0 $ for $ j=0,\ldots,k-1 $. Conclude similarly that  $ z^{\ast}_{k}=0 $ due to
$$
(x^{\ast}_{0k},y^{\ast}_{0k},u^{\ast}_{0k},a^{\ast}_{0k},b^{\ast}_{0k})=(p^{x}_{0},p^{y}_{0},p^{u}_{0},p^{a}_{0},p^{b}_{0},p^{d}_{0})=(0,0,0,0,0,0),
$$
and thus $ z^{\ast}_{j}=0 $ for all $ j=0,\ldots,k$. This contradicts the nontriviality condition in the mathematical program (MP), and hence verifies the claimed nontriviality \eqref{e:6.2} while completing the proof of the theorem. $\h$

The next theorem is the main result of this section. It provides necessary optimality conditions for problems $(P_k)$ as formulated in Section~\ref{sec:disc-app} expressed {\em entirely in terms of the given data}. We derive this result by combining necessary optimal conditions of Theorem~\ref{discr:N.C} with the coderivative calculations of Theorem~\ref{Th:co-nc} addressing the specific form of the velocity mapping $F_h$ in $(P_k)$. In this way we obtain a novel {\em discrete-time Volterra condition} as a part of primal-dual dynamic relationships for discrete approximations of controlled integro-differential sweeping processes. \vspace*{-0.1in}
\begin{theorem}[necessary conditions for discretized integro-differential sweeping control systems]\label{Thm:O.c.d} Let $\bar{z}^{k}=(\bar{x}^{k},\bar{y}^{k},\bar{u}^{k},\bar{a}^{k},\bar{b}^{k})$ be an optimal solution to the discrete-time problem $(P_{k})$, where $k\in\N$ is sufficiently large. In addition to the assumptions of Theorem~{\rm\ref{discr:N.C}}, suppose that the functions $g_{i}$, $i=1,\ldots,s$, are twice continuously differentiable around $\bar{x}^{k}_{j}-\bar{u}^{k}_{j}$, $j=0,\ldots,k-1$, with the Jacobian matrix of full rank therein. Then there exist dual elements  $(\lambda^{k},p^{k})$ as in Theorem~{\rm\ref{discr:N.C}} together with vectors $\eta_{j}^{k}\in \mathbb{R}^{s}_{+}$ as $j=0,\ldots,k$ and $\sigma ^{k}_{j}\in \mathbb{R}^{s} $ as $j=0,\ldots,k-1$ satisfying following conditions:\\[1ex]
{\sc Primal-Dual Dynamic Relationships}: for all $j=0,\ldots,k-1$ we have
\begin{equation}\label{e:6.32}
\dfrac{\bar{x}^{k}_{j+1}-\bar{x}^{k}_{j}}{h_{k}}+f_{1}(\bar{a}^{k}_{j},\bar{x}^{k}_{j})+\bar{y}^{k}_{j}+h_{k}f_{2}(\bar{b}^{k}_{j},\bar{x}^{k}_{j})=\sum\limits_{i\in I(\bar{x}^{k}_{j}-\bar{u}^{k}_{j})} \eta^{k}_{ji}\nabla g_{i}(\bar{x}^{k}_{j}-\bar{u}^{k}_{j}),
\end{equation}
\begin{equation}\label{e:6.33}
\begin{aligned}
\dfrac{p^{x}_{j+1}-p^{x}_{j}}{h_{k}}-\lambda^{k}w^{x}_{j}+h_{k}^{-1}\xi_{j}^{x}p^{d}_{j+1}= \nabla_{x} f_{1}(\bar{a}^{k}_{j},\bar{x}^{k}_{j})^{*}(\lambda^{k}(v^{x}_{j}+h_{k}^{-1}\theta^{x}_{j})-p^{x}_{j+1})\\
+h_{k}\nabla_{x} f_{2}(\bar{b}^{k}_{j},\bar{x}^{k}_{j})^{*}(\lambda^{k}(v^{x}_{j}+h_{k}^{-1}\theta^{x}_{j})-p^{x}_{j+1})
\end{aligned}
\end{equation}
\begin{equation*}
-\Big(\sum\limits_{i=1}^{d}\eta^{k}_{ji}\nabla_{x}^{2}g_{i}(\bar{x}^{k}_{j}-\bar{u}^{k}_{j})\Big)(\lambda^{k}(v^{x}_{j}+h_{k}^{-1}\theta^{x}_{j})-p^{x}_{j+1})
-\sum\limits_{i=1}^{d}\sigma^{k}_{ji}\nabla_{x}g_{i}(\bar{x}^{k}_{j}-\bar{u}^{k}_{j})),
\end{equation*}
\begin{equation}\label{e:6.34}
\dfrac{p^{y}_{j+1}-p^{y}_{j}}{h_{k}}-\dfrac{p^{d}_{j+1}}{h_{k}}=(\lambda^{k}(v^{x}_{j}+h_{k}^{-1}\theta^{x}_{j})-p^{x}_{j+1}),
\end{equation}
\begin{equation}\label{e:6.35}
\dfrac{p^{u}_{j+1}-p^{u}_{j}}{h_{k}}-\lambda^{k}w^{u}_{j}= \Big(\sum\limits_{i=1}^{d}\eta^{k}_{ji}\nabla_{u}^{2}g_{i}(\bar{x}^{k}_{j}-\bar{u}^{k}_{j})\Big)(\lambda^{k}(v^{x}_{j}+h_{k}^{-1}\theta^{x}_{j})-p^{x}_{j+1})+\sum\limits_{i=1}^{d}\sigma^{k}_{ji}\nabla_{u}g_{i}(\bar{x}^{k}_{j}-\bar{u}^{k}_{j})),
\end{equation}
\begin{equation}\label{e:6.36}
\dfrac{p^{a}_{j+1}-p^{a}_{j}}{h_{k}}-\lambda^{k}w^{a}_{j}= \nabla_{a} f_{1}(\bar{a}^{k}_{j},\bar{x}^{k}_{j})^{*}(\lambda^{k}(v^{x}_{j}+h_{k}^{-1}\theta^{x}_{j})-p^{x}_{j+1}),
\end{equation}
\begin{equation}\label{e:6.37}
\dfrac{p^{b}_{j+1}-p^{b}_{j}}{h_{k}}-\lambda^{k}w^{b}_{j}+h_{k}^{-1}\xi_{j}^{b}p^{d}_{j+1}=h_{k}\nabla_{b} f_{2}(\bar{b}^{k}_{j},\bar{x}^{k}_{j})^{*}(\lambda^{k}(v^{x}_{j}+h_{k}^{-1}\theta^{x}_{j})-p^{x}_{j+1}).
\end{equation}
{\sc Transversality Conditions}:
\begin{equation}\label{e:6.38}
-p^{x}_{k}\in\lambda^{k}\partial\varphi(\bar{x}^{k}_{k})-\sum\limits_{i=1}^{s}\eta^{k}_{ki}\nabla_{x}g_{i}(\bar{x}^{k}-\bar{u}^{k}),
\end{equation}
\begin{equation}\label{e:6.39}
-p^{uk}_{k}=\sum\limits_{i=1}^{s}\eta^{k}_{ki}\nabla_{u} g_{i}(\bar{x}^{k}_{k}-\bar{u}^{k}_{k}),\,\,\,p^{yk}_{k}=0,\,\,\,p^{ak}_{k}=0,\,\,\,p^{bk}_{k}=0 .
\end{equation}
{\sc Complementary Slackness}: for all $j=0,\ldots,k-1$ and $i=1,\ldots,s$ we have the implications
\begin{equation}\label{e:6.41}
g_{i}(\bar{x}^{k}_{j}-\bar{u}^{k}_{j})>0 \Longrightarrow\eta^{k}_{ji}=0,
\end{equation}
\begin{equation}\label{e:6.42}
\big[g_{i}(\bar{x}^{k}_{j}-\bar{u}^{k}_{j})>0,\;\mbox{ or }\;\eta^{k}_{ji}=0\;\mbox{ and }\;\langle \nabla\,g_{i}(\bar{x}^{k}_{j}-\bar{u}^{k}_{j}),\lambda^{k}(v^{x}_{j}+h_{k}^{-1}\theta^{x}_{j})-p^{x}_{j+1} \rangle>0\Big]\Longrightarrow\sigma^{k}_{ji}=0,
\end{equation}
\begin{equation}\label{e:6.43}
\big[g_{i}(\bar{x}^{k}_{j}-\bar{u}^{k}_{j})=0,\,\,\,\eta^{k}_{ji}=0,\,\,\,\langle \nabla\,g_{i}(\bar{x}^{k}_{j}-\bar{u}^{k}_{j}),\lambda^{k}(v^{x}_{j}+h_{k}^{-1}\theta^{x}_{j})-p^{x}_{j+1} \rangle<0\big]\Longrightarrow\sigma^{k}_{ji}\geq 0,
\end{equation}
\begin{equation}\label{e:6.44}
g_{i}(\bar{x}^{k}_{k}-\bar{u}^{k}_{k})>0\Longrightarrow \eta^{k}_{ki}= 0,
\end{equation}
\begin{equation}\label{e:6.45}
\eta^{k}_{ji}>0\Longrightarrow \langle \nabla\,g_{i}(\bar{x}^{k}_{j}-\bar{u}^{k}_{j}),\lambda^{k}(v^{x}_{j}+h_{k}^{-1}\theta^{x}_{j})-p^{x}_{j+1} \rangle=0.
\end{equation}
{\sc Nontriviality Condition}:
\begin{equation}\label{e:6.31}
\lambda^{k}+\sum\limits_{j=0}^{k}\lVert p^{kd}_{j} \rVert  + \lVert p^{ku}_{0} \rVert \neq 0,
\end{equation}
where the vectors $\xi^x_j$ and $\xi^b_j$ are taken from \eqref{xi}.
\end{theorem}\vspace*{-0.1in}
{\bf Proof}. Using the coderivative construction \eqref{cod}, for all $j=0,\ldots,k-1$ we rewrite the discrete Euler-Lagrange inclusions \eqref{e:6.11} of Theorem~\ref{discr:N.C} in the form
\begin{equation}\label{e:6.46}
\begin{array}{ll}
&\Big(\dfrac{p^{x}_{j+1}-p^{x}_{j}}{h_{k}}-\lambda^{k}w^{x}_{j}+h_{k}^{-1}\xi_{j}^{x}p^{d}_{j+1},\dfrac{p^{y}_{j+1}-p^{y}_{j}}{h_{k}}-\dfrac{p^{d}_{j+1}}{h_{k}},\dfrac{p^{u}_{j+1}-p^{u}_{j}}{h_{k}}-\lambda^{k}w^{u}_{j},\\&
\dfrac{p^{a}_{j+1}-p^{a}_{j}}{h_{k}}-\lambda^{k}w^{a}_{j},\dfrac{p^{b}_{j+1}-p^{b}_{j}}{h_{k}}-\lambda^{k}w^{b}_{j}+h_{k}^{-1}\xi_{j}^{b}p^{d}_{j+1}, p^{x}_{j+1}-
\lambda^{k}(v^{x}_{j}+h_{k}^{-1}\theta^{x}_{j})\Big)\\
&\in D^{*}F_{h}\Big(\bar{x}^{k}_{j},\bar{y}^{k}_{j},\bar{u}^{k}_{j},\bar{a}^{k}_{j},\bar{b}^{k}_{j},\dfrac{\bar{x}^{k}_{j+1}-\bar{x}^{k}_{j}}{-h_{k}}\Big)(\lambda^{k}(v^{x}_{j}+h_{k}^{-1}\theta^{x}_{j})
-p^{x}_{j+1}).
\end{array}
\end{equation}
It follows from the discrete dynamics \eqref{e:5.0}, representation \eqref{F1} of the velocity mapping, and the structure of the moving set in \eqref{e:mset} that there exist vectors $\eta^{k}_{j}\in \mathbb{R}^{m}_{+}$ as $j=0,\ldots,k-1$ such that the conditions in \eqref{e:6.32} and \eqref{e:6.41}  are satisfied. Employing further the coderivative calculation in \eqref{e:4.2} and \eqref{e:4.3} of Theorem~\ref{Th:co-nc} with $x:=\bar{x}^{k}_{j}$,  $y:=\bar{y}^{k}_{j}$, $u:=\bar{u}^{k}_{j}$, $a:=\bar{a}^{k}_{j}$, $b:=\bar{b}^{k}_{j}$, $w:=\dfrac{\bar{x}^{k}_{j+1}-\bar{x}^{k}_{j}}{-h_{k}}$, and $z=(\lambda^{k}(v^{x}_{j}+h_{k}^{-1}\theta^{x}_{j})-p^{x}_{j+1})$ allows us to find $\sigma^{k}_{j}$, $j=0,\ldots,k-1$, for which conditions \eqref{e:6.33}, \eqref{e:6.34},  \eqref{e:6.35}, \eqref{e:6.36}, \eqref{e:6.37}, \eqref{e:6.42}, and \eqref{e:6.43} hold. Define $\eta^{k}_{k}:=\alpha^{k}$ and observe that $\eta^{k}_{j}\in \mathbb{R}^{m}_{+}$ for all $j=0,\ldots,k$. In this way
we deduce the transversality conditions \eqref{e:6.38} and \eqref{e:6.39} from \eqref{e:6.4} and \eqref{e:6.6}, while \eqref{e:6.44} follows from \eqref{e:6.14} and the definition of $\eta^{k}_{k}$.
Note also that \eqref{e:6.45} follows directly from \eqref{e:4.3}.\vspace*{-0.05in}

It remains to verify the fulfillment of the nontriviality condition \eqref{e:6.31}. To this end, deduce first from \eqref{e:6.2} and the constructions above that
\begin{equation}\label{e:6.30}
\lambda^{k}+\lVert \eta^{k}_{k} \rVert+\sum\limits_{j=0}^{k-1}\lVert p^{kx}_{j} \rVert +\sum\limits_{j=0}^{k}\lVert p^{kd}_{j} \rVert +\lVert p^{ky}_{0} \rVert + \lVert p^{ku}_{0} \rVert+\lVert p^{ka}_{0} \rVert+\lVert p^{kb}_{0} \rVert \neq 0 .
\end{equation}
Assume now that \eqref{e:6.13} is violated, i.e., $\lambda^{k}=0$, $p^{kd}_{j}=0$ for  all $j=0,\ldots,k$, and $p^{ku}_{0}=0$. Then it follows from \eqref{e:6.7} and \eqref{e:6.9} that $p^{uk}_{j}=0$, $j=0,\ldots,k$, $p^{yk}_{j}=0$, $p^{ak}_{j}=0$, and $p^{bk}_{j}=0$ for all $j=1,\ldots,k$. Furthermore, \eqref{e:6.4} yields $\sum\limits_{i=1}^{s}\eta^{k}_{ki}\nabla_{u} g_{i}(\bar{x}^{k}_{k}-\bar{u}^{k}_{k}) =0$, and hence $p^{xk}_{k}=0$. Employing \eqref{e:6.35} tells us that
\begin{equation*}
\Big(\sum\limits_{i=1}^{d}\eta^{k}_{ji}\nabla_{u}^{2}g_{i}(\bar{x}^{k}_{j}-\bar{u}^{k}_{j})\Big)(\lambda^{k}(v^{x}_{j}+h_{k}^{-1}\theta^{x}_{j})-p^{x}_{j+1})+\sum\limits_{i=1}^{d}
\sigma^{k}_{ji}\nabla_{u}g_{i}(\bar{x}^{k}_{j}-\bar{u}^{k}_{j}))=0,\quad j=0,\ldots,k-1.
\end{equation*}
Combining the latter with \eqref{e:6.33} and $p^{xk}_{k}=0$ ensures that $p^{xk}_{j}=0$ whenever $j=0,\ldots,k-1$. To complete the proof of the theorem, we deduce from \eqref{e:6.34}, \eqref{e:6.36}, and \eqref{e:6.37} that $p^{yk}_{0}=0$, $p^{ak}_{0}=0$, and $p^{bk}_{0}=0$. This contradicts the fulfillment of \eqref{e:6.30} and hence justifies \eqref{e:6.31}. $\h$\vspace*{-0.2in}
%%%%%%%%%%%%%%%%%%%%%%%%%%%%%%%%%%%%%%%%%%%%%%%%%%%%%%%%%%%%%%%%%%%%%%%%%%%%%%%%%%%%%%%%%%%%%%%%%%%%%%%%%%%%%%%%%%%%%%%%%%%%%%%%%%%%%%%%%%%%%%%%%%%%%%%%%%%%%%%%%%%%%%%%%%%%%%%%%%%%%%%%%%%%%%%
%%%%%%%%%%%%%%%%%%%%%%%%%%%%%%%%%%%%%%%%%%%%%%%%%%%%%%%%%%%%%%%%%%%%%%%%%%%%%%%%%%%%%%%%%%%%%%%%%%%%%%%%%%%%%%%%%%%%%%%%%%%%%%%%%%%%%%%%%%%%%%%%%%%%%%%%%

\section{Necessary Conditions for Integro-Differential Processes}\label{sec:nec-sw}
\setcounter{equation}{0}\vspace*{-0.1in}
This section establishes the main result of the paper providing---for the first time in the literature---efficient necessary optimality conditions, expressed entirely via the given data, for local minimizers (in the sense of Definition~\ref{locmin} (ii)), of the original optimal control problem $(P)$ governed by the sweeping integro-differential inclusions \eqref{e:1} with the pointwise mixed state-control constraints \eqref{mixed}. The derivation of these conditions presented below is based on the results obtained in the previous sections as well as on the appropriate properties of the generalized differential constructions of variational analysis reviewed in Section~\ref{sec:2va} that allow us furnishing the passage to the limit as $k\to\infty$ from the necessary optimality conditions for discrete approximations. For simplicity, we suppose below that the running cost in \eqref{e:bolza}  does not depend on $t$.\vspace*{-0.1in}
\begin{theorem}[\bf optimality conditions for integro-differential sweeping processes]\label{Thm:N.C.I} Consider an relaxed intermediate local minimizer $\bar{z}(\cdot)=(\bar{x}(\cdot),\bar{y}(\cdot),\bar{u}(\cdot),\bar{a}(\cdot),\bar{b}(\cdot))$ of problem $(P)$ and suppose in addition to the assumptions of Theorem~{\rm\ref{Thm:O.c.d}} that the running cost in \eqref{e:bolza} is represented as
\begin{equation}\label{run-cost}
l_{0}(\bar{x},\bar{u},\bar{a},\bar{b},\dot{\bar{x}},\dot{\bar{u}},\dot{\bar{a}},\dot{\bar{b}}):=l(z,\dot{z})=l_{1}(z,\dot{x})+l_{2}(\dot{u})+l_{3}(\dot{a})+l_{4}(\dot{b}),
\end{equation}
where $l_{2}$ is differentiable on $\mathbb{R}^{n}$ with the estimates
\begin{equation*}
\lVert \nabla_{\dot{u}} l_{2}(\dot{u})\rVert\leq L \lVert \dot{u} \rVert\;\mbox{ and }\;\lVert \nabla_{\dot{u}} l_{2}(\dot{u}_{1})-\nabla_{\dot{u}} l_{2}(\dot{u}_{2})\rVert\leq L \lVert \dot{u}_{1}-\dot{u}_{2} \rVert.
\end{equation*}
Then there exist a multiplier $\lambda\geq 0$ and quintuples $p(\cdot)=(p^{x}(\cdot),p^{y}(\cdot),p^{u}(\cdot),p^{a}(\cdot),p^{b}(\cdot))\in W^{1,2}([0,T],\mathbb{R}^{3n+m+d})$, $w(\cdot)=(w^{x}(\cdot),0,w^{u}(\cdot),w^{a}(\cdot),w^{b}(\cdot))\in L^{2}([0,T],\mathbb{R}^{3n+m+d})$, and $v(\cdot)=(v^{x}(\cdot),0,v^{u}(\cdot),v^{a}(\cdot),v^{b}(\cdot))\\\in L^{2}([0,T],\mathbb{R}^{3n+m+d})$ satisfying the subdifferential inclusion
\begin{equation}\label{e:7.1}
(w(t),v(t))\in{\rm co}\,\partial l(\bar{z}(t), \dot{\bar{z}}(t))\;\mbox{ a.e. }\;t\in [0,T ],
\end{equation}
as well as measure $\gg=(\gg_{1},\ldots,\gg_{n})\in{\cal C}([0,T],\mathbb{R}^{n})^*$ for which the following conditions hold:\\
{\sc Primal-Dual Dynamic Relationships}:
\begin{equation}\label{e:7.2}
\dot{\bar{x}}(t)+f_{1}(\bar{a}(t),\bar{x}(t))+\bar{y}(t)=\sum\limits_{i=1}^{s}\eta_{i}(t)\nabla g_{i}(\bar{x}(t)-\bar{u}(t))\;\mbox{ a.e. }\;t\in [0,T];
\end{equation}
\begin{equation}\label{e:7.4}
g_{i}(\bar{x}(t)-\bar{u}(t))>0\Longrightarrow\eta_{i}(t)= 0\;\mbox{ a.e. }\;t\in [0,T];
\end{equation}
\begin{equation}\label{e:7.5}
\eta_{i}(t)>0 \Longrightarrow \langle \nabla\,g_{i}(\bar{x}(t)-\bar{u}(t)) , \lambda v^{x}(t)-q^{x}(t) \rangle =0\;\mbox{ a.e. }\;t\in [0,T],
\end{equation}
where the functions $\eta_i(\cdot)\in L^{2}([0,T],\mathbb{R}_{+})$ are uniquely determined by representation \eqref{e:7.2} for a.e.\ $t\in[0,T]$ while being well-defined at $t=T$;
\begin{equation}\label{e:7.6}
\dot{p}(t)=\lambda w(t)+\Big(\nabla_{x} f_{1}(\bar{a}(t),\bar{x}(t))^{*}(\lambda v^{x}(t)-q^{x}(t)),\lambda v^{x}(t)-q^{x}(t), 0, \nabla_{a} f_{1}(\bar{a}(t),\bar{x}(t))^{*}(\lambda v^{x}(t)-q^{x}(t)), 0 \Big);
\end{equation}
\begin{equation}\label{e:7.7}
q^{u}(t)=\lambda \nabla _{\dot{u}}l_{2}(\dot{\bar{u}}(t)),\,\,q^{a}(t)\in \lambda\, \partial_{\dot{a}}l_{3}(\dot{\bar{a}}(t)),\,\,q^{b}(t)\in \lambda\, \partial_{\dot{b}}l_{4}(\dot{\bar{b}}(t)),
\end{equation}
where $q(\cdot)=(q^{x}(\cdot),q^{y}(\cdot),q^{u}(\cdot),q^{a}(\cdot),q^{b}(\cdot)) : [0,T]\to \mathbb{R}^{3n+m+d} $ is a vector function of bounded variation, with $q^{y}(\cdot)$ being absolutely continuous on $[0,T]$, such that
left-continuous representative of $q(\cdot)$ satisfies, for all $t\in[0,T]$ except at most a countable subset, the integral equation
\begin{equation}\label{e:7.8}
q(t)=p(t)-\int\limits_{[t,T]}\Big(-d\gg(s)-\nabla_{x}f_{2}(\bar{b}(s),\bar{x}(s))q^{y}(s)\,ds,\;q^{y}(s)\,ds,\;d\gg(s),\;0,\;-\nabla_{b}f_{2}(\bar{b}(s),\bar{x}(s))q^{y}(s)\,ds\Big).
\end{equation}
{\sc Extended Volterra Condition}:
\begin{equation}\label{e:7.9}
\dot{q}^{y}(t)=\lambda v^{x}(t)-p^{x}(t)-\int\limits_{[t,T]}d\gg(s)-\int\limits_{[t,T]}\nabla_{x} f_{2}(\bar{b}(s),\bar{x}(s))q^{y}(s)\,ds+q^{y}(t)\;\mbox{ a.e. }\;t\in[0,T].
\end{equation}
{\sc Right Endpoint Conditions}:
\begin{equation}\label{e:7.10}
-p^{x}(T)+\sum\limits_{i\in I(\bar{x}(T)-\bar{u}(T))}\eta_{i}(T)\nabla g_{i}(\bar{x}(T)-\bar{u}(T))\in \lambda \partial \varphi  (\bar{x}(T)),\;p^{y}(T)=0,
\end{equation}
\begin{equation}\label{e:7.11}
-p^{u}(T)=\sum\limits_{i\in I(\bar{x}(T)-\bar{u}(T))}\eta_{i}(T)\nabla g_{i}(\bar{x}(T)-\bar{u}(T)),\;p^{a}(T)=0,\;p^{b}(T)=0,
\end{equation}
\begin{equation}\label{e:7.13}
-\sum\limits_{i\in I(\bar{x}(T)-\bar{u}(T))}\eta_{i}(T)\nabla g_{i}(\bar{x}(T)-\bar{u}(T))\in N_{C}(\bar{x}(T)-\bar{u}(T)).
\end{equation}
{\sc Nontriviality Condition}:
\begin{equation}\label{e:7.14}
\lambda +\lVert q^{u}(0) \rVert + \lVert p(T) \rVert + \int\limits_{0}^{T}\lVert q^{y}(t) \rVert\,dt>0.
\end{equation}
Furthermore, the following implications hold while ensuring the {\sc Enhanced nontriviality}:
\begin{equation}\label{ench:1}
[g_{i}(x_{0}-\bar{u}(0))>0,\;\;i=1,\ldots,s]\Longrightarrow [\lm +\lVert p(T) \rVert + \int\limits_{0}^{T}\lVert q^{y}(t) \rVert\,dt>0],
\end{equation}
\begin{equation}\label{ench:2}
[g_{i}(\bar{x}(T)-\bar{u}(T))>0,\;\;i=1,\ldots,s]\Longrightarrow [\lm +\lVert q^{u}(0) \rVert + \int\limits_{0}^{T}\lVert q^{y}(t) \rVert\,dt>0].
\end{equation}
\end{theorem}
{\bf Proof}. We proceed with passing to the limit as $k\to \infty$ in the necessary optimality conditions for discrete approximation problems obtained in Theorem~\ref{Thm:O.c.d} with taking into account the strong convergence of discrete optimal solutions established in Theorem~\ref{disc:conver}. Since some arguments in this procedure are similar to those used in \cite[Theorem~8.1]{b1} in a more special setting, we skip them for brevity while focusing on significantly new developments. Note, in particular, that the existence
of the subgradient functions $(w^{x}(\cdot),0,w^{u}(\cdot),w^{a}(\cdot),w^{b}(\cdot),v^{x}(\cdot),0,v^{u}(\cdot),v^{a}(\cdot),v^{b}(\cdot))$ satisfying \eqref{e:7.1} can be checked as in
\cite{b1}, while the existence of the uniquely defined elements $\eta_i(\cdot)\in L^{2}([0,T],\mathbb{R}_{+})$ satisfying \eqref{e:7.2} follows from representation \eqref{e:4.1} by repeating the limiting procedure of \cite[Theorem~8.1]{b1} with taking account the uniform convergence of $\{\bar{y}^{k}(t^{k}_{j})+h_{k}f_{2}(\bar{b}^{k}(t^{k}_{j}),\bar{x}^{k}(t^{k}_{j}))\}$ to $\bar{y}(t)$ on $[0,T]$ as $k\to\infty$, which is established Theorem~\ref{disc:conver}.\vspace*{-0.03in}

For each $t\in[t^{k}_{j},t^{k}_{j+1})$ with $j=0,\ldots,k-1$, consider the quintuples
\begin{equation}\label{theta}
\theta^{k}(t)=\big(\theta^{kx}(t),\theta^{ky}(t),\theta^{ku}(t),\theta^{ka}(t),\theta^{kb}(t)\big):=\Big(\dfrac{\theta^{kx}_{j}}{h_{k}},\dfrac{\theta^{ky}_{j}}{h_{k}},\dfrac{\theta^{ku}_{j}}{h_{k}},
\dfrac{\theta^{ka}_{j}}{h_{k}},\dfrac{\theta^{kb}_{j}}{h_{k}}\Big),
\end{equation}
where the components in \eqref{theta} are defined in \eqref{e:6.1a} and \eqref{e:6.1}. It easy follows from these constructions and the convergence result of Theorem~\ref{disc:conver} that the sequence $\{\theta^k(\cdot)\}$ converges to 0 strongly in $L^{2}([0,T],\mathbb{R}^{3n+m+d})$, and hence we have that $\theta^{k}(t)\to 0$ as $k\to\infty$ for a.e.\ $t\in[0,T)$ along a subsequence (without relabeling).\vspace*{-0.05in}

Having $p^k_j$ with $j=0,\ldots,k$ from Theorem~\ref{Thm:O.c.d}, construct $q^{k}(\cdot)=(q^{kx}(\cdot),q^{ky}(\cdot),q^{ku}(\cdot),q^{ka}(\cdot),q^{kb}(\cdot))$ by setting  $q^{k}(t^{k}_{j}):=p^{k}_{j}$ and then extending the quintuples piecewise linearly to $[0,T]$ for all $j=0,\ldots,k$. Define further $\sigma^{k}(t)$ and $\varsigma^{k}(t)$ on $[0,T]$ by
\begin{equation}\label{sigma}
\sigma^{k}(t):=\sigma^{k}_{j}\;\mbox{ for }\;t\in[t^{k}_{j},t^{k}_{j+1}),\,j=0,\ldots,k-1,\;\mbox{ with }\;\sigma^{k}(t^{k}_{k}):=0,
\end{equation}
\begin{equation}\label{vsigma}
\varsigma^{k}(t):=\dfrac{p^{kd}_{j+1}}{h_{k}}\;\mbox{ for }\;t\in[t^{k}_{j},t^{k}_{j+1}),\;j=0,\ldots,k-1,\;\mbox{ with }\;\varsigma^{k}(t^{k}_{k}):=0.
\end{equation}
Considering the auxiliary functions
\begin{equation}\label{vtheta}
\vartheta^{k}(t):=\max\big\{t^{k}_{j}\:\big|\;t^{k}_{j}\leq t,\,\,0\leq j\leq k\big\}\;\mbox{ for all }\;t\in[0,T],\,\,\,k\in \mathbb{N},
\end{equation}
it is obvious to see that $\vartheta^{k}(t)$ converge to $t$ uniformly in $[0,T]$ as $k\to\infty$. It follows from \eqref{e:6.33}--\eqref{e:6.37} and the constructions above that for all $t\in (t^{k}_{j},t^{k}_{j+1})$ and $j=0,\ldots,k-1$ we get
\begin{equation*}
\begin{aligned}
\dot{q}^{kx}(t)-\lambda^{k}w^{kx}(t)&=-\xi^{x}(\vartheta^{k}(t))\varsigma^{k}(t)+\nabla_{x}f_{1}\big(\bar{a}^{k}(\vartheta^{k}(t)),\bar{x}^{k}(\vartheta^{k}(t))\big)^{*}\big( \lambda^{k}(v^{kx}(t)+\theta^{kx}(t))-q^{kx}(\vartheta^{k}(t)+h_{k}) \big)\\
&+h_{k}\nabla_{x}f_{2}\big(\bar{b}^{k}(\vartheta^{k}(t)),\bar{x}^{k}(\vartheta^{k}(t))\big)^{*}\big( \lambda^{k}(v^{kx}(t)+\theta^{kx}(t))-q^{kx}(\vartheta^{k}(t)+h_{k}) \big)\\
&-\sum\limits_{i=1}^{s}\eta^{k}_{i}(t)\nabla^{2}_{x}g_{i}\big( \bar{x}^{k}(\vartheta^{k}(t))-\bar{u}^{k}(\vartheta^{k}(t)) \big)\big( \lambda^{k}(v^{kx}(t)+\theta^{kx}(t))-q^{kx}(\vartheta^{k}(t)+h_{k}) \big)\\
&-\sum\limits_{i=1}^{s}\sigma^{k}_{i}(t)\nabla_{x}g_{i}\big( \bar{x}^{k}(\vartheta^{k}(t))-\bar{u}^{k}(\vartheta^{k}(t)) \big),
\end{aligned}
\end{equation*}
\begin{equation*}
\dot{q}^{ky}(t)=\varsigma^{k}(t)+\lambda^{k}(v^{kx}(t)+\theta^{kx}(t))-q^{kx}(\vartheta^{k}(t)+h_{k}),
\end{equation*}
\begin{equation*}
\begin{aligned}
\dot{q}^{ku}(t)-\lambda^{k}w^{ku}(t)&=\sum\limits_{i=1}^{s}\eta^{k}_{i}(t)\nabla^{2}_{u}g_{i}\big( \bar{x}^{k}(\vartheta^{k}(t))-\bar{u}^{k}(\vartheta^{k}(t)) \big)\big( \lambda^{k}(v^{kx}(t)+\theta^{kx}(t))-q^{kx}(\vartheta^{k}(t)+h_{k}) \big)\\
&+\sum\limits_{i=1}^{s}\sigma^{k}_{i}(t)\nabla_{u}g_{i}\big( \bar{x}^{k}(\vartheta^{k}(t))-\bar{u}^{k}(\vartheta^{k}(t)) \big),
\end{aligned}
\end{equation*}
\begin{equation*}
\begin{aligned}
\dot{q}^{ka}(t)-\lambda^{k}w^{ka}(t)=\nabla_{a}f_{1}\big(\bar{a}^{k}(\vartheta^{k}(t)),\bar{x}^{k}(\vartheta^{k}(t))\big)^{*}\big( \lambda^{k}(v^{kx}(t)+\theta^{kx}(t))-q^{kx}(\vartheta^{k}(t)+h_{k}) \big),
\end{aligned}
\end{equation*}
\begin{equation*}
\begin{aligned}
\dot{q}^{kb}(t)-\lambda^{k}w^{kb}(t)&=-\xi^{b}(\vartheta^{k}(t))\varsigma^{k}(t)+h_{k}\nabla_{b}f_{2}\big(\bar{a}^{k}(\vartheta^{k}(t)),\bar{x}^{k}(\vartheta^{k}(t))\big)^{*}\big( \lambda^{k}(v^{kx}(t)+\theta^{kx}(t))-q^{kx}(\vartheta^{k}(t)+h_{k}) \big).
\end{aligned}
\end{equation*}\vspace*{-0.1in}

The next adjoint quintuple is $p^{k}(\cdot)=(p^{kx}(\cdot),p^{ky}(\cdot),p^{ku}(\cdot),p^{ka}(\cdot),p^{kb}(\cdot))$ defined
by
\begin{equation}\label{e:7.16}
\begin{aligned}
p^{k}(t):=q^{k}(t)&+\int\limits_{t}^{T}\Big( -\sum\limits_{i=1}^{s}\eta^{k}_{i}(\tau)\nabla^{2}_{x}g_{i}\big( \bar{x}^{k}(\vartheta^{k}(\tau))-\bar{u}^{k}(\vartheta^{k}(\tau)) \big)\big( \lambda^{k}(v^{kx}(\tau)+\theta^{kx}(\tau))-q^{kx}(\vartheta^{k}(\tau)+h_{k}) \big)\\
&-\sum\limits_{i=1}^{s}\sigma^{k}_{i}(\tau)\nabla_{x}g_{i}\big( \bar{x}^{k}(\vartheta^{k}(\tau))-\bar{u}^{k}(\vartheta^{k}(\tau)) \big) - \xi^{x}(\vartheta^{k}(\tau))\varsigma^{k}(\tau),\;\varsigma^{k}(\tau) ,\\
& \sum\limits_{i=1}^{s}\eta^{k}_{i}(\tau)\nabla^{2}_{u}g_{i}\big( \bar{x}^{k}(\vartheta^{k}(\tau))-\bar{u}^{k}(\vartheta^{k}(\tau)) \big)\big( \lambda^{k}(v^{kx}(\tau)+\theta^{kx}(\tau))-q^{kx}(\vartheta^{k}(\tau)+h_{k}) \big)\\
&+\sum\limits_{i=1}^{s}\sigma^{k}_{i}(\tau)\nabla_{u}g_{i}\big( \bar{x}^{k}(\vartheta^{k}(\tau))-\bar{u}^{k}(\vartheta^{k}(\tau)) \big) ,\;0,\;-\xi^{b}(\vartheta^{k}(\tau))\varsigma^{k}(\tau)\Big)\, d\tau
\end{aligned}
\end{equation}
for all $t\in[0,T]$. This gives us $p^{k}(T)=q^{k}(T)$ with the pointwise derivative relationship
\begin{equation}\label{e:7.17}
\begin{aligned}
\dot{p}^{k}(t)=\dot{q}^{k}(t)&-\Big( -\sum\limits_{i=1}^{s}\eta^{k}_{i}(t)\nabla^{2}_{x}g_{i}\big( \bar{x}^{k}(\vartheta^{k}(t))-\bar{u}^{k}(\vartheta^{k}(t)) \big)\big( \lambda^{k}(v^{kx}(t)+\theta^{kx}(t))-q^{kx}(\vartheta^{k}(t)+h_{k}) \big)\\
&-\sum\limits_{i=1}^{s}\sigma^{k}_{i}(t)\nabla_{x}g_{i}\big( \bar{x}^{k}(\vartheta^{k}(t))-\bar{u}^{k}(\vartheta^{k}(t)) \big) - \xi^{x}(\vartheta^{k}(t))\varsigma^{k}(t),\;\varsigma^{k}(t) ,\\
& \sum\limits_{i=1}^{s}\eta^{k}_{i}(t)\nabla^{2}_{u}g_{i}\big( \bar{x}^{k}(\vartheta^{k}(t))-\bar{u}^{k}(\vartheta^{k}(t)) \big)\big( \lambda^{k}(v^{kx}(t)+\theta^{kx}(t))-q^{kx}(\vartheta^{k}(t)+h_{k}) \big)\\
&+\sum\limits_{i=1}^{s}\sigma^{k}_{i}(t)\nabla_{u}g_{i}\big( \bar{x}^{k}(\vartheta^{k}(t))-\bar{u}^{k}(\vartheta^{k}(t)) \big) ,  0 , -\xi^{b}(\vartheta^{k}(t))\varsigma^{k}(t)\Big)\;\mbox{ a.e. }\;t\in[0,T].
\end{aligned}
\end{equation}
Furthermore, we deduce from the above that the componentwise equalities
\begin{equation}\label{e:7.18}
\begin{aligned}
\dot{p}^{kx}(t)-\lambda^{k}w^{kx}(t)&=\nabla_{x}f_{1}\big(\bar{a}^{k}(\vartheta^{k}(t)),\bar{x}^{k}(\vartheta^{k}(t))\big)^{*}\big( \lambda^{k}(v^{kx}(t)+\theta^{kx}(t))-q^{kx}(\vartheta^{k}(t)+h_{k}) \big)\\
&+h_{k}\nabla_{x}f_{2}\big(\bar{b}^{k}(\vartheta^{k}(t)),\bar{x}^{k}(\vartheta^{k}(t))\big)^{*}\big( \lambda^{k}(v^{kx}(t)+\theta^{kx}(t))-q^{kx}(\vartheta^{k}(t)+h_{k}) \big),
\end{aligned}
\end{equation}
\begin{equation}\label{e:7.19}
\dot{p}^{ky}(t)=\lambda^{k}(v^{kx}(t)+\theta^{kx}(t))-q^{kx}(\vartheta^{k}(t)+h_{k})
\end{equation}
\begin{equation}\label{e:7.20}
\dot{p}^{ku}(t)-\lambda^{k}w^{ku}(t)=0,
\end{equation}
\begin{equation}\label{e:7.21}
\dot{p}^{ka}(t)-\lambda^{k}w^{ka}(t)= \nabla_{a}f_{1}\big(\bar{a}^{k}(\vartheta^{k}(t)),\bar{x}^{k}(\vartheta^{k}(t))\big)^{*}\big( \lambda^{k}(v^{kx}(t)+\theta^{kx}(t))-q^{kx}(\vartheta^{k}(t)+h_{k}) \big),
\end{equation}
\begin{equation}\label{e:7.22}
\dot{p}^{kb}(t)-\lambda^{k}w^{kb}(t)=h_{k}\nabla_{b}f_{2}\big(\bar{a}^{k}(\vartheta^{k}(t)),\bar{x}^{k}(\vartheta^{k}(t))\big)^{*}\big( \lambda^{k}(v^{kx}(t)+\theta^{kx}(t))-q^{kx}(\vartheta^{k}(t)+h_{k}) \big)
\end{equation}
hold for every $t\in (t^{k}_{j},t^{k}_{j+1})$, $j=0,\ldots,k-1$, and $i=1,\ldots,s$.\vspace*{-0.05in}

For each $k\in\N$, define the vector measure $\gg^{k}$ on $[0,T]$ by
\begin{equation}\label{e:7.23}
\begin{aligned}
\int\limits_{B}d\,\gg^{k}&:=\int\limits_{B}\Big(\sum\limits_{i=1}^{s}\eta^{k}_{i}(t)\nabla^{2}g_{i}\big( \bar{x}^{k}(\vartheta^{k}(t))-\bar{u}^{k}(\vartheta^{k}(t)) \big)\big( \lambda^{k}(v^{kx}(t)+\theta^{kx}(t))-q^{kx}(\vartheta^{k}(t)+h_{k})\Big)\,dt \\
&+ \int\limits_{B}\Big( \sum\limits_{i=1}^{s}\sigma^{k}_{i}(t)\nabla g_{i}\big( \bar{x}^{k}(\vartheta^{k}(t))-\bar{u}^{k}(\vartheta^{k}(t))\Big)\,dt
\end{aligned}
\end{equation}
for any Borel subset $B\subset [0,T]$, where the vector function $\sigma^k(t)=(\sigma^k_1(t),\ldots,\sigma^k_s(t))$ is taken from \eqref{sigma}. Normalizing the nontriviality condition \eqref{e:6.31}, we get that
\begin{equation}\label{e:7.24}
\lambda^{k}+\lVert p^{k}(T) \rVert+\lVert q^{ku}(0) \rVert + \sum\limits_{j=0}^{k}\lVert p^{kd}_{j} \rVert =1,\,\,\,k\in\mathbb{N},
\end{equation}
which tells us that all the sequential terms in \eqref{e:7.24} are uniformly bounded. Thus we have without loss of generality that there exists $\lm\ge 0$ with $\lm^k\to\lm$ as $k\to\infty$. It follows from \eqref{vsigma} that
\begin{equation}\label{e:7.25}
\int\limits_{0}^{T}\lVert \varsigma^{k}(t) \rVert\,dt=\sum\limits_{j=0}^{k-1}\int\limits_{t^{k}_{j}}^{t^{k}_{j+1}}\lVert \varsigma^{k}(t) \rVert\,dt=\sum\limits_{j=0}^{k-1}\int\limits_{t^{k}_{j}}^{t^{k}_{j+1}}\dfrac{\lVert p^{kd}_{j+1} \rVert}{h_{k}}\,dt=\sum\limits_{j=0}^{k-1}\lVert p^{kd}_{j+1} \rVert=\sum\limits_{j=1}^{k}\lVert p^{kd}_{j} \rVert\leq 1,\quad k\in\N.
\end{equation}\vspace*{-0.1in}

Arguing further as in the proof of \cite[Theorem~6.1]{b1} gives us the following:\\[1ex]
$\bullet$ The boundedness and the uniform bounded variations of $\{q^{k}(\cdot)\}$. Hence the Helly theorem ensures the existence of a function $q(\cdot)$ with bounded variation on $[0,T]$ such that $q^{k}(t)\to q(t)$ as $k\to\infty$ for all $t\in[0,T]$.\\[1ex]
$\bullet$ The boundedness of the sequence $\{p^{k} (\cdot)\}$ in $W^{ 1,2} ([0, T ], \mathbb{R}^{3n+m+d})$, and hence its weak compactness in this space. It follows therefore from the aforementioned Mazur theorem  and basic real analysis that there exists $p(\cdot) \in W^{ 1,2} ([0, T ], \mathbb{R}^{3n+m+d} )$ such that a sequence of convex combinations of $\dot{p}^{k}(t)$ converges to $\dot{p}(t)$ for a.e.\ $t\in [0,T]$. Then we arrive at \eqref{e:7.6} by passing to the limit in \eqref{e:7.18}--\eqref{e:7.22} as $k\to\infty$.

Our next step is to show that the sequence $\{\gg^k(\cdot)\}$ is bounded in ${\cal C}([0,T],\mathbb{R}^{n})^{*}$. To proceed, take any Borel subset $B\subset [0,T] $ and deduce from \eqref{e:6.35} that
\begin{equation*}
\begin{aligned}
&\Big\Vert\int\limits_{B}\Big(\sum\limits_{i=1}^{s}\eta^{k}_{i}(t)\nabla_{u}^{2}g_{i}\big(\bar{x}^{k}(\vartheta^{k}(t))-\bar{u}^{k}
(\vartheta^{k}(t))\big)\big( \lambda^{k}(v^{kx}(t)+\theta^{kx}(t))-q^{kx}(\vartheta^{k}(t)+h_{k})\Big)\,dt\\
&+\int\limits_{B}\Big(\sum\limits_{i=1}^{s}\sigma^{k}_{i}(t)\nabla_{u}g_{i}\big( \bar{x}^{k}(\vartheta^{k}(t))-\bar{u}^{k}(\vartheta^{k}(t)) \Big)\,dt\Big\Vert\\
&\le\Big\Vert\int\limits_{0}^{T}\Big(\sum\limits_{i=1}^{s}\eta^{k}_{i}(t)\nabla_{u}^{2}g_{i}\big(\bar{x}^{k}(\vartheta^{k}(t))
-\bar{u}^{k}(\vartheta^{k}(t))\big)\big(\lambda^{k}(v^{kx}(t)+\theta^{kx}(t))-q^{kx}(\vartheta^{k}(t)+h_{k})\Big)\,dt\\
&+\int\limits_{0}^{T}\Big( \sum\limits_{i=1}^{s}\sigma^{k}_{i}(t)\nabla_{u} g_{i}\big( \bar{x}^{k}(\vartheta^{k}(t))-\bar{u}^{k}(\vartheta^{k}(t))\Big)\,dt\Big\Vert\\
&=\Big\Vert\sum\limits_{j=0}^{k-1}\int\limits_{t_{j}}^{t_{j+1}}\Big(\sum\limits_{i=1}^{s}\eta^{k}_{i}(t)\nabla_{u}^{2}g_{i}\big( \bar{x}^{k}(\vartheta^{k}(t))-\bar{u}^{k}(\vartheta^{k}(t))\big)\big(\lambda^{k}(v^{kx}(t)+\theta^{kx}(t))-q^{kx}(\vartheta^{k}(t)
+h_{k})\Big)\,dt \\
&+\sum\limits_{j=0}^{k-1}\int\limits_{t_{j}}^{t_{j+1}}\Big( \sum\limits_{i=1}^{s}\sigma^{k}_{i}(t)\nabla_{u} g_{i}\big( \bar{x}^{k}(\vartheta^{k}(t))-\bar{u}^{k}(\vartheta^{k}(t))\Big)\,dt\Big\Vert=\Big\Vert\sum\limits_{j=0}^{k-1}\int\limits_{t_{j}}^{t_{j+1}}\Big
(\dfrac{p^{u}_{j+1}-p^{u}_{j}}{h_{k}}-\lambda^{k}w^{u}_{j}\Big)\,dt\Big\Vert\\
&=\Big\Vert\sum\limits_{j=0}^{k-1}(p^{u}_{j+1}-p^{u}_{j}-h_{k}
\lambda^{k}w^{u}_{j})\Big\Vert\le\sum\limits_{j=0}^{k-1}\Vert p^{u}_{j+1}-p^{u}_{j}\Vert+\lambda^{k}\sum\limits_{j=0}^{k-1}\Vert h_{k}w^{u}_{j} \Vert\\
&=\sum\limits_{j=0}^{k-1} \Vert q^{u}(t_{j+1})-q^{u}(t_{j}) \Vert + \lambda^{k}\sum\limits_{j=0}^{k-1}\Vert h_{k}w^{u}_{j} \Vert.
\end{aligned}
\end{equation*}
Then from \eqref{e:7.24}, the imposed structure of the running cost \eqref{run-cost} with the Lipschitzian functions therein, and the uniform bounded variations of $\{q^{u}(\cdot)\} $ we deduce that the sequence of
\begin{equation*}
\begin{aligned}
&\int\limits_{B}\Big(\sum\limits_{i=1}^{s}\eta^{k}_{i}(t)\nabla_{u}^{2}g_{i}\big( \bar{x}^{k}(\vartheta^{k}(t))-\bar{u}^{k}(\vartheta^{k}(t)) \big)\big( \lambda^{k}(v^{kx}(t)+\theta^{kx}(t))-q^{kx}(\vartheta^{k}(t)+h_{k})\Big)\,dt\\
&+ \int\limits_{B}\Big( \sum\limits_{i=1}^{s}\sigma^{k}_{i}(t)\nabla_{u} g_{i}\big( \bar{x}^{k}(\vartheta^{k}(t))-\bar{u}^{k}(\vartheta^{k}(t))\Big)\,dt
\end{aligned}
\end{equation*}
is uniformly bounded on $[0,T]$. In the same way it follows from \eqref{e:6.33}, \eqref{e:7.24}, and the uniform bounded variations of $ \{q^{x}(\cdot)\} $ that the sequence of
\begin{equation*}
\begin{aligned}
&\int\limits_{B}\Big(\sum\limits_{i=1}^{s}\eta^{k}_{i}(t)\nabla_{x}^{2}g_{i}\big( \bar{x}^{k}(\vartheta^{k}(t))-\bar{u}^{k}(\vartheta^{k}(t)) \big)\big( \lambda^{k}(v^{kx}(t)+\theta^{kx}(t))-q^{kx}(\vartheta^{k}(t)+h_{k})\Big)\,dt\\
&+ \int\limits_{B}\Big( \sum\limits_{i=1}^{s}\sigma^{k}_{i}(t)\nabla_{x} g_{i}\big( \bar{x}^{k}(\vartheta^{k}(t))-\bar{u}^{k}(\vartheta^{k}(t))\Big)\,dt
\end{aligned}
\end{equation*}
is also uniformly bounded on $[0,T]$. This verifies the boundedness in ${\cal C}([0,T],\mathbb{R}^{n})^{*}$ of the sequence
$\{\gg^{k}\}$.\vspace*{-0.05in}

Thus we get from the weak$^*$ sequential compactness of the unit ball in ${\cal C}([0,T],\mathbb{R}^{n})^{*}$  that there exists
a measure $\gg\in{\cal C}([0,T],\mathbb{R}^{n})^{*}$ such that $\{\gg^k\}$ weak$^{*}$ converges to $\gg$ along a subsequence (without relabeling).
This allows us to derive from \eqref{e:7.23} and the construction of the measures $\gg^k$ in \eqref{e:7.23} that for all $t\in [0,T]$ the following convergence holds:
\begin{equation*}
\begin{aligned}
&\int\limits_{[t,T]}\Big(\sum\limits_{i=1}^{s}\eta^{k}_{i}(t)\nabla^{2}_{x}g_{i}\big( \bar{x}^{k}(\vartheta^{k}(t))-\bar{u}^{k}(\vartheta^{k}(t)) \big)\big( \lambda^{k}(v^{kx}(t)+\theta^{kx}(t))-q^{kx}(\vartheta^{k}(t)+h_{k})\Big)\,dt \\
&+ \int\limits_{[t,T]}\Big( \sum\limits_{i=1}^{s}\sigma^{k}_{i}(t)\nabla_{x}g_{i}\big( \bar{x}^{k}(\vartheta^{k}(t))-\bar{u}^{k}(\vartheta^{k}(t))\Big)\,dt\to\int\limits_{[t,T]}d\gg(s)\;\mbox{ as }\;k\to\infty.
\end{aligned}
\end{equation*}
Employing now \eqref{e:6.26}, the definition of $\varsigma^{k}(\cdot)$ in \eqref{vsigma}, as well as the convergence $q^{ky}(t)\to q^{y}(t)$ and $\theta^{ky}(t)\to 0$ for a.e. $t\in[0,T]$ as $k\to\infty$ implies that
\begin{equation*}
\varsigma^{k}(t)\to q^{y}(t)\;\mbox{ a.e. }\;t\in[0,T]\;\mbox{ as }\;k\to\infty.
\end{equation*}
Observe by \eqref{e:5.1} and \eqref{e:6.1} that $\dfrac{\theta^{ky}_{j}}{h_{k}}=\dfrac{1}{h_{k}}\displaystyle\int\limits_{t^{k}_{j}}^{t^{k}_{j+1}}(f_{2}(\bar{b}^{k}_{j},\bar{x}^{k}_{j})
-f_{2}(\bar{b}(t),\bar{x}(t)))\,dt$, which leads us to the estimates
\begin{equation*}
\begin{aligned}
\lVert\varsigma^{k}(t) \rVert&\le\lVert q^{ky}(\vartheta^{k}(t)+h_{k}) \rVert +\dfrac{\lambda^{k}}{h_{k}}\int\limits_{t^{k}_{j}}^{t^{k}_{j+1}}\lVert f_{2}(\bar{b}^{k}_{j},\bar{x}^{k}_{j})-f_{2}(\bar{b}(t),\bar{x}(t))\rVert\,dt\\
&\leq  \lVert q^{ky}(\vartheta^{k}(t)+h_{k}) \rVert +L_{1}\dfrac{\lambda^{k}}{h_{k}}\int\limits_{t^{k}_{j}}^{t^{k}_{j+1}}\lVert \bar{x}^{k}_{j}-\bar{x}(t)  \rVert \,dt+L_{2}\dfrac{\lambda^{k}}{h_{k}}\int\limits_{t^{k}_{j}}^{t^{k}_{j+1}}\lVert \bar{b}^{k}_{j}-\bar{b}(t)  \rVert \,dt.
\end{aligned}
\end{equation*}
By the optimality of $\bar{z}^{k}$ in $(P_k)$ and the constraints in \eqref{e:5.2} we get
\begin{equation*}
\lVert \varsigma^{k}(t) \rVert\leq \lVert q^{ky}(\vartheta^{k}(t)+h_{k}) \rVert + \lambda^{k} (L_{1}+L_{2})\dfrac{\varepsilon}{2}.
\end{equation*}
The boundedness of $\{q^{k}(\cdot)\}$ and the normalization condition \eqref{e:7.24} yield the boundedness of $\{\vartheta^{k}(\cdot)\}$. Then it follows from the Lebesgue dominate convergence theorem that
\begin{equation}\label{conv}
\int\limits_{[t,T]}\varsigma^{k}(s)\,ds\to\int\limits_{[t,T]}q^{y}(s)\,ds\;\mbox{ as }\;k\to\infty\;\mbox{ for all }\;t\in[0,T],
\end{equation}
and that the function $q^{y}(\cdot)$ belongs to $L^{1}([0,T],\mathbb{R}^{n})$. Furthermore, for all $t\in [0,T]$ we have
\begin{equation*}
\begin{aligned}
&\;\quad\Big\lVert \int\limits_{[t,T]}  \xi^{x}(\vartheta^{k}(s))\varsigma^{k}(s)\,ds-\int\limits_{[t,T]}  \nabla_{x} f_{2}(\bar{b}(s),\bar{x}(s))q^{y}(s)\,ds  \Big\rVert\\
&=\Big\lVert \int\limits_{[t,T]} \nabla_{x} f_{2}(\bar{b}^{k}(\vartheta^{k}(s)),\bar{x}^{k}(\vartheta^{k}(s)))\varsigma^{k}(s)  \,ds-\int\limits_{[t,T]}  \nabla_{x} f_{2}(\bar{b}(s),\bar{x}(s))q^{y}(s)\,ds  \Big\rVert\\
&\leq \Big\lVert \int\limits_{[t,T]} \nabla_{x} f_{2}(\bar{b}^{k}(\vartheta^{k}(s)),\bar{x}^{k}(\vartheta^{k}(s)))\vartheta^{k}(s)  \,ds-\int\limits_{[t,T]} \nabla_{x} f_{2}(\bar{b}(s),\bar{x}(s))\varsigma^{k}(s)\,ds  \Big\rVert\\
&+\Big\lVert  \int\limits_{[t,T]}  \nabla_{x} f_{2}(\bar{b}(s),\bar{x}(s))\varsigma^{k}(s)\,ds - \int\limits_{[t,T]}  \nabla_{x} f_{2}(\bar{b}(s),\bar{x}(s))q^{y}(s)\,ds\Big\rVert.
\end{aligned}
\end{equation*}
Note first that $\Big\lVert \displaystyle\int\limits_{[t,T]}\nabla_{x} f_{2}(\bar{b}^{k}(\vartheta^{k}(s)),\bar{x}^{k}(\vartheta^{k}(s)))\varsigma^{k}(s)\,ds-\displaystyle\int\limits_{[t,T]}\nabla_{x} f_{2}(\bar{b}(s),\bar{x}(s))\varsigma^{k}(s)\,ds\Big\rVert\to 0$ as $k\to\infty$ due to the uniform convergence $\bar{b}^{k}(\cdot)\to\bar{b}(\cdot)$ and $\bar{x}^{k}(\cdot)\to\bar{x}(\cdot)$ on $[0,T]$ as $k\to\infty$, continuous differentiability of $f_{2}(\cdot,\cdot)$, and the uniform boundedness of $\displaystyle\int\limits_{0}^{T}\lVert \varsigma^{k}(t) \rVert \,dt$ over $k\in\N$ that follows from \eqref{e:7.25}. For the second term in the above estimate we get $\Big\lVert \displaystyle\int\limits_{[t,T]}  \nabla_{x} f_{2}(\bar{b}(s),\bar{x}(s))\varsigma^{k}(s)\,ds - \displaystyle\int\limits_{[t,T]} \nabla_{x} f_{2}(\bar{b}(s),\bar{x}(s))q^{y}(s)\,ds  \Big\rVert \to 0$  as $k\to \infty$  by the convergence of  $\varsigma^{k}(\cdot)\to q^{y}(\cdot)$ as $k\to\infty$ in $L^{1}([0,T],\mathbb{R}^{n})$. Combining this yields
\begin{equation*}
\int\limits_{[t,T]}  \xi^{x}(\vartheta^{k}(s))\varsigma^{k}(s)\,ds \to \int\limits_{[t,T]}  \nabla_{x} f_{2}(\bar{b}(s),\bar{x}(s))q^{y}(s)\,ds \;\mbox{ as }\;k \to \infty .
\end{equation*}
In the same way we also verify that
\begin{equation*}
\int\limits_{[t,T]} \xi^{b}(\vartheta^{k}(s))\varsigma^{k}(s)\,ds \to \int\limits_{[t,T]} \nabla_{b} f_{2}(\bar{b}(s),\bar{x}(s))q^{y}(s)\,ds\;\mbox{ as }\;k \to \infty
\end{equation*}
and hence arrive at the adjoint relationships \eqref{e:7.8} by passing to the limit in \eqref{e:7.16} as $k\to\infty$, where the justification of the ``except a countable subset" can be done similarly to \cite[p.\ 325]{v}. The passage to the limit in \eqref{e:7.11} as $k\to\infty$ brings us to the integral equation
\begin{equation*}
q^{y}(t)=p^{y}(t)-\int\limits_{[t,T]}q^{y}(s)\,ds\;\mbox{ a.e. }\;t\in [0,T],
\end{equation*}
which tells us that $q^{y}(\cdot)$ is an absolutely continuous function satisfying
\begin{equation}\label{e:7.26}
\dot{q}^{y}(t)=\dot{p}^{y}(t)+q^{y}(t)\;\mbox{ a.e. }\;t\in [0,T] .
\end{equation}
It follows from \eqref{e:7.6} and \eqref{e:7.8} that
\begin{equation}\label{e:7.28}
\dot{p}^{y}(t)=\lambda v^{x}(t)-q^{x}(t)\;\mbox{ a.e. }\;t\in [0,T],\;\mbox{ and}
\end{equation}
\begin{equation}\label{e:7.27}
q^{x}(t)=p^{x}(t)+\int\limits_{[t,T]}d\,\gg(s)+\int\limits_{[t,T]} \nabla_{x} f_{2}(\bar{b}(s),\bar{x}(s))q^{y}(s)\,ds\;\mbox{ a.e. }\;t\in[0,T].
\end{equation}
Substituting \eqref{e:7.27} into \eqref{e:7.28} and then \eqref{e:7.28} into \eqref{e:7.26}, we obtain the extended Volterra condition \eqref{e:7.9}. Furthermore, it follows from \eqref{e:6.26}--\eqref{e:6.29} and the definition of $\vartheta(t)$ in \eqref{vtheta}
that for each $k\in\N$ we have
\begin{equation}\label{e:7.29}
q^{ku}(\vartheta^{k}(t)+h_{k})=\lambda^{k}(\theta^{ku}(t)+v^{ku}(t)),
\end{equation}
\begin{equation}\label{e:7.30}
q^{ka}(\vartheta^{k}(t)+h_{k})=\lambda^{k}(\theta^{ka}(t)+v^{ka}(t)),\;\mbox{ and }\;q^{kb}(\vartheta^{k}(t)+h_{k})=\lambda^{k}(\theta^{kb}(t)+v^{kb}(t)).
\end{equation}
Involving now \eqref{e:7.1}, the assumptions on $l_{2}$, $l_{3}$, $l_{4}$, and the Lebesgue dominated convergence theorem  gives us both conditions in \eqref{e:7.7} by passing to the limit in \eqref{e:7.29} and \eqref{e:7.30}. The proof of the complementarity slackness conditions \eqref{e:7.4}, \eqref{e:7.5} and the endpoint conditions in \eqref{e:7.10}--\eqref{e:7.13} is similar to the one in \cite[Theorem~6.1]{b1}, and we skip it for brevity. Taking into account the convergence in \eqref{conv} and the relationship
\begin{equation*}
\int\limits_{0}^{T}\lVert \varsigma^{k}(t) \rVert\,dt=\sum\limits_{j=0}^{k-1}\lVert p^{kd}_{j+1} \rVert,
\end{equation*}
we arrive at the nontriviality condition \eqref{e:7.14} by passing to the limit in \eqref{e:7.24} as $k\to\infty$. The verification of the enhanced nontriviality conditions in \eqref{ench:1} and \eqref{ench:2} under the imposed additional assumptions can be easily proved while arguing by contraposition. This therefore completes the proof of the theorem. $\h$\vspace*{-0.2in}

%%%%%%%%%%%%%%%%%%%%%%%%%%%%%%%%%%%%%%%%%%%%%%%%%%%%%%%%%%%%%%%%%%%%%%%%%%%%%%%%%%%%%%%%%%%%%%%%%%%%%%%%%%%%%%%%%%%%%%%%%%%%%%%%%%%%%%%%%%%%%%%%%%%%%%%%%%%%%%%%%%%%%%%%%%%%%%%%%%%%%%%%%%%%%%%%%%%%%%%%%%%%%%%%%%%%%%%%%%%%%%%%%%%%%%%%%%%%%%%%%%%%%%%%%%%%%%%%%%%%%%%%%%%%%%%%%%%
\section{Applications to Control of Non-Regular Electrical Circuits}\label{sec:exa}
\setcounter{equation}{0}\vspace*{-0.1in}

This section is entirely devoted to applications of the necessary optimality conditions for integro-differential sweeping control problems obtained in Theorem~\ref{Thm:N.C.I} to controlled models that appear in non-regular electrical circuits with ideal diodes. However, the dynamics in such models has been described via (uncontrolled) integro-differential sweeping processes (see, e.g., \cite{abb,bou}), we are not familiar with formulations of any optimization and/or control problems for such systems. This is done in the first two examples below in different frameworks. The third example of its own interest provides a complete solution of a two-dimensional optimal control problem of the type modeled above by using the obtained necessary optimality conditions.\vspace*{-0.1in}
\begin{figure}[h!]
\begin{center}
\includegraphics[width=6in]{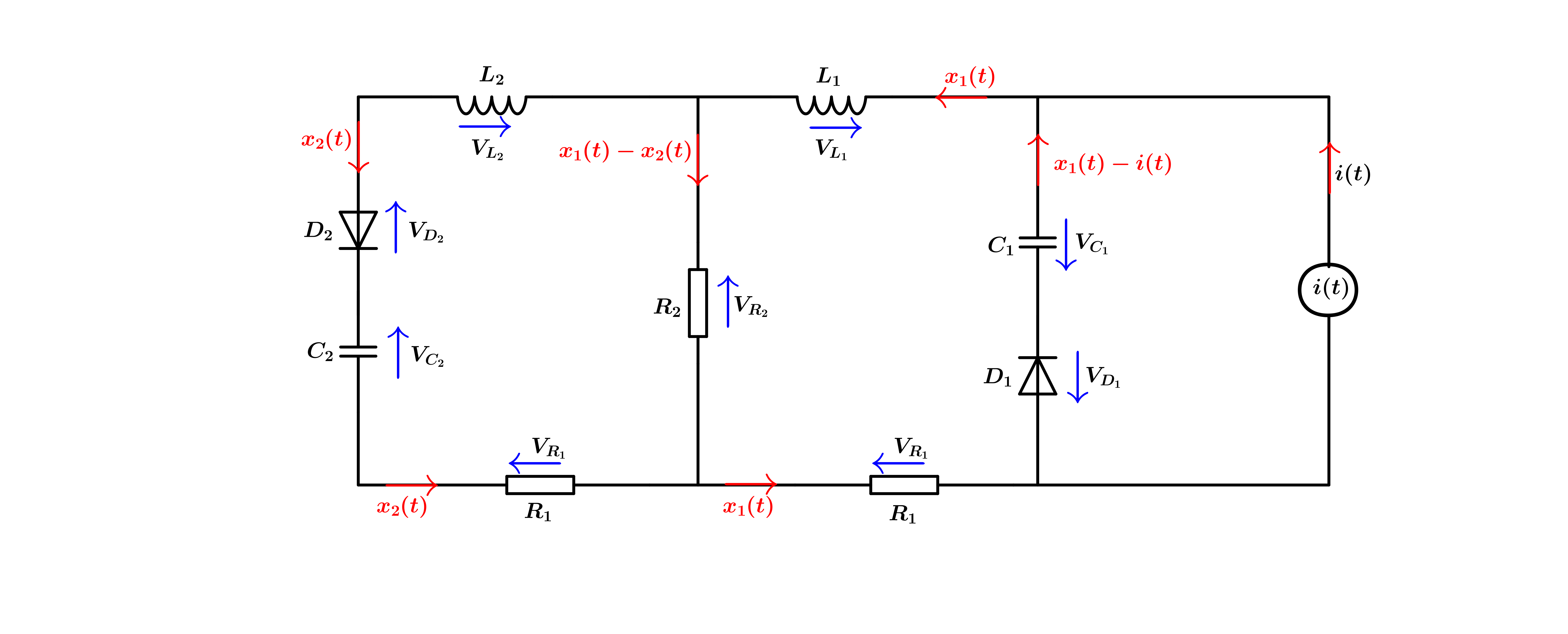}
\caption{Electrical circuit with resistors,Inductances, capacitors and ideal diodes (RLCD).}\label{fig2}
\end{center}
\end{figure}
\begin{figure}[h!]
\begin{center}
\includegraphics[width=6in]{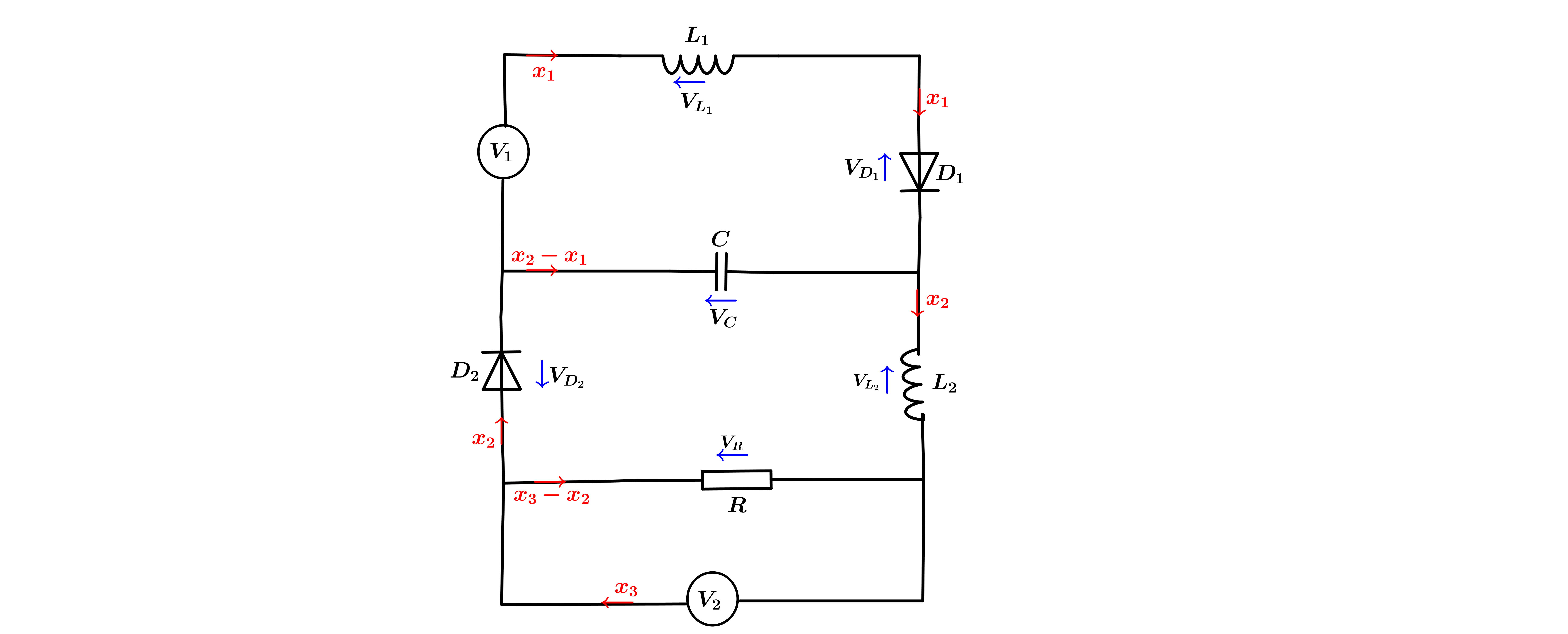}
\caption{(RLCD) circuit with controlled voltage source.}\label{fig3}
\end{center}
\end{figure}
\begin{example}[\bf optimal control of non-regular electrical circuits with controlled current source]\label{ex1} Consider the electrical system depicted in Figure~\ref{fig2} that is composed of two resistors $R_1\geq 0$ and $R_2\geq 0$ with voltage/current laws $V_{R_j}=R_j x_j$ ($j=1,2$), two inductors $L_1\geq 0$ and $L_2\geq 0$ with voltage/current laws $V_{L_j}=L_j \dot x_j$ ($j=1,2$), two capacitors $C_1>0$ and $C_2>0$ with voltage/current laws $V_{C_j}=\frac{1}{C_j} \int x_j(t)\,dt$ $(j=1,2$), two ideal diodes with characteristics $0\leq -V_{D_j}\perp i_j\geq 0$ and an absolutely continuous current source $i:[0, T]\rightarrow R$. Using Kirchhoff's laws, we have
$$\left\{
\begin{array}{l}
V_{R_1}+V_{R_2}+V_{L_1}+V_{C_1}=-V_{D_1}\in -N_{\R_+}(x_1-i),\\
V_{R_1}-V_{R_2}+V_{L_2}+V_{C_2}=-V_{D_2}\in-N_{\R_+}(x_2).
\end{array}
\right.
$$
Therefore, the dynamics of this circuit is given by
\begin{equation}\label{eq2.3}
\begin{array}{ll}
\overbrace{
\begin{pmatrix}
-\dot x_1(t)\\
-\dot x_2(t)
\end{pmatrix}}^{-\dot x(t)}&\in N_{C(t)}(x(t))+\overbrace{\begin{pmatrix}
\frac{R_1 +R_2}{L_1}& -\frac{R_2}{L_1}\\
-\frac{R_2}{L_2}&\frac{R_1 +R_2}{L_2}
\end{pmatrix}}^{A_1}\overbrace{
\begin{pmatrix}
x_1 (t)\\
x_2 (t)
\end{pmatrix}}^{x(t)}\\
&+\disp\int\limits_{0}^{t}\bigg[\overbrace{\begin{pmatrix}
\frac{1}{L_1 C_1}& 0\\
0 &\frac{1}{L_2 C_2}
\end{pmatrix}}^{A_2}\overbrace{
\begin{pmatrix}
x_1(s)\\
x_2(s)
\end{pmatrix}}^{x(s))}+\overbrace{
\begin{pmatrix}
\frac{1}{L_1 C_1}i(s)\\
0
\end{pmatrix}}^{b(s)}\bigg] \,ds.
\end{array}
\end{equation}
Put $u(t):=(i(t),0)^{*}$, $b(t)=(\frac{1}{L_1 C_1}i(t),0)^{*}$ indicating vector columns and denote $C(t)=C(u(t)):=u(t)+[0,\infty)\times [0,\infty)$, $f_{1}(x):=A_1 x$, $f_{2}(b,x):=A_2 x+b$, and $x(0):=(i(0),0)^{*}$. Thus \eqref{eq2.3} can be rewritten in the
form of the controlled sweeping dynamics \eqref{e:1} as
\begin{equation*}
-\dot{x}(t)\in N_{C(u(t))}(x(t))+f_{1}(x(t))+\displaystyle\int\limits_{0}^{t}f_{2}(b(s), x(s))\,ds\;\mbox{ a.e. }\;t\in[0, T],\;x(0)=x_0\in C(0)
\end{equation*}
with control functions acting in the moving set and the integral part of the dynamics. The cost functional in this problem is naturally formulated by: minimize
\begin{equation*}
J[x,u,b]:=\frac{\lambda_{T}}{2}(x_1(T)-i(0))^{2}+\frac{\lambda_{P}}{2}\displaystyle\int\limits_{0}^{T}(x_1(t)-i(0))^{2}\,dt
+\frac{\lambda_{I}}{2}\displaystyle\int\limits_{0}^{T}b_{1}^{2}(t)\,dt,
\end{equation*}
where $\lambda_{T}$, $\lambda_{Q}$, and $\lambda_{I}$ are nonnegative constants not equal to zero simultaneously. The necessary optimality conditions obtained in Theorem~\ref{Thm:N.C.I} can be readily applied to this class of optimal control problems.
\end{example}\vspace*{-0.1in}

The next model describes a class of dynamic processes for non-regular electric circuits with an input controlled voltage source, which is dual of the current source considered in Example~\ref{ex1}. In this model we have control actions entering the dynamic perturbation term of the sweeping process.\vspace*{-0.1in}
\begin{example}[\bf optimal control of non-regular electrical circuit with controlled voltage source]\label{ex2} $\,$
Consider the electrical system shown in Figure~\ref{fig3} that is composed of a resistor $R \geq 0$ with voltage/current law $V_{R}=R (x_3-x_2)$, two inductors $L_1\geq 0$, $L_2\geq 0$ with voltage/current laws $V_{L_j}=L_j \dot x_j$ ($j=1,2$), a capacitor $C>0$ with voltage/current law $V_{C}=\dfrac{1}{C}\displaystyle\int (x_2(t)-x_1(t))\,dt$, two ideal diodes with characteristics $0\leq -V_{D_j}\perp x_j\geq 0$, and two voltage sources $v_j(t)$ ($j=1,2$). Employing Kirchhoff's laws tells us that
$$\left\{
\begin{array}{l}
V_{L_1}-V_{C}-v_{1}(t)=-V_{D_1}\in -N_{\mathbb{R}_+}(x_1),\\
V_{L_2}+V_{C}-V_{R}=-V_{D_2}\in -N_{\mathbb{R}_+}(x_2),\\
V_{R}-v_{2}(t)=0.
\end{array}
\right.
$$
In the system above, we substitute the third equation into the second one and get
$$
-\begin{pmatrix}
\dot x_1(t)\\
\dot x_2(t)
\end{pmatrix}\in N_{\R^{2}_{+}}\begin{pmatrix}
 x_1(t)\\
x_2(t)
\end{pmatrix}+\overbrace{\begin{pmatrix}
-\frac{1}{L_1}& 0\\
0&-\frac{1}{L_2}
\end{pmatrix}}^{A_1}
\begin{pmatrix}
v_1 (t)\\
v_2 (t)
\end{pmatrix}+
\int\limits_{0}^{t}  \overbrace{\begin{pmatrix}
\frac{1}{L_1 C}& -\frac{1}{L_1 C}\\
-\frac{1}{L_2 C} &\frac{1}{L_2 C}
\end{pmatrix}}^{A_2}
\begin{pmatrix}
x_1(s)\\
 x_2(s)
\end{pmatrix} \,ds.
$$
The goal is to start from a given state $(x_{1}(0),x_{2}(0))$ and reach the other state
$(x_{1}(T),x_{2}(T))$ with minimizing the input energy. Denoting $C:=\{x=(x_{1},x_{2})\in \mathbb{R}^{2}\;|\;g_{1}(x)=x_{1}\geq 0,\,\,\,g_{2}(x)=x_{2}\geq 0\}$, $u(t)=(0,0)$, $x_0 \in C $, $a(t)=(a_{1} (t),a_{2} (t))$ with $(a_1(t),a_2(t)):=(v_{1} (t),v_{2} (t)$, $f_{1}(a,x):=A_{1} a $, and $f_{2}(b,x):=A_{2} x $, the formulated problem can be written in the form of our basic problem $(P)$ as follows:
\begin{equation}\label{eq2.6bis}
\mbox{minimize }\;J[x,a]:= \dfrac{x_{2}^{2}(T)}{2}+\frac{1}{2}\int\limits_{0}^{T}[ a_{1}^{2}(t)+ a_{2}^{2}(t)]\,dt
\end{equation}
subject to the controlled integro-differential sweeping process
\begin{equation}\label{eq2.7}
-\dot{x}(t)\in N_{C}(x(t))+f_{1}(a(t),x(t))+\displaystyle\int\limits_{0}^{t}f_{2}(x(s))\,ds\;\mbox{  a.e. }\;[0, T],\;
x(0)=x_0\in C.
\end{equation}
It follows from the existence result of Theorem~\ref{exist:opti} that the optimal control problem formulated in \eqref{eq2.6bis} and \eqref{eq2.7} admits an optimal solution. Furthermore, the structure of the cost functional in \eqref{eq2.6bis} together with Proposition~\ref{exist:solu} yields the uniqueness of the optimal pair $(\ox(\cdot),\oa(\cdot))$ in problem \eqref{eq2.6bis}, \eqref{eq2.7}. This allows us to find the optimal solution to the formulated optimal control problems by using the necessary optimality conditions Theorem~\ref{Thm:N.C.I}, provided that a feasible solution determined by these conditions is unique.
\end{example}\vspace*{-0.1in}

Finally, we present a numerical example that demonstrates how to use the necessary optimality conditions of Theorem~\ref{Thm:N.C.I} to solve a particular case of the optimal control problem formulated in Example~\ref{ex2} with the given data of controlled integro-differential sweeping inclusion. The corresponding example for the first model described in this section can be constructed similarly, while we skip it for brevity.\vspace*{-0.1in}

\begin{example}[\bf calculating optimal solutions of two-dimensional sweeping control model]\label{ex3} $\,$ Consider the optimal control problem defined in \eqref{eq2.6bis} and \eqref{eq2.7} with the following data:
\begin{equation}\label{pv}
\begin{array}{ll}
T=1,\;x_0=(1,1),\;g(x)=(g_{1}(x),g_{2}(x)):=(x_{1},x_{2}),\;l_{0}(x,a):=\dfrac{a_{1}^{2}+a_{2}^{2}}{2}.\;\varphi(x)=\dfrac{x^{2}_{2}}{2},\\
A_1: =\begin{pmatrix}
-1& 0\\
0 & -1
\end{pmatrix},\;\mbox{ and }\;A_2 =\begin{pmatrix}
1& -1\\
-1 & 1
\end{pmatrix}.
\end{array}
\end{equation}
We seek for solutions to problem \eqref{pv} such that
\begin{equation}\label{hyp}
\bar{x}(t)>0\;\mbox{ for all }\;t\in [0,1)\;\mbox{ and }\;\bar{x}(1)\in{\rm bd}\,C.
\end{equation}
Depending on whether the diodes are blocking (off) or conducting (on), condition \eqref{hyp} provides the following three possibilities that are listed below as {\em modes}:\\[1ex]
$\bullet$ {\sc Mode~1}: {\em For all $t\in [0,1)$ both diodes are on}, i.e., $x_1 > 0$ $(V_{D_1}=0)$ and $x_2 >0$ $(V_{D_2}=0)$. Furthermore, {\em at the ending time $T=1$ the first diode is off while the second one is on}, i.e., $x_1(1) =0$ and $x_2(1) >0$ $(V_{D_2}=0) $.\\[1ex]
$\bullet$ {\sc Mode 2}: {\em For all $t\in [0,1)$ both diodes are on}, i.e., $x_1 > 0$ $(V_{D_1}=0)$ and $x_2 >0$ $(V_{D_2}=0)$, while {\em at the ending time $T=1$ the first diode is on and the second one is off}, i.e., $x_1(1) > 0$ $(V_{D_1}=0)$ and $x_2(1)=0.$\\[1ex]
$\bullet$ {\sc Mode 3}: {\em For all $t\in [0,1)$ both diodes are on}, i.e., $x_1 > 0$ $(V_{D_1}=0)$ and $x_2 >0$ $(V_{D_2}=0)$, while {\em at the ending time $T=1$ both diodes are off}, i.e., $x_1(1)=0$ and $x_2(1)=0.$\vspace*{-0.05in}

Applying the necessary optimality conditions of Theorem~\ref{Thm:N.C.I} gives us a number $\lambda\geq 0$, functions $\eta_{i}(\cdot)\in L^{2}([0,1],\mathbb{R}_{+})$ as $i=1,2$ well-defined at $t =1$, triples $p(\cdot)=(p^{x}(\cdot),p^{y}(\cdot),p^{a}(\cdot))\in W^{1,2}([0,T];\R^6)$ and $q(\cdot)=(q^{x}(\cdot),q^{y}(\cdot),q^{a}(\cdot))$ with values in $\R^6$ and of bounded variations on $[0,T]$, as well as a pair $w(\cdot)=(w^x(\cdot),w^a(\cdot))\in L^2([0,T],\R^4)$ such that
\begin{enumerate}
\item $ w(t)=(0,0,\bar{a}(t)),\,\;v(t)=(0,0,0)\quad \text{a.e.}\quad t\in[0,1].$
\item $\begin{cases}\dot{\bar{x}}_{1}(t)-\bar{a}_{1}(t)+\disp\int\limits_{0}^{t}(\bar{x}_{1}(s)-\bar{x}_{2}(s))\,ds=\eta_{1}(t)\;
\mbox{ a.e. }\;t\in[0,1], &\\
\dot{\bar{x}}_{2}(t)-\bar{a}_{2}(t)-\disp\int\limits_{0}^{t}(\bar{x}_{1}(s)-\bar{x}_{2}(s))\,ds=\eta_{2}(t)\;\mbox{ a.e. }\;t\in[0,1] .
 \end{cases}$
\item $\bar{x}_{i}(t)>0 \Longrightarrow \eta_{i}(t)=0\;\mbox{ a.e. }\;t\in [0,1],\;\;i=1,2 .$
\item $\eta_{i}(t)>0 \Longrightarrow q_{i}^{x}(t)=0\;\mbox{ a.e. }\;t\in [0,1),\;i=1,2$.
\item $(\dot{p}^{x}(t),\dot{p}^{y}(t),\dot{p}^{a}(t))=((0,0),-q^{x}(t),\lambda \bar{a}(t)-A_{1}q^{x}(t))\;\mbox{ a.e. }\;t\in[0,1] $.
\item $ q^{a}(t)=(0,0)\;\mbox{ a.e. }\;t\in[0,1]$.
\item $(q^{x}(t),q^{y}(t),q^{a}(t))=(p^{x}(t),p^{y}(t),p^{a}(t))-\disp\int\limits_{[t,1]}(-d\gamma(s)-A_{2}q^{y}(s)\,ds,q^{y}(s)\,ds,(0,0))\;\mbox{ a.e. }\;t\in[0,1]$.
\item
$\begin{cases}\dot{q}_{1}^{y}(t)=-p_{1}^{x}(t)-\disp\int\limits_{[t,1]}d\gamma_{1}(s)-\disp\int\limits_{[t,1]}(q_{1}^{y}(s)-q_{2}^{y}(s))\,ds
+q_{1}^{y}(t)\;\mbox{ a.e. }\;t\in[0,1],\\
\dot{q}_{2}^{y}(t)=-p_{2}^{x}(t)-\disp\int\limits_{[t,1]}d\gamma_{2}(s)+\disp\int\limits_{[t,1]}(q_{1}^{y}(s)-q_{2}^{y}(s))\,ds+q_{2}^{y}(t)\;
 \mbox{ a.e. }\;t\in[0,1] .
\end{cases}$
\item $-p_{1}^{x}(1)+\eta_{1}(1)=0,\,\,-p_{2}^{x}(1)+\eta_{2}(1)=\lambda\bar{x}_{2}(1)$.
\item $-\eta(1)\in N_{C}(\bar{x}(1)).$\\
\item $p^y(1)=p^a(1)=(0,0)$.
\item $\lm+|\eta_1(1)|+|\eta_2(1)|+\disp\int\limits^1_0\|q^y(t)\|\,dt> 0$.
\end{enumerate}
It follows directly from items~5--7 that
\begin{equation*}
p^{x}(t)\equiv p^{x}(1)\;\mbox{ and }\;\lambda\bar{a}(t)=-q^{x}(t)\;\mbox{ for all }\;t\in[0,1].
\end{equation*}
Proceeding now similarly to \cite[Example~1]{ca3} gives us $\displaystyle\int\limits_{[t,1]}d\,\gg(s)=\gg(\{1\})$. Then
\begin{equation}\label{e:8.8}
q^{x}(t)=p^{x}(1)+\gg(\{1\})+\disp\int\limits_{[t,1]}A_{2}q^{y}(s)\,ds,\quad t\in[0,1].
\end{equation}
Solving the integro-differential system in item~8 with $q^{y}(1)=(0,0)$, we get on $[0,1]$ that
\begin{equation*}
\begin{cases}
\disp\int\limits_{[t,1]}q_{1}^{y}(s)\,ds=\dfrac{p^{x}_{1}(1)+\gg_{1}(\{1\})+p^{x}_{2}(1)+\gg_{2}(\{1\})}{2}e^{t-1}+\dfrac{p^{x}_{1}(1)
+\gg_{1}(\{1\})-p^{x}_{2}(1)-\gg_{2}(\{1\})}{6}e^{-(t-1)}\\
+\dfrac{p^{x}_{1}(1)+\gg_{1}(\{1\})-p^{x}_{2}(1)-\gg_{2}(\{1\})}{12}e^{2t-2}-\dfrac{p^{x}_{1}(1)+\gg_{1}(\{1\})-p^{x}_{2}(1)
-\gg_{2}(\{1\})}{4}-\dfrac{p^{x}_{1}(1)+\gg_{1}(\{1\})+p^{x}_{2}(1)+\gg_{2}(\{1\})}{2}t,\\\\
\disp\int\limits_{[t,1]}q_{2}^{y}(s)\,ds=\dfrac{p^{x}_{1}(1)+\gg_{1}(\{1\})+p^{x}_{2}(1)+\gg_{2}(\{1\})}{2}e^{t-1}-\dfrac{p^{x}_{1}(1)
+\gg_{1}(\{1\})-p^{x}_{2}(1)-\gg_{2}(\{1\})}{6}e^{-(t-1)}\\
-\disp\dfrac{p^{x}_{1}(1)+\gg_{1}(\{1\})-p^{x}_{2}(1)-\gg_{2}(\{1\})}{12}e^{2t-2}+\dfrac{p^{x}_{1}(1)+\gg_{1}(\{1\})-p^{x}_{2}(1)
-\gg_{2}(\{1\})}{4}-\dfrac{p^{x}_{1}(1)+\gg_{1}(\{1\})+p^{x}_{2}(1)+\gg_{2}(\{1\})}{2}t.
\end{cases}
\end{equation*}
Then it follows from the above relationships, equation \eqref{e:8.8}, and the definition of $A_{2}$ that
\begin{equation*}
\begin{cases} q_{1}^{x}(t)=p^{x}_{1}(1)+\gg_{1}(\{1\})+\dfrac{p^{x}_{1}(1)+\gg_{1}(\{1\})-p^{x}_{2}(1)-\gg_{2}(\{1\})}{3}e^{-(t-1)}\vspace*{0.3cm}\\
+\dfrac{p^{x}_{1}(1)+\gg_{1}(\{1\})-p^{x}_{2}(1)-\gg_{2}(\{1\})}{6}e^{2t-2} -\dfrac{p^{x}_{1}(1)+\gg_{1}(\{1\})-p^{x}_{2}(1)-\gg_{2}(\{1\})}{2} &\vspace*{0.5cm}\\
q_{2}^{x}(t)=p^{x}_{2}(1)+\gg_{2}(\{1\})- \dfrac{p^{x}_{1}(1)+\gg_{1}(\{1\})-p^{x}_{2}(1)-\gg_{2}(\{1\})}{3}e^{-(t-1)}\vspace*{0.3cm}\\
-\dfrac{p^{x}_{1}(1)+\gg_{1}(\{1\})-p^{x}_{2}(1)-\gg_{2}(\{1\})}{6}e^{2t-2} +\dfrac{p^{x}_{1}(1)+\gg_{1}(\{1\})-p^{x}_{2}(1)-\gg_{2}(\{1\})}{2}.
\end{cases}
\end{equation*}
This leads us to the representations
\begin{equation*}
\begin{cases} \lambda\bar{a}_{1}(t)=-q_{1}^{x}(t)=-p^{x}_{1}(1)-\gg_{1}(\{1\})-  \dfrac{p^{x}_{1}(1)+\gg_{1}(\{1\})-p^{x}_{2}(1)-\gg_{2}(\{1\})}{3}e^{-(t-1)}\vspace*{0.3cm}\\
-\dfrac{p^{x}_{1}(1)+\gg_{1}(\{1\})-p^{x}_{2}(1)-\gg_{2}(\{1\})}{6}e^{2t-2} +\dfrac{p^{x}_{1}(1)+\gg_{1}(\{1\})-p^{x}_{2}(1)-\gg_{2}(\{1\})}{2}, &\vspace*{0.5cm}\\
\lambda\bar{a}_{2}(t)=-q_{2}^{x}(t)=-p^{x}_{2}(1)-\gg_{2}(\{1\})+ \dfrac{p^{x}_{1}(1)+\gg_{1}(\{1\})-p^{x}_{2}(1)-\gg_{2}(\{1\})}{3}e^{-(t-1)}\vspace*{0.3cm}\\
+\dfrac{p^{x}_{1}(1)+\gg_{1}(\{1\})-p^{x}_{2}(1)-\gg_{2}(\{1\})}{6}e^{2t-2} -\dfrac{p^{x}_{1}(1)+\gg_{1}(\{1\})-p^{x}_{2}(1)-\gg_{2}(\{1\})}{2},
\end{cases}
\end{equation*}
for all $t\in[0,1]$. Suppose now for simplicity that $\lambda>0$ and denote $v_{1}:=\dfrac{p^{x}_{1}(1)+\gg_{1}(\{1\})}{\lambda}$, $v_{2}:=\dfrac{p^{x}_{2}(1)+\gg_{2}(\{1\})}{\lambda}$. Then the nontriviality condition in item~12 holds automatically, and we have
\begin{equation}\label{e:8.9}
\begin{cases} \bar{a}_{1}(t)=-v_{1}-  \dfrac{v_{1}-v_{2}}{3}e^{-(t-1)}
-\dfrac{v_{1}-v_{2}}{6}e^{2t-2} +\dfrac{v_{1}-v_{2}}{2},&\vspace*{0.5cm}\\
\bar{a}_{2}(t)=-v_{2}+ \dfrac{v_{1}-v_{2}}{3}e^{-(t-1)}+\dfrac{v_{1}-v_{2}}{6}e^{2t-2}-\dfrac{v_{1}-v_{2}}{2}.
\end{cases}
\end{equation}
Substituting these expressions into item~2 and then using item~3 yield
\begin{equation*}
\begin{cases}\dot{\bar{x}}_{1}(t)+v_{1}+\dfrac{v_{1}-v_{2}}{3}e^{-(t-1)}+\dfrac{v_{1}-v_{2}}{6}e^{2t-2}-\dfrac{v_{1}-v_{2}}{2}
+\disp\int\limits_{0}^{t}(\bar{x}_{1}(s)-\bar{x}_{2}(s))\,ds=0,&\vspace*{0.5cm}\\
\dot{\bar{x}}_{2}(t)+v_{2}-\dfrac{v_{1}-v_{2}}{3}e^{-(t-1)}-\dfrac{v_{1}-v_{2}}{6}e^{2t-2}+\dfrac{v_{1}-v_{2}}{2}
-\disp\int\limits_{0}^{t}(\bar{x}_{1}(s)-\bar{x}_{2}(s))\,ds=0.
\end{cases}
\end{equation*}
Solving the obtained integro-differential system with $\bar{x}(0)=(1,1)$ gives us
\begin{equation}\label{e:8.11}
\begin{cases}\bar{x}_{1}(t)=1-\dfrac{v_{1}+v_{2}}{2}t -\dfrac{v_{1}-v_{2}}{18}e^{2t-2}+\dfrac{v_{1}-v_{2}}{9}e^{1-t}+\dfrac{1}{3}\cos(\sqrt{2}t)\Big(\dfrac{v_{1}-v_{2}}{6}e^{-2}
-\dfrac{v_{1}-v_{2}}{3}e\Big)\vspace*{0.3cm}\\-\dfrac{1}{\sqrt{2}}\sin(\sqrt{2}t)\Big(\dfrac{v_{1}-v_{2}}{18}e^{-2}+\dfrac{2(v_{1}-v_{2})}{9}e\Big),&\vspace*{0.5cm}\\
\bar{x}_{2}(t)=1-\dfrac{v_{1}+v_{2}}{2}t+\dfrac{v_{1}-v_{2}}{18}e^{2t-2}-\dfrac{v_{1}-v_{2}}{9}e^{1-t}-\dfrac{1}{3}\cos(\sqrt{2}t)
\Big(\dfrac{v_{1}-v_{2}}{6}e^{-2}-\dfrac{v_{1}-v_{2}}{3}e\Big)\vspace*{0.3cm}\\+\dfrac{1}{\sqrt{2}}\sin(\sqrt{2}t)\Big(\dfrac{v_{1}-v_{2}}{18}
e^{-2}+\dfrac{2(v_{1}-v_{2})}{9}e\Big),&\vspace*{0.5cm}\\
\bar{x}_{3}(t)=\bar{a}_{2}(t)+\bar{x}_{2}(t).
\end{cases}
\end{equation}
By the second condition in \eqref{hyp} we have the following three possibilities :
\begin{enumerate}
\item[\bf(i)] If $\bar{x}_{1}(1)=0$, then $ v_{1}=v_{2}-\dfrac{v_{2}-1}{c}$, where
$c:=\dfrac{4}{9}-\dfrac{\cos(\sqrt{2})}{9}\Big(\dfrac{e^{-2}}{2}-e\Big)+\dfrac{\sin(\sqrt{2})}{9\sqrt{2}}\Big(\dfrac{e^{-2}}{2}+2e\Big)$. It is easy to deduce from \eqref{e:8.9} and \eqref{e:8.11} that
\begin{equation*}
\begin{cases} \bar{a}_{1}(t)=-v_{2}+  \dfrac{v_{2}-1}{3c}e^{-(t-1)}+\dfrac{v_{2}-1}{6c}e^{2t-2} +\dfrac{v_{2}-1}{2c},&\vspace*{0.5cm}\\
\bar{a}_{2}(t)=-v_{2}- \dfrac{v_{2}-1}{3c}e^{-(t-1)}-\dfrac{v_{2}-1}{6c}e^{2t-2} +\dfrac{v_{2}-1}{2c}.
\end{cases}
\end{equation*}
\begin{equation*}
\begin{cases}\bar{x}_{1}(t)=1-v_{2}t+\dfrac{v_{2}-1}{c}\Big[ \dfrac{t}{2}+\dfrac{e^{2t-2}}{18}-\dfrac{e^{1-t}}{9}-\dfrac{1}{9}\cos(\sqrt{2}t)\Big(\dfrac{e^{-2}}{2}-e\Big)+\dfrac{1}{9\sqrt{2}}\sin(\sqrt{2}t)
\Big(\dfrac{e^{-2}}{2}+2e\Big)\Big] ,&\vspace*{0.5cm}\\
\bar{x}_{2}(t)=1-v_{2}t+\dfrac{v_{2}-1}{c}\Big[\dfrac{t}{2}-\dfrac{e^{2t-2}}{18}+\dfrac{e^{1-t}}{9}+\dfrac{1}{9}\cos(\sqrt{2}t)
\Big(\dfrac{e^{-2}}{2}-e\Big)-\dfrac{1}{9\sqrt{2}}\sin(\sqrt{2}t)\Big(\dfrac{e^{-2}}{2}+2e\Big)\Big],&\vspace*{0.5cm} \\
\bar{x}_{3}(t)=\bar{a}_{2}(t)+\bar{x}_{2}(t).
\end{cases}
\end{equation*}
Observe further that $\bar{x}_{2}(t)+\bar{x}_{1}(t)=2-2v_{2}t+\dfrac{v_{2}-1}{c}t$. Then we have that $\bar{x}_{2}(1)=2-2v_{2}+\dfrac{v_{2}-1}{c}$,
and that the cost functional in \eqref{eq2.6bis} reduces to
\begin{equation*}
\begin{aligned}
J[\bar{x},\bar{a}]&=\dfrac{1}{2}\Big(2-2v_{2}+\dfrac{v_{2}-1}{c}\Big)^{2}+\dfrac{1}{2}\disp\int\limits_{0}^{1}\Big[-v_{2}+\dfrac{v_{2}-1}{c}
\Big(\dfrac{1}{2}+\dfrac{e^{-(t-1)}}{3}+\dfrac{e^{2t-2}}{6}\Big)\Big]^{2}\,dt\\
&+\dfrac{1}{2}\disp\int\limits_{0}^{1}\Big[-v_{2}+\dfrac{v_{2}-1}{c}\Big(\dfrac{1}{2}-\dfrac{e^{-(t-1)}}{3}-\dfrac{e^{2t-2}}{6}\Big)\Big]^{2}\,dt.
\end{aligned}
\end{equation*}
This function clearly achieves its absolute minimum at $ v_{2}= 0.5988481275 $. Hence
\begin{equation*}
\begin{cases} \bar{a}_{1}(t)=-0.5988481275-  \dfrac{0.4011518725}{3c}e^{-(t-1)}
-\dfrac{0.4011518725}{6c}e^{2t-2} -\dfrac{0.4011518725}{2c},&\vspace*{0.5cm}\\
\bar{a}_{2}(t)=-0.5988481275+ \dfrac{0.4011518725}{3c}e^{-(t-1)}
+\dfrac{0.4011518725}{6c}^{2t-2}-\dfrac{0.4011518725}{2c},
\end{cases}
\end{equation*}
\begin{equation*}
\begin{cases}\bar{x}_{1}(t)=1-0.5988481275.t-\dfrac{0.4011518725}{c}\Big[ \dfrac{t}{2}+\dfrac{e^{2t-2}}{18}-\dfrac{e^{1-t}}{9}-\dfrac{1}{9}\cos(\sqrt{2}t)\Big(\dfrac{e^{-2}}{2}-e\Big)\vspace*{0.3cm}\\
+\dfrac{1}{9\sqrt{2}}\sin(\sqrt{2}t)\Big(\dfrac{e^{-2}}{2}+2e\Big)\Big],&\vspace*{0.5cm}\\
\bar{x}_{2}(t)=1-0.5988481275.t-\dfrac{0.4011518725}{c}\Big[ \dfrac{t}{2}-\dfrac{e^{2t-2}}{18}+\dfrac{e^{1-t}}{9}+\dfrac{1}{9}\cos(\sqrt{2}t)\Big(\dfrac{e^{-2}}{2}-e\Big)\vspace*{0.3cm}\\
-\dfrac{1}{9\sqrt{2}}\sin(\sqrt{2}t)\Big(\dfrac{e^{-2}}{2}+2e\Big)\Big], &\vspace*{0.5cm}\\
\bar{x}_{3}(t)=\bar{a}_{2}(t)+\bar{x}_{2}(t)\;\mbox{ for all }\;t\in[0,1].
\end{cases}
\end{equation*}
\item[\bf(ii)] If $\bar{x}_{2}(1)=0$, then $ v_{1}=v_{2}+\dfrac{v_{2}-1}{c}$, where $c:=-\dfrac{5}{9}-\dfrac{\cos(\sqrt{2})}{9}\Big(\dfrac{e^{-2}}{2}-e\Big)+\dfrac{\sin(\sqrt{2})}{9\sqrt{2}}\Big(\dfrac{e^{-2}}{2}+2e\Big)$.\\
Arguing similar to the above Case (i), we get that
\begin{equation*}
\begin{cases} \bar{a}_{1}(t)=-1.056787399- \dfrac{0.056787399}{3c}e^{-(t-1)}
-\dfrac{0.056787399}{6c}e^{2t-2} -\dfrac{0.056787399}{2c},&\vspace*{0.5cm}\\
\bar{a}_{2}(t)=-1.056787399+ \dfrac{0.056787399}{3c}e^{-(t-1)}
+\dfrac{0.056787399}{6c}e^{2t-2} -\dfrac{0.056787399}{2c},
\end{cases}
\end{equation*}
\begin{equation*}
\begin{cases}\bar{x}_{1}(t)=1-1.056787399.t+\dfrac{0.056787399}{c}\Big[ -\dfrac{t}{2}-\dfrac{e^{2t-2}}{18}+\dfrac{e^{1-t}}{9}+\dfrac{1}{9}\cos(\sqrt{2}t)\Big(\dfrac{e^{-2}}{2}-e\Big)\vspace*{0.3cm}\\
-\dfrac{1}{9\sqrt{2}}\sin(\sqrt{2}t)\Big(\dfrac{e^{-2}}{2}+2e\Big)\Big],&\vspace*{0.5cm}\\
\bar{x}_{2}(t)=1-1.056787399.t+\dfrac{0.056787399}{c}\Big[ -\dfrac{t}{2}+\dfrac{e^{2t-2}}{18}-\dfrac{e^{1-t}}{9}-\dfrac{1}{9}\cos(\sqrt{2}t)\Big(\dfrac{e^{-2}}{2}-e\Big)\vspace*{0.3cm}\\
+\dfrac{1}{9\sqrt{2}}\sin(\sqrt{2}t)\Big(\dfrac{e^{-2}}{2}+2e\Big)\Big],&\vspace*{0.5cm}\\
\bar{x}_{3}(t)=\bar{a}_{2}(t)+\bar{x}_{2}(t)\;\mbox{ for all }\;t\in[0,1].
\end{cases}
\end{equation*}
\item[\bf(iii)] Let $ \bar{x}_{1}(1)=0 $ and $ \bar{x}_{2}(1)=0 $. Then it follows from Case (i) that $ v_{1}=v_{2}-\dfrac{v_{2}-1}{c}$ and $0=2-2v_{2}+\dfrac{v_{2}-1}{c}$, where $c:=\dfrac{4}{9}-\dfrac{\cos(\sqrt{2})}{9}\Big(\dfrac{e^{-2}}{2}-e\Big)+\dfrac{\sin(\sqrt{2})}{9\sqrt{2}}\Big(\dfrac{e^{-2}}{2}+2e\Big)$.
Hence $ v_{1}=v_{2}=1 $.\\
Furthermore, we deduce from \eqref{e:8.9} and \eqref{e:8.11} the following expressions valid for all $t\in[0,1]$:
\begin{equation*}
\begin{cases}
\bar{a}_{1}(t)=\bar{a}_{2}(t)=-1,&\vspace*{0.5cm}\\
\bar{x}_{1}(t)=\bar{x}_{2}(t)=1-t,&\vspace*{0.5cm}\\
\bar{x}_{3}(t)=\bar{a}_{2}(t)+\bar{x}_{2}(t)=-t.
\end{cases}
\end{equation*}
\end{enumerate}
Thus the obtained necessary optimality conditions from Theorem~\ref{Thm:N.C.I} with $\lm\ne 0$ allows us to fully compute the unique optimal solution of the integro-differential sweeping control problem formulated in \eqref{eq2.6bis} and \eqref{eq2.7} with the initial data given in \eqref{pv} such that condition \eqref{hyp} is satisfied.
\end{example}
\end{document}